\documentclass[11pt,a4paper,final]{article}
\usepackage[left=2.00cm, right=2.00cm, top=1.00cm, bottom=2.00cm]{geometry}
\usepackage[]{inputenc}
\usepackage{amsmath}
\usepackage{amsfonts}
\usepackage{amssymb}
\usepackage{graphicx}

\usepackage{moreverb}
\usepackage{changepage}
\usepackage{cite}
\usepackage{cancel}
\usepackage{graphicx}
\usepackage{caption}
\usepackage{subcaption}
\usepackage{slashed}
\usepackage{accents}							
\newlength{\dhatheight}
\newcommand{\doublehat}[1]{%
    \settoheight{\dhatheight}{\ensuremath{\hat{#1}}}%
    \addtolength{\dhatheight}{-0.25ex}%
    \hat{\vphantom{\rule{5pt}{\dhatheight}}%
    \smash{\hat{#1}}}}

\DeclareMathAlphabet{\mathpzc}{OT1}{pzc}{m}{it} 

\newcommand\BibTeX{{\rmfamily B\kern-.05em \textsc{i\kern-.025em b}\kern-.08em
T\kern-.1667em\lower.7ex\hbox{E}\kern-.125emX}}

\begin{document}
\title{An Efficient Semi-Implicit Method for Three-Dimensional Non-Hydrostatic Flows in Compliant Arterial Vessels}

\author{Francesco Fambri\\
Department of Physics, University of Trento, I-38123 Trento, Italy \\
francesco.fambri@unitn.it\\
\and
Michael Dumbser \;\;and\;\; Vincenzo Casulli\\
Dep. of Civil, Env. and Mech. Engineering, University of Trento, I-38123 Trento, Italy } 

\date{May 19,  2014} 
\maketitle

\begin{abstract}
Blood flow in arterial systems can be described by the three-dimensional Navier-Stokes equations within a time-dependent spatial domain that accounts for the elasticity of the arterial walls. In this article blood is treated  as an incompressible  Newtonian fluid that flows through compliant vessels of general cross section. A three-dimensional semi-implicit finite difference and finite volume model is derived so that numerical stability is obtained at a low computational cost on a staggered grid.  The key idea of the method consists in a splitting of the pressure into a hydrostatic and a non-hydrostatic part, where first a small quasi-one-dimensional nonlinear system is solved for the hydrostatic pressure and only in a second step the fully three-dimensional non-hydrostatic pressure is computed from a three-dimensional nonlinear system as a correction to the hydrostatic one. The resulting algorithm is robust, efficient, locally and globally mass conservative and applies to hydrostatic and non hydrostatic flows in one, two and  three space dimensions. These features are illustrated on nontrivial test cases for flows in tubes with circular or elliptical cross section where the exact analytical solution is known. Test cases of steady and pulsatile flows in uniformly curved rigid and elastic tubes are presented.  Wherever possible, axial velocity development and secondary flows are shown and compared with previously published results. 
\end{abstract}

\textbf{keywords}-- blood flow; compliant arteries; moving boundaries; non-hydrostatic; semi-implicit; three-dimensional; curved blood vessels; general  cross-sections;  secondary flows; finite difference; finite volume.

\section{Introduction}
Blood flow in medium to large arterial systems can be accurately described by the three-dimensional Navier-Stokes equations within a time-dependent spatial domain. These equations, however, are too complex to be efficiently solved over a large network of systemic arteries, and alternatively, simplified one-dimensional equations are often used \cite{Sherwin:2003,Sherwin:2003b,MuellerJCP,MuellerCNM}. 
Two dimensional models, like the ones proposed in \cite{Canic:2005} and \cite{Casulli:2012}, have the aim of filling the gap between highly complex and computationally  expensive three-dimensional fluid-structure interaction models  \cite{Formaggia:2001,Formaggia:2007,Formaggia:2009,Lohner2002,Lohner2003,Lohner2010} and reduced one-dimensional equations, where the missing information about the velocity profiles has to be derived from empirical closure relations.  
Recently a very efficient semi-implicit method for two-dimensional axially symmetric blood flow in compliant arterial systems has been introduced in  \cite{Casulli:2012}, which has been subsequently extended also to the simulation of substance transport in \cite{Tavelli:2013}. The present three-dimensional model aims at extending the above-mentioned two-dimensional model to the third spatial dimension in order to give accurate and 
detailed three-dimensional information about the pressure and the velocity field. Non-hydrostatic corrections of the pressure along the cross section  are necessary to give an accurate description of secondary flows, and consequently to simulate blood flow through curved sections in the vascular tree. This is very important because asymmetric irregularities of blood flow due to curved sections, junctions \cite{Katristis:2007}, stenosis \cite{Ikbal:2012}, and turbulence instabilities \cite{Prado:2008} can cause the leading contributions to the wall shear stress.

In general fully non-hydrostatic problems show no favorite direction. This implies that a fine spatial discretization is needed along any spatial coordinate. However, certain non-hydrostatic flow problems may still contain a preferential direction, such as the gravity direction for 
gravity-driven free surface flows or the axial direction for blood flow in compliant arteries. In the 1990ies, very efficient semi-implicit 
fractional step methods have been developed for gravity-driven three-dimensional non hydrostatic free-surface flows. 
In \cite{Casulli:1999, Casulli:2002}, the hydrostatic component of the pressure is 
computed in a first fractional step and in a second fractional step a non-hydrostatic correction is determined in order to obtain the fully 
three-dimensional non-hydrostatic pressure field. The method is relatively simple, numerically stable for large Courant numbers and particularly 
suitable for three-dimensional problems in which the non-hydrostatic correction is small compared to the hydrostatic components. 
For such methods, the computational efficiency is the larger the more the problem is hydrostatic. 

Encouraged by the success of those methods, a new fractional method is designed in this paper in order to face the problem of three-dimensional  non-hydrostatic blood flow in compliant vessels with rather general cross section. At first the model is developed in cylindrical coordinates, 
then curvature is introduced. The resulting approach is fast, simple and provides full details of the radial and angular velocity profiles and 
of the non-hydrostatic pressure field. It will be shown by computational results that the new numerical method is also able to model secondary flows correctly.  

\section{Governing equations}

\begin{figure}[btp] 
\centering 
			\begin{subfigure}[b] {0.35\textwidth}
			\centering
			\includegraphics[width=\textwidth]{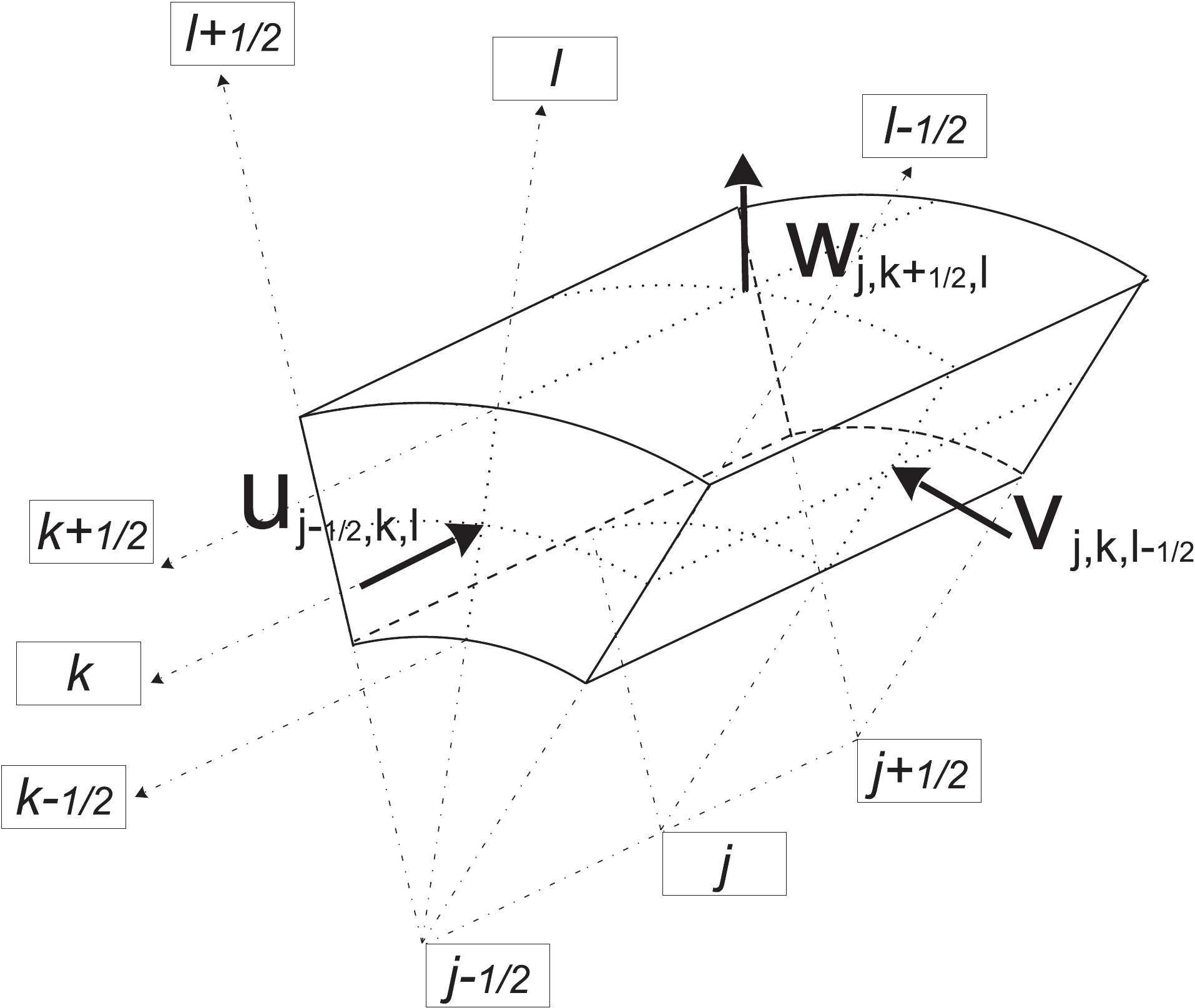}
			\subcaption{} 
			\label{fig:cap4VolElem}
			\end{subfigure}
			\hspace{5mm}
			\begin{subfigure}[b] {0.45\textwidth}
			\centering
			\includegraphics[width=\textwidth]{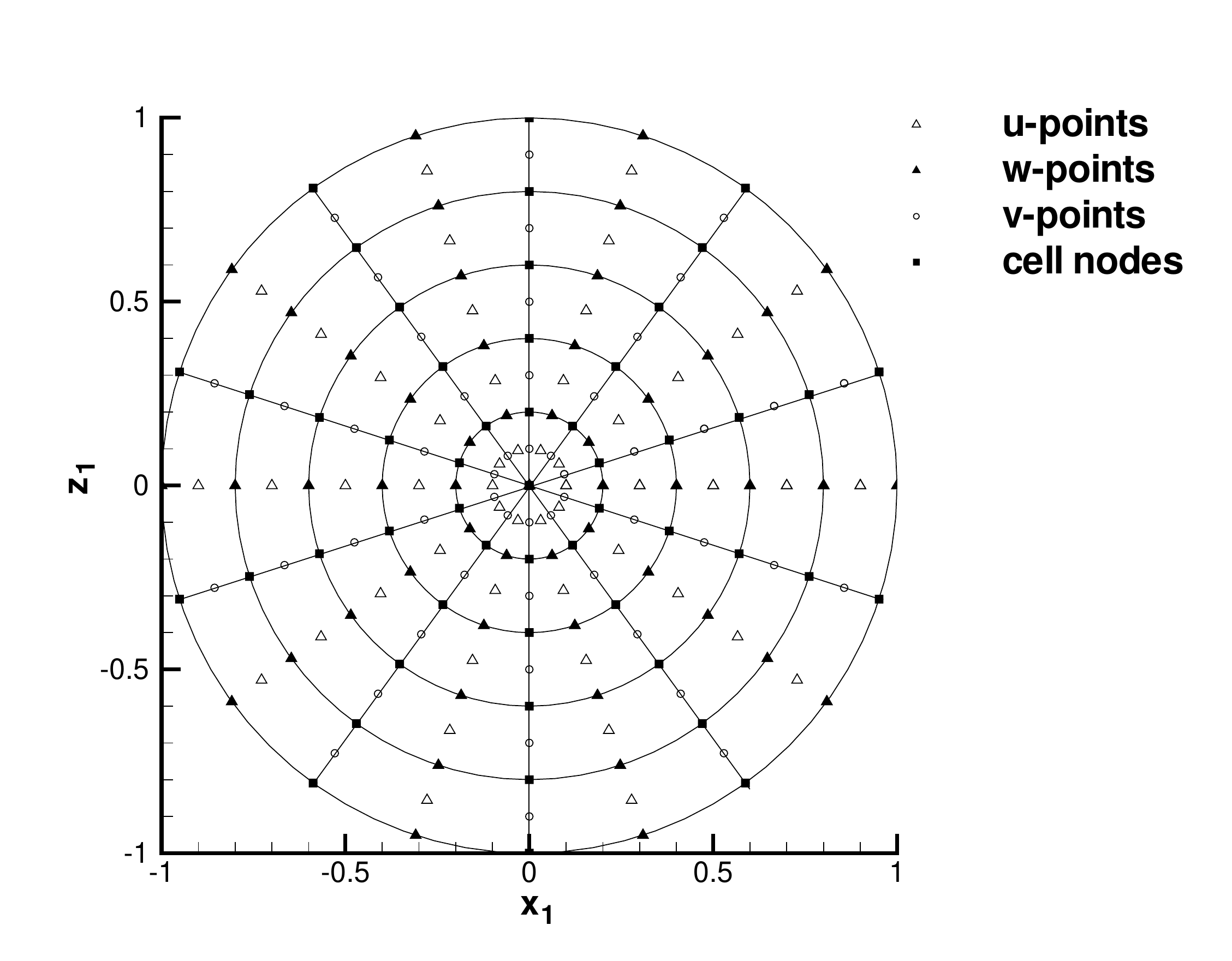}
			\subcaption{} 
			\label{fig:cap4section}
			\end{subfigure}
			\vspace{5mm}
			\begin{subfigure}[b]{0.35\textwidth}
			\centering
			\includegraphics[width=\textwidth]{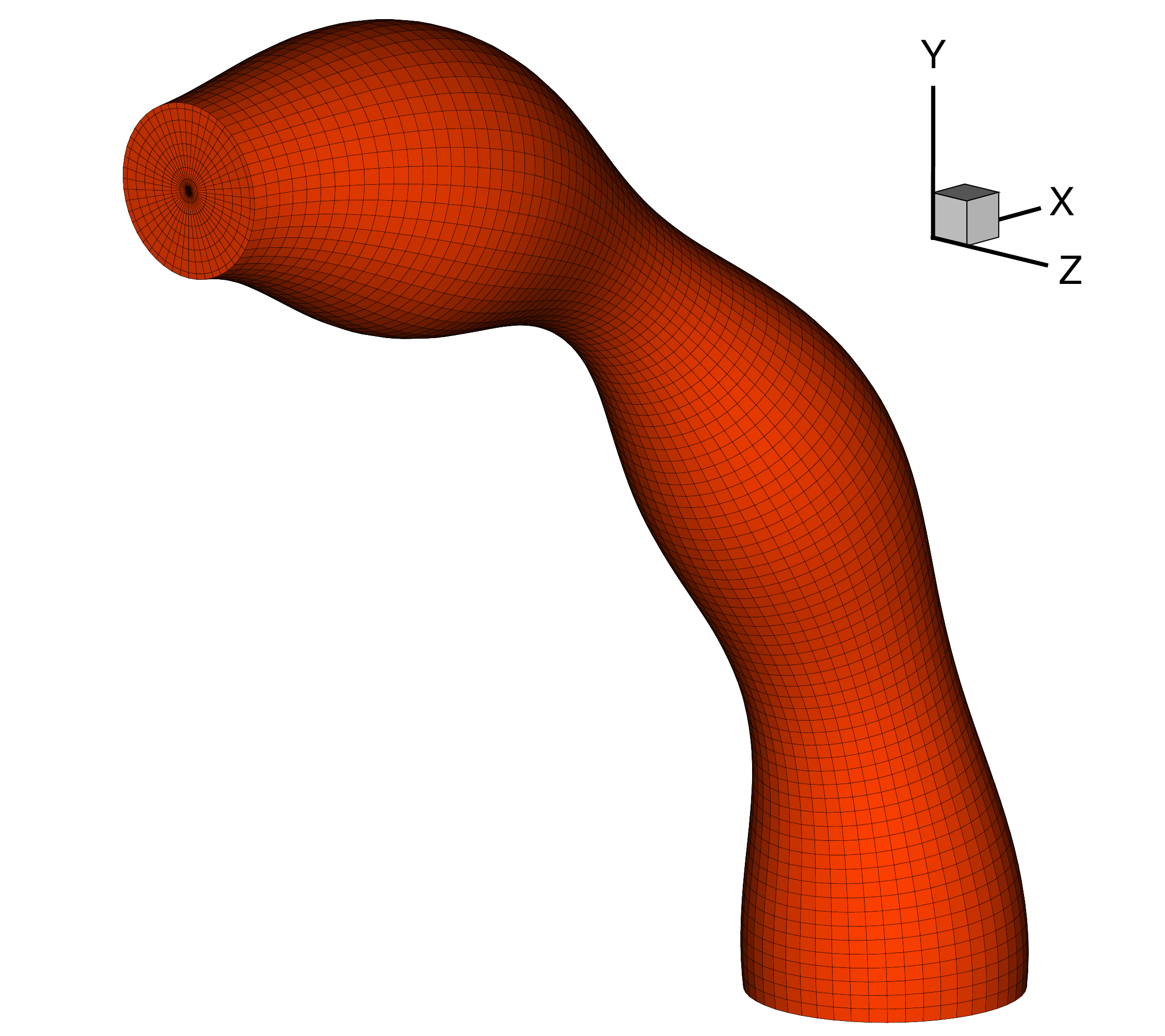}
			\subcaption{} 
			\label{fig:cap5FineMesh}
			\end{subfigure}  \hspace{5mm}
			\begin{subfigure}[b]{0.45\textwidth}
			\centering
			\includegraphics[width=\textwidth]{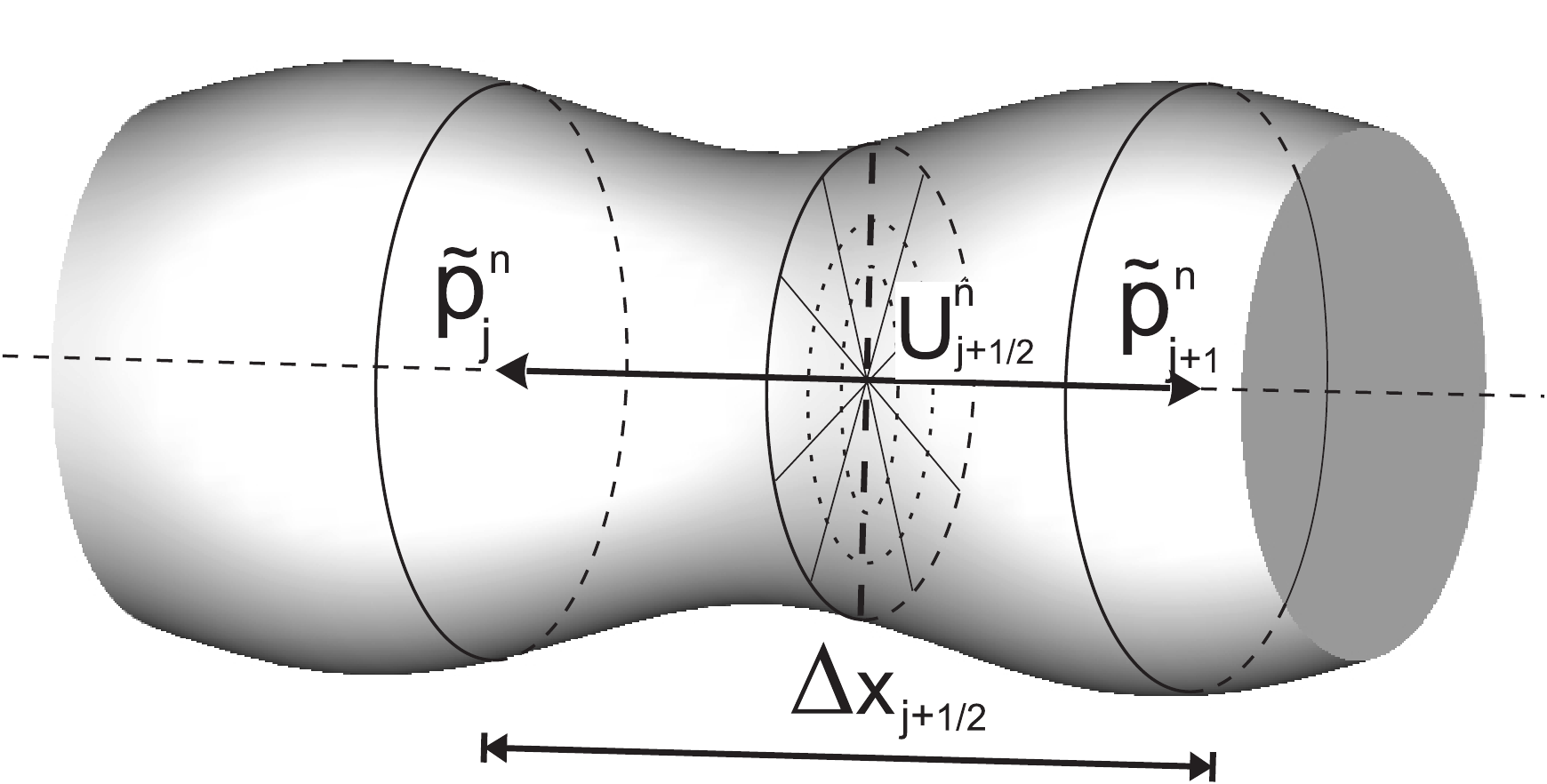}
			\subcaption{} 
			\label{fig:cap4StaggRep}
			\end{subfigure}
\caption{Geometry of the numerical grid: (\ref{fig:cap4VolElem}) Typical geometry of a control volume of the model; (\ref{fig:cap4section}) Axial projection of the nodes and velocity points of the $\mathit{j}$-th segment; (\ref{fig:cap5FineMesh}) Example of  a curved numerical grid; (\ref{fig:cap4StaggRep}) Staggered representation of the $(\mathit{j}+\frac{1}{2})$-th segment;   }\label{fig:cap4NumModel2}
\end{figure}

In this article the three-dimensional momentum and continuity equations in cylindrical coordinates are considered, 
\begin{align}
\frac{d u}{d t}   &= - \frac{\partial p}{\partial x}+ \nu \left[ \frac{\partial^2 u}{\partial x^2} + \frac{1}{z}\frac{\partial}{\partial z} \left(z \frac{\partial u}{\partial z} \right) + \frac{1}{z^2}\frac{\partial^2 u}{\partial \varphi^2} \right],   \label{eq:cap5NSu}\\
\frac{d v}{d t} &= - \frac{1}{ z} \frac{\partial p}{\partial \varphi}+  \nu \left[ \frac{\partial^2 v}{\partial x^2}+ \frac{\partial}{\partial z} \left( \frac{1}{z}\frac{\partial}{\partial z}  \left(z v \right) \right) + \frac{1}{z^2}\frac{\partial^2 v}{\partial \varphi^2} + \frac{2}{z^2} \frac{\partial w}{\partial \varphi} \right], \label{eq:cap5NSv}\\
\frac{d w}{d t} &=  -  \frac{\partial p}{\partial z}+\nu \left[ \frac{\partial^2 w}{\partial x^2} + \frac{\partial}{\partial z} \left( \frac{1}{z}\frac{\partial}{\partial z}  \left(z w \right) \right) + \frac{1}{z^2}\frac{\partial^2 w}{\partial \varphi^2} - \frac{2}{z^2} \frac{\partial v}{\partial \varphi} \right], \label{eq:cap5NSw}\\
0&=\frac{\partial}{\partial x}\left( z u\right) +  \frac{\partial}{\partial z} \left( z w\right) +  \frac{\partial v}{\partial \varphi}.
\label{eqn.conti} 
\end{align}
where $(x,z,\varphi)$ are the axial position along the vessel, the radial coordinate and the radial angle, respectively;   $u(x, z, \varphi , t)$, $v(x, z, \varphi , t)$ and $w(x, z, \varphi , t)$ are the velocity components in the  $x$-, $\varphi$- and $z$- directions, respectively; 
$p(x, z, \varphi , t)$ is the normalized pressure field defined as the pressure divided by the density $\rho$ of the fluid, that is assumed to be constant; $\nu$ is the kinematic viscosity coefficient; the left hand sides are the material derivatives of the respective velocity components. %

In this paper the (normalized) pressure $\mathit{p}$ is written as a sum of the \emph{hydrostatic} component $\tilde{\mathit{p}}$, that 
is only a function of the axial coordinate and time, and the fully three-dimensional \emph{non-hydrostatic} component $\mathit{q}$, having
\begin{equation}
p \left( x, z, \varphi , t\right) = \tilde{p}\left(x, t \right) + q \left(x,z, \varphi, t \right).
\label{eq:cap5pressure}
\end{equation}
Note that for curved vessels the pressure field can actually be written as the sum of a hydrostatic pressure and some perturbative corrections, 
which can be expressed as a power series in the curvature ratio, see \cite{Berger:1983}. Hence, the choice of splitting the total pressure 
into a hydrostatic one and a non-hydrostatic correction as done in Eq. \eqref{eq:cap5pressure} is justified. For flow problems in curved vessels 
the algorithm designed in the following will be therefore particularly efficient when corrections caused by curvature are actually  
perturbative. 

In order to close the problem, an \emph{equation of state} relating the arterial radius $R(x,\varphi,t)$ to the unknown pressure at the wall boundaries $p(x,R,\varphi,t)$ is needed. To this purpose, a typical choice is the law of Laplace, see \cite{Canic:2006} for example,  
\begin{gather}
p(x,R,\varphi,t) = p_{\text{ext}}(x,\varphi,t) + \beta(x,\varphi,t) \left( R(x,\varphi,t) - R_0(x,\varphi)\right),\label{eq:cap5laplace}
\end{gather}
where $p_{\text{ext}}$ is a prescribed external pressure function, $\beta$ is a positive and in general variable \emph{rigidity coefficient} that can be obtained 
for example by the solution of inverse problems, see for example \cite{Perego,Wall2013}, and $R_0(x,\varphi)$ is the equilibrium radius, the dependence of which on 
the axial and angular coordinate allows to represent rather general cross-sections. According to \cite{Sherwin:2003}, the coefficient $\beta$ is related to basic material 
properties by the formula 
\begin{equation}
	\beta = \frac{h_0 E}{(1-\nu^2) R_0^2},  
	\label{eqn.beta.def} 
\end{equation}
where $h_0$ is the wall thickness, $E$ is the Young modulus and $\nu$ is the Poisson ratio of the vessel wall. 
The only restriction on the cross section in our present model is that the equilibrium 
radius must be a \textit{single-valued function} of $x$ and $\varphi$. Often, the geometry of large to medium-scale blood vessels is obtained from in vivo CT images, but 
usually the vessel geometry obtained in that manner is \textit{not} the equilibrium geometry, but already a \textit{deformed} configuration. A special technique to obtain 
the corresponding stresses, a so-called \textit{prestressing} method for patient-specific geometries and material parameters can be found in \cite{Wall2010b}. 

Integrating the continuity equation over the cross sectional area and using a kinematic condition at the moving boundary leads to the following equation for the \emph{moving vessel boundary}
\begin{gather}
\frac{ \partial A}{\partial t} + \frac{\partial}{\partial x} \left(\int_0^{2\pi} \int_0^R{z u \, dz\,d \varphi} \right) = 0,\label{eq:walls}
\end{gather}
where $A$ is the area of the cross section.

The no-slip boundary conditions and the kinematic radial condition at the vessel wall can be summarized as 
\begin{align}
u\left(x,R,\varphi,t \right) &= v\left(x,R,\varphi,t \right) = 0;\label{eq:cap5BCuv}\\
w\left(x,R,\varphi,t \right) &= \frac{\partial R }{\partial t}.\label{eq:cap5BCw}
\end{align} 

\section{A semi-implicit finite difference-finite volume model}

In non-hydrostatic flow models for incompressible fluids the pressure terms must always be discretized implicitly, due to the elliptic nature of the PDE that governs the non-hydrostatic pressure. 
However, if a fully implicit method was used, the resulting numerical scheme would be far too complex and difficult to control, mainly because of its nonlinearities. 

The \emph{criteria} for choosing the precise terms that need to be discretized implicitly or explicitly are clearly indicated in \cite{Casulli:1990} for the two-dimensional shallow water equations, with the aim of improving  efficiency, robustness and to ensure stability. It is shown how the characteristic analysis of the governing equations leads to a particular semi-implicit discretization. The resulting algorithm was soon extended to three-dimensional problems \cite{Casulli:1992, Casulli:1994}, non-hydrostatic flows \cite{Casulli:1999} and it was further adapted to orthogonal unstructured grids, which admit higher flexibility in fitting complex geometries \cite{Casulli:2000,BoscheriDumbser} and are suitable for general irregular domains of arbitrary scale \cite{Casulli:2002}.

Within the aforementioned framework, it has been shown that an implicit discretization of the pressure gradients in the momentum equations (\ref{eq:cap5NSu})-(\ref{eq:cap5NSw}) and of the velocity in the moving boundary equation (\ref{eq:walls}) can avoid the CFL restriction on the gravity wave speed. Stability has been further ensured by an implicit discretization of the viscosity and the wall friction \cite{Casulli:1994}. Further to that, the discretization in the moving boundary equation (\ref{eq:walls}) gives rise to a particular nonlinearity. This kind of nonlinearity can be efficiently solved by an iterative Newton-type algorithm as pointed out in later works, 
see \cite{Casulli:2009} and \cite{Casulli:2011b}. A detailed convergence analysis of this solution algorithm for the resulting mildly nonlinear systems is illustrated in \cite{Brugnano:2009,BrugnanoCasulli,BrugnanoSestini} and \cite{Casulli:2012b}. The resulting algorithm is  very efficient, highly accurate, and stable for large time-steps. 

\subsection{Unstructured staggered grid}

To simulate arterial flows, one assumes that the arterial system consists in a set of interconnected arterial branches where the flow is governed by equations (\ref{eq:cap5NSu})-(\ref{eq:cap5BCw}). For a single branch of reference the axial domain is discretized  in $N_x$ non overlapping segments, the radial domain up to $N_z$ rings and the angular domain in $N_{\varphi}$ slices  (see Figure \ref{fig:cap5FineMesh}). The resulting control volumes are at most $N_x  N_z  N_{\varphi}$ three-dimensional annular sectors. Velocity components are defined on the barycentres of the faces of such control volumes, in such a way that the fluxes in the local mass conservation equations are well defined. Figures \ref{fig:cap4VolElem} and \ref{fig:cap4section} summarize efficiently the chosen  convention of the indexing on the chosen staggered grid. 

The discrete arterial radius at the $i$-th axial location, $l$-th angular slice and time level $t_n$ is denoted with $R_{i,l}^n$. 
Specific details about the chosen radial mesh are available in \cite{Casulli:2012}.

When axial intervals are labeled with only the axial index, the intervals are referred to the central position, along the axis of the vessel. For simplicity, the angular mesh is chosen to be uniform, having $\Delta \varphi = 2\pi/N_\varphi$. 

\subsection{Semi-implicit discretization}
A consistent semi-implicit finite difference discretization of the momentum equations for each control volume $(j,k,l)$ reads 
\begin{align}
\frac{u_{j+\frac{1}{2},k,l}^{n+1}-u_{j+\frac{1}{2},k,l}^{n,L}}{\Delta t}  &= - \frac{ p^{n+\theta}_{j+1,k,l} - p^{n+\theta}_{j,k,l}}
{\Delta x_{j+\frac{1}{2},k,l}} 
+ \nu \left( \mathcal{L}^{ux}_{j+\frac{1}{2},k,l}[u^{n,L}] + \mathcal{L}^{uz\varphi}_{j+\frac{1}{2},k,l}[u^{n+\theta}] \right), \label{eq:cap5NumAxial}
\end{align}
\begin{align}
\frac{v_{j,k,l+\frac{1}{2}}^{n+1}-v_{j,k,l+\frac{1}{2}}^{n,L}}{\Delta t}  &=  - \frac{ p^{n+\theta}_{j,k,l+1} - p^{n+\theta}_{j,k,l}}
{z_{j,k,l+\frac{1}{2}}^n \Delta \varphi} 
+ \nu \left( \mathcal{L}^{vx\varphi}_{j,k,l+\frac{1}{2}}[v^{n,L},w^{n,L}] + \mathcal{L}^{vz\varphi}_{j,k,l+\frac{1}{2}}[v^{n+\theta}] \right),\label{eq:cap5NumTang}\\
\frac{w_{j,k+\frac{1}{2},l}^{n+1}-w_{j,k+\frac{1}{2},l}^{n,L}}{\Delta t}  &= - \frac{ p^{n+\theta}_{j,k+1,l} - p^{n+\theta}_{j,k,l}}
{\Delta z_{j,k+\frac{1}{2},l}^n} 
+ \nu \left( \mathcal{L}^{wx\varphi}_{j,k+\frac{1}{2},l}[w^{n,L},v^{n,L}] + \mathcal{L}^{wz\varphi}_{j,k+\frac{1}{2},l}[w^{n+\theta}] \right),\label{eq:cap5NumRad}
\end{align}
where $\Delta t$ is the time-step size; $\theta$ is an implicitness factor to be taken in the range $\frac{1}{2} \leq \theta \leq 1 $  (
\cite{Casulli:1994} 
); $p^{n+\theta} = \theta p^{n+1} + (1-\theta) p^n$;  $u_{j+1/2,k,l}^{n,L}$, $v_{j,k,l+1/2}^{n,L}$ and $w_{j,k+1/2,l}^{n,L}$ denote respectively the axial, tangential and radial velocity components interpolated at time $t_n$ at the end of the Lagrangian trajectory in the three-dimensional space. 
In particular, one has 
\begin{align*}
\mathcal{L}^{ux}_{j+\frac{1}{2},k,l}[u^{n,L}]  &=  \frac{\frac{u^{n,L}_{j+\frac{3}{2},k,l} - u^{n,L}_{j+\frac{1}{2},k,l}}{\Delta x_{j+1,k,l}} - \frac{u^{n,L}_{j+\frac{1}{2},k,l} -u^{n,L}_{j-\frac{1}{2},k,l-1}}{\Delta x_{j,k,l}}}{\Delta x_{j+\frac{1}{2},k,l}},  
\end{align*}
\begin{align}
\mathcal{L}^{uz\varphi}_{j+\frac{1}{2},k,l}[u^{n+\theta}]  &=  
  \frac{1}{z^2_{j + \frac{1}{2},k}}\frac{u^{n+\theta}_{j + \frac{1}{2},k,l+1} - 2u^{n+\theta}_{j + \frac{1}{2},k,l} +u^{n+\theta}_{j + \frac{1}{2},k,l-1}}{\Delta \varphi^2} + \nonumber \\ 
  &+ \frac{z^n_{j + \frac{1}{2},k+\frac{1}{2}} \frac{u^{n+\theta}_{j + \frac{1}{2},k+1,l} - u^{n+\theta}_{j + \frac{1}{2},k,l}}{\Delta z^n_{j + \frac{1}{2},k+\frac{1}{2}}}- z^n_{j + \frac{1}{2},k-\frac{1}{2}} \frac{u^{n+\theta}_{j + \frac{1}{2},k,l} - u^{n+\theta}_{j + \frac{1}{2},k-1,l}}{\Delta z^n_{k-\frac{1}{2}}}}{z^n_{j + \frac{1}{2},k} \Delta z^n_{j + \frac{1}{2},k}} ,\label{eq:cap5ViscOperU}
\end{align}
\begin{align*}
\mathcal{L}^{vx\varphi}_{j,k,l+\frac{1}{2}}[v^{n,L},w^{n}] & =    \frac{\frac{v^{n,L}_{j+1,k,l + \frac{1}{2}} - v^{n,L}_{j,k,l + \frac{1}{2}}}{\Delta x_{j,k,l+1}} - \frac{v^{n,L}_{j,k,l + \frac{1}{2}}- v^{n,L}_{j-1,k,l + \frac{1}{2}}}{\Delta x_{j,k,l}}}{\Delta x_{j,k,l+\frac{1}{2}}} + \nonumber\\
& + \frac{w_{j,k+\frac{1}{2},l+1}+w_{j,k-\frac{1}{2},l+1}-w_{j,k+\frac{1}{2},l}-w_{j,k-\frac{1}{2},l}}{z^2_{j,k}\Delta \varphi}, \nonumber 
\end{align*}
\begin{align} 
\mathcal{L}^{vz\varphi}_{j,k,l+\frac{1}{2}}[v^{n+\theta}] & =   \frac{ \frac{z^n_{j,k+1}v^{n+\theta}_{j,k+1,l + \frac{1}{2}} - z^n_{j,k}v^{n+\theta}_{j,k,l + \frac{1}{2}}}{z^n_{j,k+\frac{1}{2}}\Delta z^n_{j,k+\frac{1}{2}}}- \frac{ z^n_{j,k}v^{n+\theta}_{j,k,l + \frac{1}{2}} - z^n_{j,k-1}v^{n+\theta}_{j,k-1,l + \frac{1}{2}}}{z^n_{j,k-\frac{1}{2}} \Delta z^n_{k-\frac{1}{2}}}}{\Delta z^n_{j,k}} + \nonumber\\
&  +  \frac{1}{z^2_{j,k}}\frac{v^{n+\theta}_{j,k,l+\frac{3}{2}} - 2v^{n+\theta}_{j,k,l+\frac{1}{2}} +v^{n+\theta}_{j,k,l-\frac{1}{2}}}{\Delta \varphi^2} ,\label{eq:cap5ViscOperV}
\end{align} 
\begin{align*}
\mathcal{L}^{wx\varphi}_{j,k+\frac{1}{2},l}[w^{n,L},v^{n}]  &=  \frac{\frac{w^{n,L}_{j+1,k+\frac{1}{2},l} - w^{n,L}_{j,k+\frac{1}{2},l}}{\Delta x_{j,k+1,l}} - \frac{  w^{n,L}_{j,k+\frac{1}{2},l}-w^{n,L}_{j-1,k+\frac{1}{2},l}}{\Delta x_{j,k}}}{\Delta x_{j,k+\frac{1}{2},l}} + \nonumber\\
& -  \frac{v_{j,k,l+\frac{1}{2}}+v_{j,k+1,l+\frac{1}{2}}-v_{j,k,l-\frac{1}{2}}-v_{j,k+1,l-\frac{1}{2}}}{z^2_{j,k+\frac{1}{2}}\Delta \varphi}
\end{align*} 
\begin{align}
\mathcal{L}^{wz\varphi}_{j,k+\frac{1}{2},l}[w^{n+\theta}]  &=  \frac{ \frac{z^n_{j,k+\frac{3}{2}}w^{n+\theta}_{j,k+\frac{3}{2},l} - z^n_{j,k+\frac{1}{2}}w^{n+\theta}_{j,k+\frac{1}{2},l}}{z^n_{j,k+1}\Delta z^n_{j,k+1}}- \frac{ z^n_{j,k+\frac{1}{2}}w^{n+\theta}_{j,k+\frac{1}{2},l} - z^n_{j,k-\frac{1}{2}}w^{n+\theta}_{j,k-\frac{1}{2},l}}{z^n_{j,k} \Delta z^n_{k}}}{\Delta z^n_{j,k+\frac{1}{2}}} + \nonumber\\
&  +  \frac{1}{z^2_{j,k+\frac{1}{2}}}\frac{w^{n+\theta}_{j,k+\frac{1}{2},l+1} - 2w^{n+\theta}_{j,k+\frac{1}{2},l} +w^{n+\theta}_{j,k+\frac{1}{2},l-1}}{\Delta \varphi^2} . \label{eq:cap5ViscOperW}
\end{align}

Alternative explicit schemes for discretizing the advective terms can be employed. For example, an accurate description of rapidly varying flows can be achieved by using  a conservative formulation \cite{Stelling:2003}.

In (\ref{eq:cap5NumAxial})-(\ref{eq:cap5NumRad}) the operators $\nu \mathcal{L}^{uz\varphi}_{j+\frac{1}{2},k,l}[u^{n+\theta}]$, $\nu \mathcal{L}^{vz\varphi}_{j,k,l+\frac{1}{2}}[v^{n+\theta}]$ and $\nu \mathcal{L}^{wz\varphi}_{j,k+\frac{1}{2},l}[w^{n+\theta}]$ denote the implicit discretization of the viscous force acting respectively on $u$, $v$ and $w$ in the radial direction, 
while $\nu \mathcal{L}^{ux}_{j+\frac{1}{2},k,l}[u^{n,L}]$, $\nu \mathcal{L}^{vx\varphi}_{j,k,l+\frac{1}{2}}[v^{n,L},w^{n,L}]$ and $\nu\mathcal{L}^{wx\varphi}_{j,k+\frac{1}{2},l}[w^{n,L},v^{n,L}]$ denote discretization of the remaining viscous terms and are taken explicitly at the footpoint of the Lagrangian trajectory.

Concerning the discrete boundary conditions, the radial component of the velocity $w$ at the wall is the time derivative of the local radius, see Eq. (\ref{eq:cap5BCw}), that is related to the local pressure through the equation of state (\ref{eq:cap5laplace}). Hence, one obtains the boundary conditions 
\begin{align}
u_{j+\frac{1}{2},K(j)+\frac{1}{2},l}^{n+1} &= v_{j,K(j)+\frac{1}{2},l+\frac{1}{2}}^{n+1} = 0, \label{eq:NumBCuv}\\
w^{n+1}_{j,K(j,l)+\frac{1}{2},l} &= \frac{1}{ \beta} \frac{p_{j,K,l}^{n+1}-p_{j,K,l}^n}{\Delta t}. \label{eq:cap5wR}
\end{align}
where $K$ is the local radial index such that $z_{i,K+1/2,l}^n = R_{i,l}^n$, with $z_{i,1/2,l}^n = 0$.
At each hydrostatic pressure point $j$, the semi-implicit \emph{finite volume} approximation of the moving boundary equation (\ref{eq:walls}) reads 
\begin{align}
V_j \left( p^{n+1} \right) &= V_j \left( p^n  \right) - \Delta t \sum_{l=1}^{N_{\varphi}} \left[ \sum_{k=1}^{K^n_{j+\frac{1}{2},l}}a^n_{j+\frac{1}{2},k,l} u^{n+\theta}_{j+\frac{1}{2},k,l}  - \sum_{k=1}^{K^n_{j-\frac{1}{2},l}}a^n_{j-\frac{1}{2},k,l}  u^{n+\theta}_{j-\frac{1}{2},k,l} \right],\label{eq:cap5NumVolume}\\
V_j\left( p  \right) &= \frac{\pi}{2N_{\varphi}} \sum_{l=1}^{N_{\varphi}}\left\{ \Delta x_{j+\frac{1}{2}} \left[ R_{j+\frac{1}{2},l} \left(p_{j,K,l}\right)\right]^2  +\Delta x_{j-\frac{1}{2}} \left[R_{j-\frac{1}{2},l} \left(p_{j,K,l}\right)\right]^2\right\},\label{eq:cap5volume}
\end{align}
where $V_j$ is the nonlinear volume function, $a^n_{j+\frac{1}{2},k,l} = z^n_{j+\frac{1}{2},k,l} \Delta z^n_{j+\frac{1}{2},k,l}\Delta \varphi$ is the area of the axial  surface, 
$R_{j+\frac{1}{2},l}(p_{j,K,l})$ is computed through the equation of state (\ref{eq:cap5laplace}). 

Let us write the velocity and pressure fields as the sum of the hydrostatic component  ($\tilde{\mathbf{v}}$ and $\tilde{p}$) and the non-hydrostatic correction ($\delta\mathbf{v}$ and $\mathit{q}$). At the first fractional step a 
\emph{hydrostatic problem} is solved, defined by equations (\ref{eq:cap5NumAxial})-(\ref{eq:cap5NumRad}) and (\ref{eq:cap5NumVolume}) by neglecting all the implicit non-hydrostatic terms and by using a purely explicit discretization of non-hydrostatic viscous terms, having
\begin{align}
\frac{\tilde{u}_{j+\frac{1}{2},k,l}^{n+1}-u_{j+\frac{1}{2},k,l}^{n,L}}{\Delta t}  &= - \frac{ \tilde{p}^{n+\theta}_{j+1}- \tilde{p}^{n+\theta}_{j}}{\Delta x_{j+\frac{1}{2},k,l}}  
- \left(1-\theta' \right)  \frac{ q^{n}_{j+1,k,l} - q^{n}_{j,k,l}}{\Delta x_{j+\frac{1}{2},k,l}}  + \nonumber \\ 
 &  + \nu \left( \mathcal{L}^{ux}_{j+\frac{1}{2},k,l}\left[ u^{n,L} \right] + \mathcal{L}^{uz\varphi}_{j+\frac{1}{2},k,l}\left[ \tilde{u}^{n+\theta}\right] +  \mathcal{L}^{uz\varphi}_{j+\frac{1}{2},k,l}\left[ \delta u^{n,L} \right]\right)  ,\quad \quad\label{eq:cap5NumAxialH}
\end{align}
\begin{align}
\frac{\tilde{v}_{j,k,l+\frac{1}{2}}^{n+1}-v_{j,k,l+\frac{1}{2}}^{n,L}}{\Delta t}  &=   - \left(1-\theta' \right)  \frac{ q^{n}_{j,k,l+1} - q^{n}_{j,k,l}}{z^n_{j,k,l+\frac{1}{2}}\Delta \varphi} + \nonumber \\ 
 &+ \nu \left( \mathcal{L}^{vx\varphi}_{j,k,l+\frac{1}{2}}\left[ v^{n,L},w^{n,L} \right] + \mathcal{L}^{vz\varphi}_{j,k,l+\frac{1}{2}}\left[ \tilde{v}^{n+\theta}  \right]  + \mathcal{L}^{vz\varphi}_{j,k,l+\frac{1}{2}}\left[ \delta v^{n,L} \right]\right),\label{eq:cap5NumTangH}\\
\frac{\tilde{w}_{j,k+\frac{1}{2},l}^{n+1}-w_{j,k+\frac{1}{2},l}^{n,L}}{\Delta t}  &=   - \left(1-\theta' \right) \frac{ q^{n}_{j,k+1,l} - q^{n}_{j,k,l}}{\Delta z^n_{j,k+\frac{1}{2},l}}+ \nonumber \\ 
 &+  \nu \left( \mathcal{L}^{wx\varphi}_{j,k+\frac{1}{2},l}\left[ w^{n,L},v^{n,L}\right] + \mathcal{L}^{wz\varphi}_{j,k+\frac{1}{2},l}\left[ \tilde{w}^{n+\theta} \right]  + \mathcal{L}^{wz\varphi}_{j,k+\frac{1}{2},l}\left[ \delta w^{n,L} \right] \right) ,\label{eq:cap5NumRadH}
\end{align}
\begin{align}
V_j \left( \tilde{p}^{n+1}+q^{n} \right) &= V_j \left( p^n  \right) - \Delta t \sum_{l=1}^{N_{\varphi}} \left[ \sum_{k=1}^{K^n_{j+\frac{1}{2},l}}a^n_{j+\frac{1}{2},k,l} \left( \theta \tilde{u}^{n+1}_{j+\frac{1}{2},k,l}  + \left(1-\theta\right) u^{n}_{j+\frac{1}{2},k,l} \right) - \right. \nonumber\\
& \left. \hspace{32mm} \sum_{k=1}^{K^n_{j-\frac{1}{2},l}}a^n_{j-\frac{1}{2},k,l} \left( \theta \tilde{u}^{n+1}_{j-\frac{1}{2},k,l}  + \left(1-\theta\right) u^{n}_{j-\frac{1}{2},k,l} \right)\right].\label{eq:cap5NumVolumeH}
\end{align}
It is important to notice that equations (\ref{eq:cap5NumAxialH})-(\ref{eq:cap5NumVolumeH}) depend on the non-hydrostatic pressure and velocity fields through explicit terms, while on the other hand implicit terms are hydrostatic. In the system above a different implicitness factor $\theta'$ can be used to discretize the non-hydrostatic terms 
in time. 
Then, in the second fractional step the \emph{non-hydrostatic problem} is solved, defined as the \emph{difference} between the original system of equations and the hydrostatic problem. 

\begin{figure}[bt] 
\centering
{\includegraphics[width=0.50\linewidth]{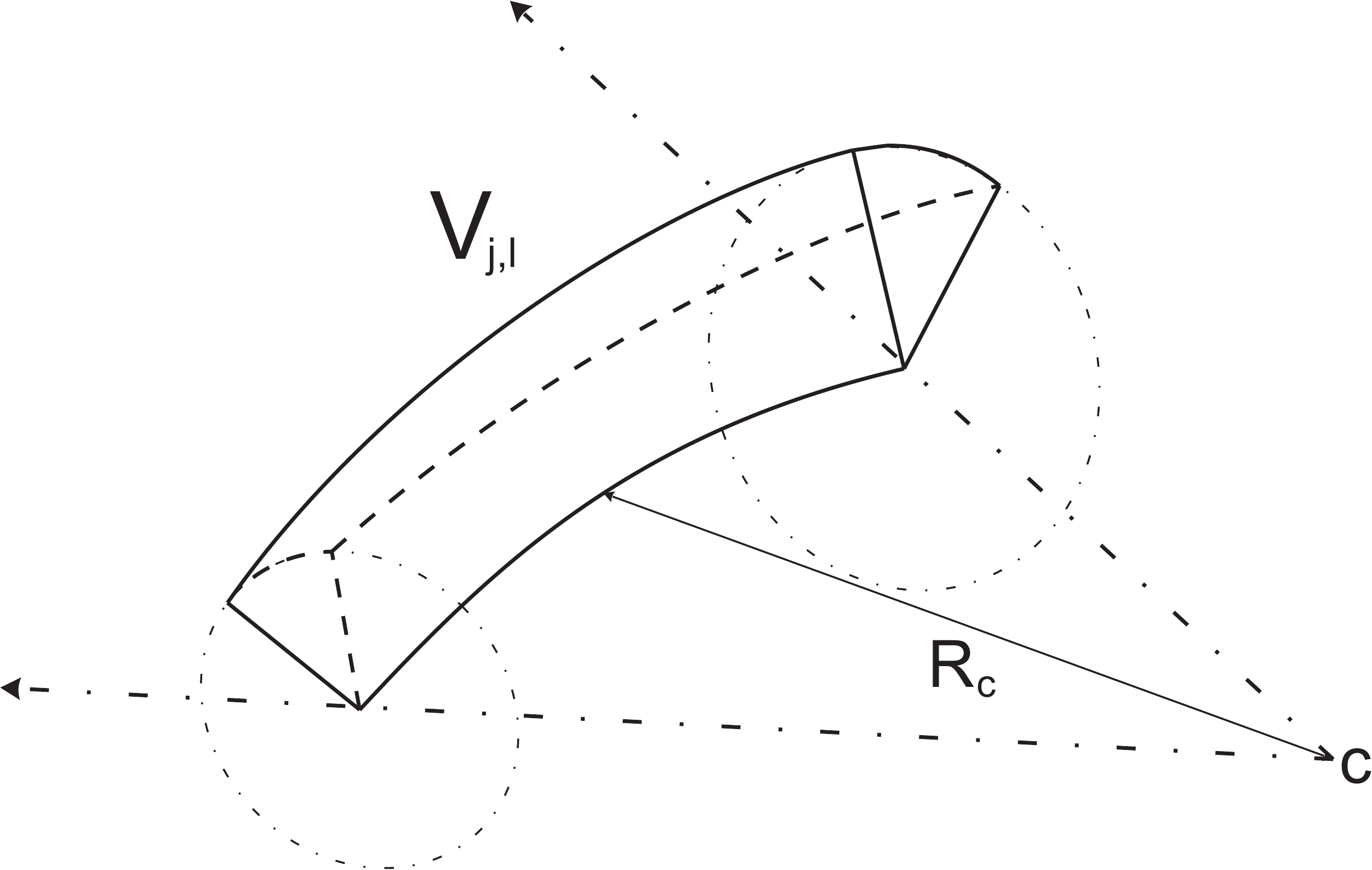}} 
\caption{Example of a curved slice. $\mathit{C}$ is the center of curvature. }\label{fig:cap5CurvedVol}
\end{figure}

Then, local and global mass conservation are ensured by a conservative finite-volume discretization of the continuity equation, i.e. 
\begin{align}
 \left(a^n_{j+\frac{1}{2},k,l}u^{n+1}_{j+\frac{1}{2},k,l} - a^n_{j-\frac{1}{2},k,l}u^{n+1}_{j-\frac{1}{2},k,l} \right)  +  \left(\mathit{b}^n_{j,k,l+\frac{1}{2}}v^{n+1}_{j,k,l+\frac{1}{2}} - \mathit{b}^n_{j,k,l-\frac{1}{2}}v^{n+1}_{j,k,l-\frac{1}{2}}  \right)  + & \nonumber \\
 \left(\mathit{c}^n_{j,k+\frac{1}{2},l}w^{n+1}_{j,k+\frac{1}{2},l} - \mathit{c}^n_{j,k-\frac{1}{2},l}w^{n+1}_{j,k-\frac{1}{2},l}  \right)  = &\; 0, \label{eq:cap5NumCont}
\end{align}
for each control volume $V_{j,k,l}$ below the radius, where $\mathit{b}^n_{j,k,l+1/2}$ and $\mathit{c}^n_{j,k+1/2,l}$ are the angular and radial surfaces of $V_{j,k,l}$, respectively, and at the walls a \emph{local} semi-implicit finite-volume expression of the kinematic boundary condition reads 
\begin{align}
V_{j,l} \left( p^{n+1} \right) &= V_{j,l} \left( p^n\right)  - \Delta t \sum_{k=1}^{K^n_{i,l}} \left[ \left( a^n_{j+\frac{1}{2},k,l} u^{n+\theta'}_{j+\frac{1}{2},k,l} - a^n_{j-\frac{1}{2},k,l} u^{n+\theta'}_{j-\frac{1}{2},k,l}  \right)    \right. \nonumber\\ 
  &  +  \left. \left(\mathit{b}^n_{j,k,l+\frac{1}{2}}v^{n+\theta'}_{j,k,l+\frac{1}{2}} - \mathit{b}^n_{j,k,l-\frac{1}{2}}v^{n+\theta'}_{j,k,l-\frac{1}{2}}  \right)  \right]\label{eq:cap5NumVolume2} 
\end{align}
where $V_{j,l}$ is the $\mathit{l}$-th slice of the volume of the $\mathit{j}$-th  axial segment $V_{j}$ (see Figure \ref{fig:cap5CurvedVol}).

\subsection{Solution algorithm}
\paragraph{First fractional step.} Equations (\ref{eq:cap5NumVolumeH}) and (\ref{eq:cap5NumAxialH}) constitute a mildly nonlinear system of at most $N_zN_xN_\varphi+N_x$ equations. This system has to be solved at each time step in order to calculate the new hydrostatic field variables $\tilde{u}^{n+1}_{j+1/2,k,l}$ and $\tilde{p}^{n+1}_j$ throughout the flow domain.

Upon multiplication by $a^n_{j+1/2,k,l} \Delta t \, \Delta \varphi$ and after including the boundary conditions (\ref{eq:NumBCuv}), equations (\ref{eq:cap5NumVolume})-(\ref{eq:cap5NumAxialH}) are first written in a compact matrix form as

\begin{align}
\doublehat{\mathbb{M}}^n_{j+\frac{1}{2}} \cdot {\tilde{\mathbf{U}}}^{n+1}_{j+\frac{1}{2}} &= \hat{\mathbf{G}}^n_{j+\frac{1}{2}} - \theta \Delta t \Delta \varphi \left[ \tilde{p}^{n+1}_{j+1} - \tilde{p}^{n+1}_{j}\right] \tilde{{\mathbf{A}}}^n_{j+\frac{1}{2}}\label{eq:cap4NumSystAxial} \\
V_j\left( \tilde{p}^{n+1} + q^n \right) &= V_j \left( {p}^n \right) - \Delta t \left[ \left(\hat{\mathbf{A}}^n_{j+\frac{1}{2}}\right)^T  \left( \theta {\tilde{\mathbf{U}}}^{n+1}_{j+\frac{1}{2}} +\left(1-\theta\right) {\mathbf{U}}^{n}_{j+\frac{1}{2}} \right)    \right. \nonumber \\
& \left. \hspace{22mm} -  \left({\hat{\mathbf{A}}^n_{j-\frac{1}{2}}}\right)^T  \left( \theta {\tilde{\mathbf{U}}}^{n+1}_{j-\frac{1}{2}} +\left(1-\theta\right) {\mathbf{U}}^{n}_{j-\frac{1}{2}} \right)   \right] 
\label{eq:cap4NumSystVolume}  
\end{align}
The four vectors ${\tilde{\mathbf{U}}}_{j+1/2}^{n+1}$, $\hat{\mathbf{A}}_{j+1/2}^n$, $\tilde{{\mathbf{A}}}_{j+1/2}^n$, $\hat{\mathbf{G}}_{j+1/2}^n$ and the matrix  $\doublehat{\mathbb{M}}^n_{j+1/2}$ are defined as follows: ${\tilde{\mathbf{U}}}_{j+1/2}^{n+1}$ collects all the discrete axial velocity components on the surface $j+1/2$; $\hat{\mathbf{A}}_{j+1/2}^{n}$ and  $\tilde{{\mathbf{A}}}_{j+1/2}^{n}$ collect the areas of each axial surface $a_{j+1/2,k,l}^n$ and the ratio $a_{j+1/2,k,l}^n/\Delta x_{j+1/2,k,l}$,  respectively; $\hat{\mathbf{G}}_{j+1/2}^n$ collects all the known explicit terms of equation (\ref{eq:cap5NumAxialH}); $\doublehat{\mathbb{M}}^n_{j+1/2}$ collects all the coefficients of the resulting linear system. 

The matrix $\doublehat{\mathbb{M}}_{j+1/2}$  
is \emph{pentadiagonal}, \emph{symmetric} and \emph{positive definite}. This means that  $\doublehat{\mathbb{M}}_{j+1/2}$ is \emph{non-singular} and invertible.  These properties are very important, because this kind of systems is efficiently solved by a conjugate-gradient method \cite{cgmethod}. 
Then, by multiplying equation (\ref{eq:cap4NumSystAxial}) formally with $\doublehat{\mathbb{M}}_{j+1/2}^{-1}$ from the left, yields  
\begin{gather}
{\tilde{\mathbf{U}}}^{n+1}_{j+\frac{1}{2}} =  \doublehat{\mathbb{M}}_{j+\frac{1}{2}}^{-1} \cdot \hat{\mathbf{G}}^n_{j+\frac{1}{2}} - \theta \Delta t \Delta \varphi \left[ \tilde{p}^{n+1}_{j+1} - \tilde{p}^{n+1}_{j}\right]  \doublehat{\mathbb{M}}_{j+\frac{1}{2}}^{-1} \cdot \tilde{{\mathbf{A}}}^n_{j+\frac{1}{2}} \label{eq:cap4NumSystAxialMinv}
\end{gather}
A substantial improvement of efficiency in solving the pentadiagonal linear systems of the form $\doublehat{\mathbb{M}}_{j+1/2} \hat{\xi} = \hat{\Lambda}$ of at most $N_zN_\varphi$ equations (where $\hat{\Lambda}$ is either $\hat{\mathbf{G}}^n_{j+1/2}$ or $\tilde{\mathbf{A}}^n_{j+1/2}$ is reached by firstly solving $N_\varphi$ independent tridiagonal systems $\hat{\mathbf{M}}'_{j+1/2,\cdot,l} \xi'_l  = \Lambda_l$  of $K^n_{j+1/2,l}$ equations obtained by neglecting the \emph{implicit}, \emph{angular} derivative operators.
If $\hat{\xi}$ is the tensor of rank two solution of the original system defined along the $(j+1/2)$-th cross section  (two-dimensional problem), then  $\xi'_l$ is the solution of the approximated system defined on the $l$-th slice of the $(j+1/2)$-th cross section (one-dimensional problem). Then if $\hat{\xi}'_{j+1/2}$ is the vector that 
collects all the approximate solutions $\xi'_l$ for $l=1$, $2$, $\cdots$, $N_\varphi$, the original problem reduces to solving the linear system  
\begin{gather}
\doublehat{\mathbb{M}}_{j+\frac{1}{2}}\cdot \hat{\delta\xi} = \hat{\Lambda}- \doublehat{\mathbb{M}}_{j+\frac{1}{2}}\cdot \hat{\xi}',\label{eq:perturbative}
\end{gather}
and the solution of the original two dimensional system is given by $\hat{\xi} = \hat{\xi}' + \hat{\delta \xi}$. This \emph{perturbative} way of solving the original system becomes particularly useful when axial symmetry is  fulfilled, because the one-dimensional approximate solution $\hat{\xi}'$ becomes \emph{exact} and thus no iteration is needed in solving the two dimensional problem (\ref{eq:perturbative})  
because $\hat{\delta \xi}=0$. 
and the solution becomes very fast.  

Formal substitution of ${\tilde{\mathbf{U}}}^{n+1}_{j+1/2}$ of equation (\ref{eq:cap4NumSystAxialMinv}) into  (\ref{eq:cap4NumSystVolume}) yields the discrete wave equation for the hydrostatic pressure 
\begin{align}
V_j\left( \tilde{p}^{n+1} + q^n \right) -\Delta t^2 \theta^2 \Delta \varphi \left[  \left( \hat{\mathbf{A}}_{j+\frac{1}{2}}^n \right)^T \cdot 
  \doublehat{\mathbb{M}}_{j+\frac{1}{2}}^{-1} \cdot \tilde{{\mathbf{A}}}_{j+\frac{1}{2}}^n    \left(\tilde{p}^{n+1}_{j+1} - \tilde{p}^{n+1}_{j} \right)  - \right. & \nonumber\\
	\left. \left( \hat{\mathbf{A}}_{j-\frac{1}{2}}^n \right)^T \cdot 
  \doublehat{\mathbb{M}}_{j-\frac{1}{2}}^{-1} \cdot \tilde{{\mathbf{A}}}_{j-\frac{1}{2}}^n    \left(\tilde{p}^{n+1}_{j} - \tilde{p}^{n+1}_{j-1} \right)	\right] &= d_j^n \label{eq:cap4NumSystP}
\end{align}
where
\begin{align*}
d_j^n &= V_j \left( {p}^{n} \right)  - \Delta t\left[ \left( \hat{\mathbf{A}}^n_{j+\frac{1}{2}} \right)^T \cdot  \left( \theta   \doublehat{\mathbb{M}}_{j+\frac{1}{2}}^{-1} \cdot \hat{\mathbf{G}}_{j+\frac{1}{2}}^n   +\left(1-\theta\right) \hat{\mathbf{U}}_{j+\frac{1}{2}}^n \right)  - \right.\\
& \left. \hspace{25mm} \left(\hat{\mathbf{A}}^n_{j-\frac{1}{2}}\right)^T \cdot  \left( \theta   \doublehat{\mathbb{M}}_{j-\frac{1}{2}}^{-1} \cdot \hat{\mathbf{G}}_{j-\frac{1}{2}}^n   +\left(1-\theta\right) \hat{\mathbf{U}}_{j-\frac{1}{2}}^n \right)\right].
\end{align*}
Equation (\ref{eq:cap4NumSystP}) can be assembled into a sparse, \emph{mildly nonlinear} system of at most $N_x$ equations for $\tilde{p}^{n+1}_j$, $j=1$, $2$, $\ldots, N_x$. This system is efficiently solved by a Newton-type iterative algorithm whose details are given in \cite{Casulli:2009,Casulli:2012b}.

Once the new hydrostatic pressure $\tilde{p}_j^{n+1}$ has been computed, the hydrostatic axial velocities are readily obtained from equation (\ref{eq:cap4NumSystAxial}), which now represents a set of $N_x$ independent, \emph{linear} pentadiagonal systems of at most $N_zN_\varphi$ equations, for all $j=1,2,...,N_x$.

Equations (\ref{eq:cap5NumTangH}) constitute essentially a linear system of at most $N_zN_{\varphi}N_x$ equations. It is decomposed and solved as a system of $N_x$ independent, linear \emph{pentadiagonal} systems
that can be rewritten in compact matrix-vector form as 
\begin{gather}
\doublehat{\mathbb{N}}^n_{j} \cdot {\tilde{\mathbf{V}}}^{n+1}_{j} = \hat{\mathbf{H}}^n_{j}, \label{eq:cap4NumSystTang}
\end{gather}
where the two vectors ${\tilde{\mathbf{V}}}_{j}^{n+1}$, $\hat{\mathbf{H}}_{j}^n$ and the matrix $\doublehat{\mathbb{N}}_{j}$ are defined as follows: $\hat{\tilde{\mathbf{V}}}_{j}^{n+1}$ collects all the discrete tangential 
velocity components; $\hat{\mathbf{H}}_{j}^{n}$ collects all the known explicit terms of equation (\ref{eq:cap5NumTangH}) and $\doublehat{\mathbb{N}}$ collects all the coefficient of the resulting linear system.

Also $\doublehat{\mathbb{N}}$ is pentadiagonal. Similarly to  $\doublehat{\mathbb{M}}$, the matrix $\doublehat{\mathbb{N}}$ is 
invertible and the tangential velocities are readily obtained from equation (\ref{eq:cap4NumSystTang}) 
, for all $j=1,2,...,N_x$. 

Now, the boundary condition (\ref{eq:cap5wR}) needs to be rewritten in hydrostatic form as 
\begin{gather}
\tilde{w}^{n+1}_{j,K+\frac{1}{2},l} = \frac{\tilde{R}_{j,l}^{n+1} - R_{j,l}^n}{\Delta t} = \frac{1}{ \beta} \frac{\tilde{p}_j^{n+1}-\tilde{p}_j^n}{\Delta t}.\label{eq:cap5wRH}
\end{gather}
where $\tilde{R}^{n+1}=\tilde{R}(\tilde{p}^{n+1} + q^n)$ is obtained from the equation of state (\ref{eq:cap5laplace}).
When the mildly non-linear system is solved, axial and tangential components are computed and boundary condition (\ref{eq:cap5wR}) is determined. Then, equations (\ref{eq:cap5NumRadH}) constitute a linear system of at most $N_zN_{\varphi}N_x$ equations. 
With the same procedure, 
equation (\ref{eq:cap5NumRad}) can be rewritten in compact form as 
\begin{gather}
\doublehat{\mathbb{O}}^n_{j} \cdot {\tilde{\mathbf{W}}}^{n+1}_{j} = \hat{\mathbf{L}}^n_{j}, \label{eq:cap5NumSystRad}
\end{gather}
where the two vectors $\hat{\mathbf{W}}_{j}$, $\hat{\mathbf{L}}_{j}$ and the matrix $\doublehat{\mathbb{O}}_{j}$ are defined as follows: ${\tilde{\mathbf{W}}}^{n+1}_{j}$ collects all the hydrostatic radial velocity components;  $\hat{\mathbf{L}}^n_{j}$ collects all known terms of equation (\ref{eq:cap5NumRadH}) and $\doublehat{\mathbb{O}}^n_{j}$ collects the coefficients of the resulting linear system. 

$\doublehat{\mathbb{O}}$ is a symmetric positive-definite matrix
, hence the radial velocities are readily obtained from equation (\ref{eq:cap5NumSystRad}),
 for all $j=1,2,...,N_x$.  
In the present formulation, the aforementioned perturbative approach is also used for obtaining ${\tilde{\mathbf{V}}}^{n+1}$ and ${\tilde{\mathbf{W}}}^{n+1}$.

\paragraph{Second fractional step.} 
The non-hydrostatic momentum equations, obtained from the difference between equations (\ref{eq:cap5NumAxial})-(\ref{eq:cap5NumRad}) and (\ref{eq:cap5NumAxialH})-(\ref{eq:cap5NumRadH}), can be summarized in the sequent system of equations for the components of the velocity field
\begin{align}
u_{j+\frac{1}{2},k,l}^{n+1}  &= \tilde{u}_{j+\frac{1}{2},k,l}^{n+1}  - \theta' \Delta t\frac{ q^{n+1}_{j+1,k,l} - q^{n+1}_{j,k,l}}{\Delta x_{j+\frac{1}{2},k,l}},\label{eq:cap5NumAxialNH2}\\
v_{j,k,l+\frac{1}{2}}^{n+1}  &= \tilde{v}_{j,k,l+\frac{1}{2}}^{n+1} - \theta' \Delta t \frac{ q^{n+1}_{j,k,l+1} - q^{n+1}_{j,k,l}}{z^n_{j,k,l+\frac{1}{2}}\Delta \varphi},\label{eq:cap5NumTangNH2}\\
w_{j,k+\frac{1}{2},l}^{n+1}  &= \tilde{w}_{j,k+\frac{1}{2},l}^{n+1}  - \theta'  \Delta t \frac{ q^{n+1}_{j,k+1,l} - q^{n+1}_{j,k,l}}{\Delta z^n_{j,k+\frac{1}{2},l}},\label{eq:cap5NumRadNH2}
\end{align}
that is closed by a finite volume discretization of equations (\ref{eq:cap5NumCont})-(\ref{eq:cap5NumVolume2}).
By using (\ref{eq:cap5NumCont}), equation (\ref{eq:cap5NumVolume2}) reduces to
\begin{align}
V_{j,l} \left( p^{n+1} \right) &= \delta_{j,l}  - \theta' \Delta t \left[  \left( a^n_{j+\frac{1}{2},K,l} u^{n+1}_{j+\frac{1}{2},K,l} - a^n_{j-\frac{1}{2},K,l} u^{n+1}_{j-\frac{1}{2},K,l}  \right) +\right. \nonumber\\
 &   + \left(\mathit{b}^n_{j,K,l+\frac{1}{2}}v^{n+1}_{j,K,l+\frac{1}{2}} - \mathit{b}^n_{j,K,l-\frac{1}{2}}v^{n+1}_{j,K,l-\frac{1}{2}}  \right) -  \left. \mathit{c}^n_{j,K-\frac{1}{2},l}w^{n+1}_{j,K-\frac{1}{2},l}   \right], \label{eq:cap5NumVolume3}
\end{align}
with
\begin{align}
\delta_{j,l} &=  V_{j,l} \left( p^{n} \right) - \left( 1- \theta'\right) \Delta t \sum_{k=1}^{K^n_i} \left[ \left( a^n_{j+\frac{1}{2},k,l} u^{n}_{j+\frac{1}{2},k,l} - a^n_{j-\frac{1}{2},k,l} u^{n}_{j-\frac{1}{2},k,l}  \right)    \right. \nonumber\\ 
  &  +  \left. \left(\mathit{b}^n_{j,k,l+\frac{1}{2}}v^{n}_{j,k,l+\frac{1}{2}} - \mathit{b}^n_{j,k,l-\frac{1}{2}}v^{n}_{j,k,l-\frac{1}{2}}  \right)  \right]. \label{eq:cap5NumVolume4}
\end{align}

At this point, if the expression for the new velocities (\ref{eq:cap5NumAxialNH2})-(\ref{eq:cap5NumRadNH2}) is substituted into (\ref{eq:cap5NumCont} ) and (\ref{eq:cap5NumVolume3}), the result is a system of equations for the  non-hydrostatic pressure $\mathit{q}$. The following discrete equations are obtained: 
\begin{align}
&\theta' \Delta t \left[\frac{a^n_{j+\frac{1}{2},k,l}}{\Delta x_{j+\frac{1}{2},k,l}} \left( q^{n+1}_{j+1,k,l}-q^{n+1}_{j,k,l} \right) -\frac{a^n_{j-\frac{1}{2},k,l}}{\Delta x_{j-\frac{1}{2},k,l}} \left( q^{n+1}_{j,k,l}-q^{n+1}_{j-1,k,l} \right)   \right. + \nonumber \\
&+\frac{\mathit{b}^n_{j,k,l+\frac{1}{2}}}{ z_{j,k} \Delta \varphi} \left( q^{n+1}_{j,k,l+1}-q^{n+1}_{j,k,l} \right) -\frac{\mathit{b}^n_{j,k,l-\frac{1}{2}}}{z_{j,k} \Delta \varphi} \left( q^{n+1}_{j,k,l}-q^{n+1}_{j,k,l-1} \right)   + \nonumber \\ 
& 	+\left. \frac{\mathit{c}^n_{j,k+\frac{1}{2},l}}{\Delta z_{j+\frac{1}{2},k}} \left( q^{n+1}_{j,k+1,l}-q^{n+1}_{j,k,l} \right) -\frac{\mathit{c}^n_{j,k-\frac{1}{2},l}}{\Delta z_{j,k-\frac{1}{2}}} \left( q^{n+1}_{j,k,l}-q^{n+1}_{j,k-1,l} \right)    \right] = \nonumber \\
	 &= \phantom{\left[!^!\right]}  a^n_{j+\frac{1}{2},k,l}\tilde{u}^{n+1}_{j+\frac{1}{2},k,l}   + \mathit{b}^n_{j,k,l+\frac{1}{2}}\tilde{v}^{n+1}_{j,k,l+\frac{1}{2}} + \mathit{c}^n_{j,k+\frac{1}{2},l}\tilde{w}^{n+1}_{j,k+\frac{1}{2},l}  + \nonumber \\
& \phantom{\left[!^!\right]} - a^n_{j-\frac{1}{2},k,l}\tilde{u}^{n+1}_{j-\frac{1}{2},k,l}  - \mathit{b}^n_{j,k,l-\frac{1}{2}}\tilde{v}^{n+1}_{j,k,l-\frac{1}{2}}  - \mathit{c}^n_{j,k-\frac{1}{2},l}\tilde{w}^{n+1}_{j,k-\frac{1}{2},l} . \label{eq:cap5q}
\end{align}
for $k = 1$, $2$, $K^n_{j,l}-1$, and
\begin{align}
&\theta' \Delta t \left[\frac{a^n_{j+\frac{1}{2},K,l}}{\Delta x_{j+\frac{1}{2},K,l}} \left( q^{n+1}_{j+1,K,l}-q^{n+1}_{j,K,l} \right) -\frac{a^n_{j-\frac{1}{2},K,l}}{\Delta x_{j-\frac{1}{2},K,l}} \left( q^{n+1}_{j,K,l}-q^{n+1}_{j-1,K,l} \right)   \right. + \nonumber \\
& +\frac{\mathit{b}^n_{j,K,l+\frac{1}{2}}}{ z_{j,K} \Delta \varphi} \left( q^{n+1}_{j,K,l+1}-q^{n+1}_{j,K,l} \right) -\frac{\mathit{b}^n_{j,K,l-\frac{1}{2}}}{z_{j,K} \Delta \varphi} \left( q^{n+1}_{j,K,l}-q^{n+1}_{j,K,l-1} \right)   + \nonumber \\ 
& 	\left. -\frac{\mathit{c}^n_{j,K-\frac{1}{2},l}}{\Delta z_{j,K-\frac{1}{2}}} \left( q^{n+1}_{j,K,l}-q^{n+1}_{j,K-1,l} \right)    \right] - \frac{ V_{j,l} \left( \tilde{p}^{n+1}+ q^{n+1} \right)}{\theta' \Delta t }= \nonumber \\
	 &= - \frac{ \delta_{j,l}}{\theta' \Delta t } +  \left[ a^n_{j+\frac{1}{2},K,l}\tilde{u}^{n+1}_{j+\frac{1}{2},K,l}   + \mathit{b}^n_{j,K,l+\frac{1}{2}}\tilde{v}^{n+1}_{j,K,l+\frac{1}{2}}  + \right. \nonumber \\
&\left. - a^n_{j-\frac{1}{2},K,l}\tilde{u}^{n+1}_{j-\frac{1}{2},K,l}  - \mathit{b}^n_{j,K,l-\frac{1}{2}}\tilde{v}^{n+1}_{j,K,l-\frac{1}{2}}  - \mathit{c}^n_{j,K-\frac{1}{2},l}\tilde{w}^{n+1}_{j,K-\frac{1}{2},l} \right]. \label{eq:cap5qV}
\end{align}
Equations (\ref{eq:cap5q})-(\ref{eq:cap5qV}) constitute a \emph{seven-diagonal mildly nonlinear system} of $N_xN_zN_\varphi $ unknowns that is efficiently solved by a Newton-type iterative method. Since the Jacobian of the nonlinear system (\ref{eq:cap5q})-(\ref{eq:cap5qV}) is 
again symmetric and positive definite, a conjugate gradient method can be used at each Newton step. 

That is how the non-hydrostatic pressure is computed. Subsequently, the velocity components are directly updated accordingly through (\ref{eq:cap5NumAxialNH2})-(\ref{eq:cap5NumRadNH2}). Finally, the new radii are obtained from 
the equation of state (\ref{eq:cap5laplace}), the new indices $K^{n+1}_{j,l}$ and the spatial position of the grid nodes are updated accordingly. 

Note that when axial symmetry is assumed, by setting $N_{\varphi}=1$, then a consistent \emph{two-dimensional}, \emph{tangentially averaged}, \emph{non hydrostatic} model is obtained. If then, non hydrostatic effects become negligible, the current formulation reduces exactly to the \emph{two-dimensional} \emph{hydrostatic} model of \cite{Casulli:2012}, and a 
\emph{one dimensional} \emph{sectionally averaged} model, commonly used to simulate blood flow in complex  arterial systems (see, e.g., \cite{Olufsen:2000,Sherwin:2003,Sherwin:2003b,MuellerJCP}), 
is automatically contained as well. 

The fact that a consistent two-dimensional or one-dimensional scheme can be derived from the proposed numerical method as particular case is a  valuable feature of the present formulation. 
This property leads to a general algorithm that can solve non-hydrostatic three-dimensional, as well as hydrostatic three-, two- and one-dimensional problems as a particular case. 

More precisely, when the method proposed in this article is applied to a typical model of the human cardiovascular system that includes large and small arteries, from complex to simplified geometries, a substantial computational  simplification is achieved because blood flow in nearly straight large and medium sized arteries are properly modeled within a two-dimensional approach, whereas small and very small arteries are automatically represented by a  simple and consistent one-dimensional formulation. This means that those simplified two- or one-dimensional branches of the full system get their own two- or one-dimensional representation, without any special or artificial  treatment at the interface. 
For other successful numerical approaches documented in literature that concern the coupling of three-dimensional and one-dimensional models for blood flow in arterial systems
see \cite{Formaggia2012} and for particular boundary conditions used in biomechanics the reader is referred to \cite{Vergara2008,Porpora,Wall2013}. 



\section{Numerical tests}

The aforementioned numerical scheme for three-dimensional non-hydrostatic blood flow in compliant vessels is first applied to three test problems with known analytical solution for circular and elliptical cross sections. Then,  steady and pulsatile flows in a uniformly curved vessel are chosen as non-hydrostatic test problems: a comparison with experimental and numerical results of  \cite{Timite:2010} is presented. Finally, an estimate of the efficiency  of the presented numerical method is given by measuring the computing time needed for the simulation of one of the test problems.

\subsection{Steady flow in an elastic tube}

\begin{table}[t]
\begin{center}
\begin{tabular}{c c c c c c c c } 
\hline
$N_x$ & $N_z$ & $N_\varphi$ & $N_t$ &  $\epsilon^{u,N_t}_{L_2}$ & $\mathpzc{O}^u_{L_2}$ &  $\epsilon^{R,N_t}_{L_2}$ & $\mathpzc{O}^R_{L_2}$\\
\hline
$\phantom{0}50$ &  $\phantom{0}25$ & $\phantom{0}30$ & $100$ & $3.9260E-04$ &       & $7.2111E-05$ & \\
$100$ 					&  $\phantom{0}50$ & $\phantom{0}30$ & $100$ & $7.7510E-05$ & $2.3$ & $1.4446E-05$ & $2.3$\\
$200$ 					&    $100$         & $\phantom{0}30$ & $100$ & $1.5831E-05$ & $2.3$ & $3.3021E-06$ & $2.1$ \\
$400$ 					&    $200$         & $\phantom{0}30$ & $100$ & $3.6767E-06$ & $2.1$ & $7.8251E-07$ & $2.1$ \\
$800$ 					&    $400$         & $\phantom{0}30$ & $100$ & $8.6881E-07$ & $2.1$ & $1.9160E-07$ & $2.0$ \\%
\hline
\end{tabular}
\end{center}
\caption{Numerical convergence results in $L_2$ norm for $u(x,z)$ and for $R(x)$ for the steady flow problem through an elastic vessel.}
\label{tab:cap4Steady}
\end{table}

\begin{figure}[bt] 
\centering 
			\begin{subfigure}[b]{0.45\textwidth}
			\centering
			\includegraphics[width=\textwidth]{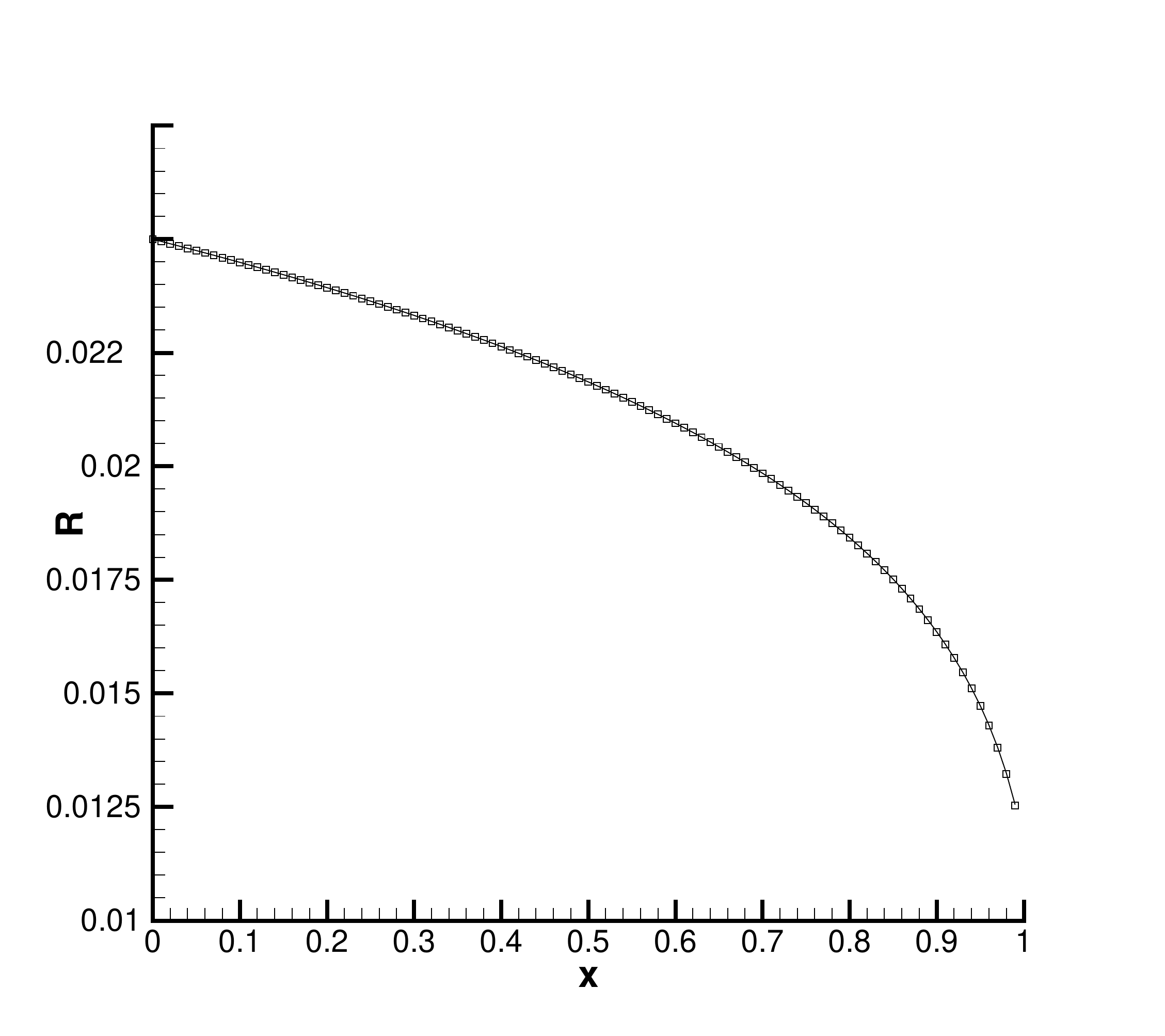}
			\end{subfigure}\;
			\begin{subfigure}[b]{0.45\textwidth}
			\centering
			\includegraphics[width=\textwidth]{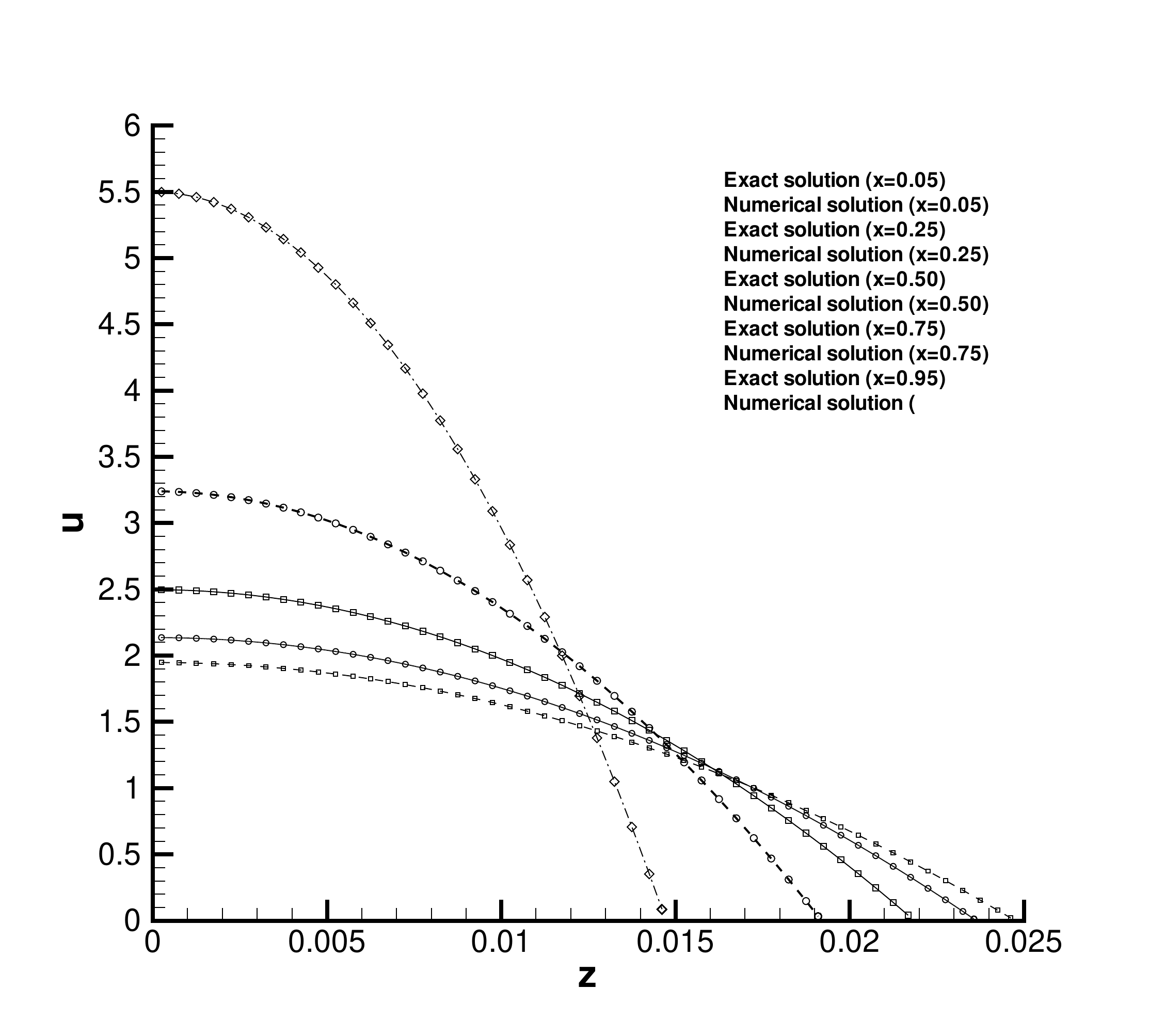}
			\end{subfigure}
\caption{Comparison of the numerical results obtained at time $t_e=10.0$ with the exact steady solution. Shape of the elastic tube (left) and selected velocity profiles at different axial positions (right). The velocity profile of only a fixed angular coordinate ($\varphi=180^{\circ}$) is shown, because of the axial symmetry of the numerical solution.}\label{fig:cap4Steady}
\end{figure}

At first, the correctness of the behaviour of the compliant wall within the proposed algorithm is tested by the convergence to the steady flow through an elastic tube with circular cross section, starting with a fluid at rest. 
Assuming axial symmetry and the pressure to be hydrostatic ($\mathit{q}=0$ everywhere) and neglecting the nonlinear convective terms in the momentum equation (\ref{eq:cap5NSu}), an exact steady solution of (\ref{eq:cap5NSu})-(\ref{eq:cap5BCw}) 
is given according to \cite{Fung:2010,Casulli:2012} by 
\begin{equation}
 u(x,z) = \frac{2Q}{\pi R^4(x)} \left(R^2(x)-z^2\right), \qquad 
 R(x) = \sqrt[5]{R_0^5-\frac{40 \nu Q}{\pi \beta} x}. 
\label{eqn.steady.uR}  
\end{equation} 
The chosen parameters for the present test are $\nu = 10^{-3}$, $\beta = 2500$, $p_{\text{ext}} = 0$, $R_0 = 0.025$, $Q=0.001875$, and $L=1$. Numerically, a transient solution is generated over a sufficiently long time interval, from the starting time $t = 0$, with initial conditions 
\begin{gather*}
u(x,z,\varphi,0) = 0 \;\;\; \text{and} \;\;\; R(x,\varphi,0) = R_0.
\end{gather*}
Then, for times $t > 0$, the boundary conditions are given by specifying the exact parabolic velocity profile at the inlet ($x = 0$), and the exact pressure according at the outlet ($x = L$). 

The computational domain is discretized with $N_x = 100$ segments in the axial direction, $N_z = 50$ rings are used along the radial direction to discretize the reference radius $R_0$ and  $N_\varphi=30$ uniform slices are used to discretize the angular coordinate in the interval $\left[\right.0,2\pi\left.\right)$. Assuming that the steady state is reached  at the final time $t_e= 10$, by using $\theta = 1$
, the simulation is advanced for $N_t = 100$ time steps with a time-step size $\Delta t = t_e / N_t$. 
 The resulting tube radius obtained at $t = t_e$ and some representative velocity profiles at different axial locations are illustrated in Figure \ref{fig:cap4Steady}. Symmetry is reproduced \emph{exactly}, hence, results are shown in two dimensional space at a fixed angular coordinate $\varphi=180^\circ$.  The classical parabolic Hagen-Poiseuille profile for the velocity is well reproduced, and an overall excellent agreement between the numerical results and the exact solution is clearly shown. 

The analytical solution of the present test problem is sufficiently smooth and, consequently, the order of accuracy of the proposed algorithm can be numerically determined by successively refining the spatial grid size. To this purpose, the discrete $\mathit{L}_2$ error norms for the axial velocity and for the radius, respectively, are evaluated with the three dimensional extension of the formula given in \cite{Casulli:2012}.
The errors $\epsilon^{u,n}_{\mathit{L}_2}$ and $\epsilon^{R,n}_{\mathit{L}_2}$ are computed by using a sequence of successively refined meshes obtained with $N_x= 100$, $200$, $400$, and  $800$ and $N_z = 50$, $100$, $200$,  and  $400$, respectively. Because the test problem depends only on the axial and radial coordinates, the angular mesh and the time step are kept constant.  
The convergence results listed in Table \ref{tab:cap4Steady} indicate that the designed second order of accuracy, $\mathpzc{O}^u_{\mathit{L}_2}$ and $\mathpzc{O}^R_{\mathit{L}_2}$, is achieved for this steady test problem.

\subparagraph{Non-hydrostatic corrections.} It is interesting to see what happens if non-hydrostatic corrections are included in the simulation of the steady problem above. 
By using $N_t = 100$ time steps with a time-step size $\Delta t = t_e / N_t$ and $\theta = 1$, the errors $\epsilon^{u,n}_{L_2}$ and $\epsilon^{R,n}_{L_2}$ are computed by using a sequence of successively refined meshes obtained with $N_x= 25$, $50$, $100$, $200$, and  $400$, and $N_z = 10$, $20$, $40$, $80$ and $160$, respectively. The time-step size and the angular discretization number are kept constant $\Delta t = 0.1$ and $N_\varphi=10$.

\begin{figure}[btp]
\centering 
			\includegraphics[width=0.9\textwidth]{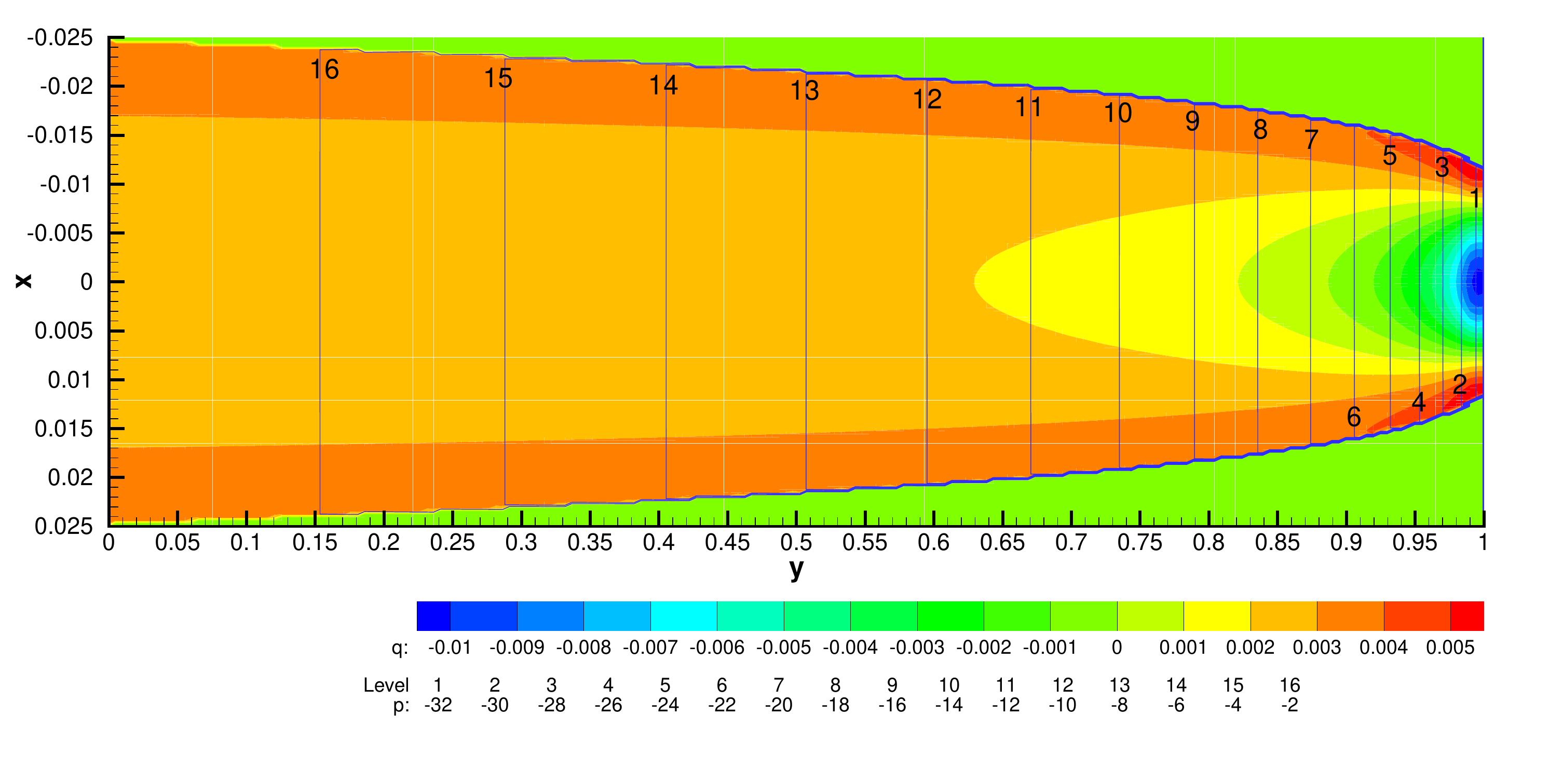}\\
\caption{Non-hydrostatic pressure distribution along an axial section for the steady state flow in an elastic tube. $N_x=200$, $N_z=80$, $N_\varphi=20$, $N_t=100$ and $t_e=20.0$. $\mathit{x}$ and $\mathit{y}$ are Cartesian coordinates}\label{fig:cap5SteadyQc} 
\end{figure}

\begin{figure}[btp]
\centering 
			\begin{tabular}{cc} 
			\includegraphics[width=0.4\textwidth]{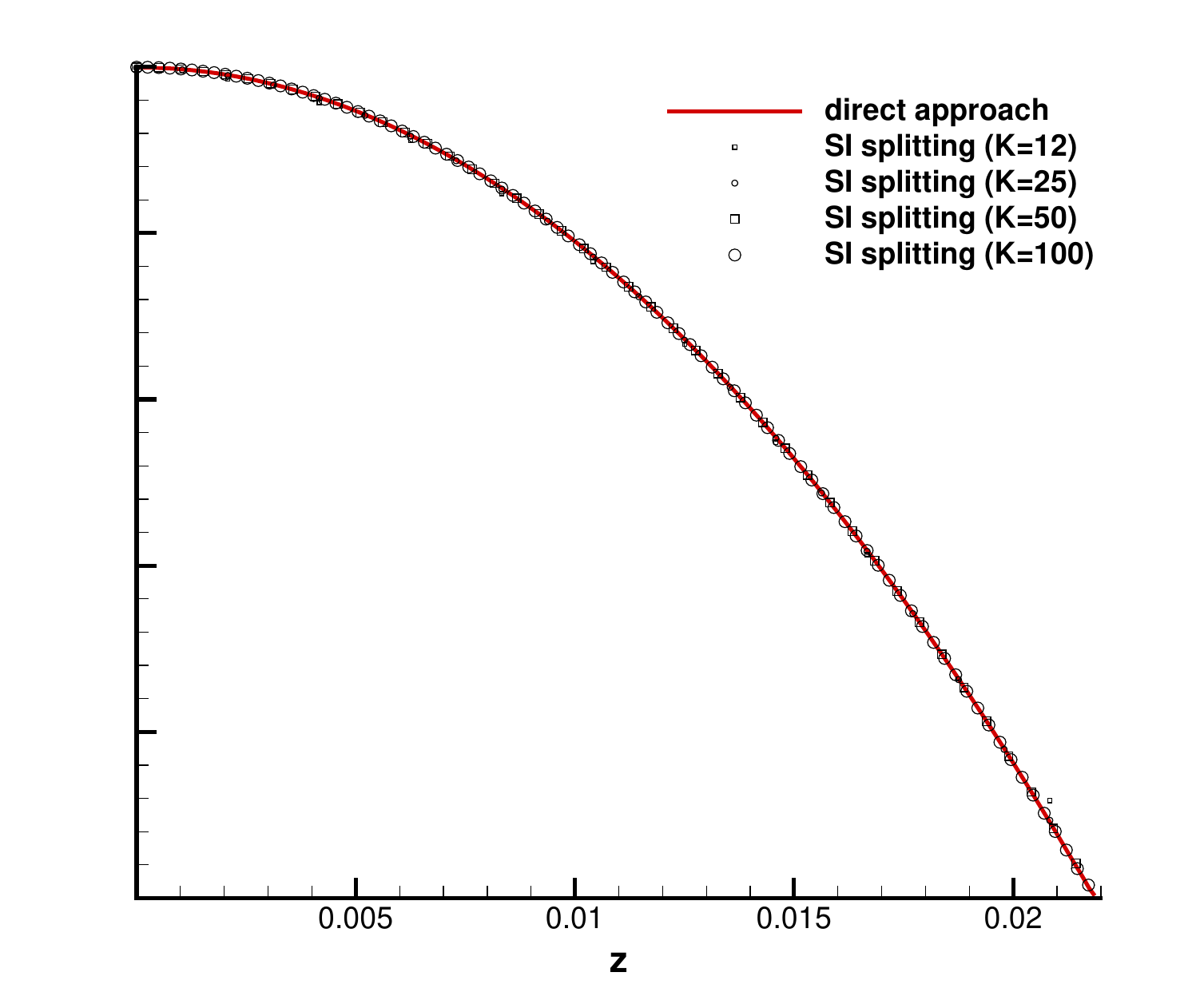} & 
			\includegraphics[width=0.4\textwidth]{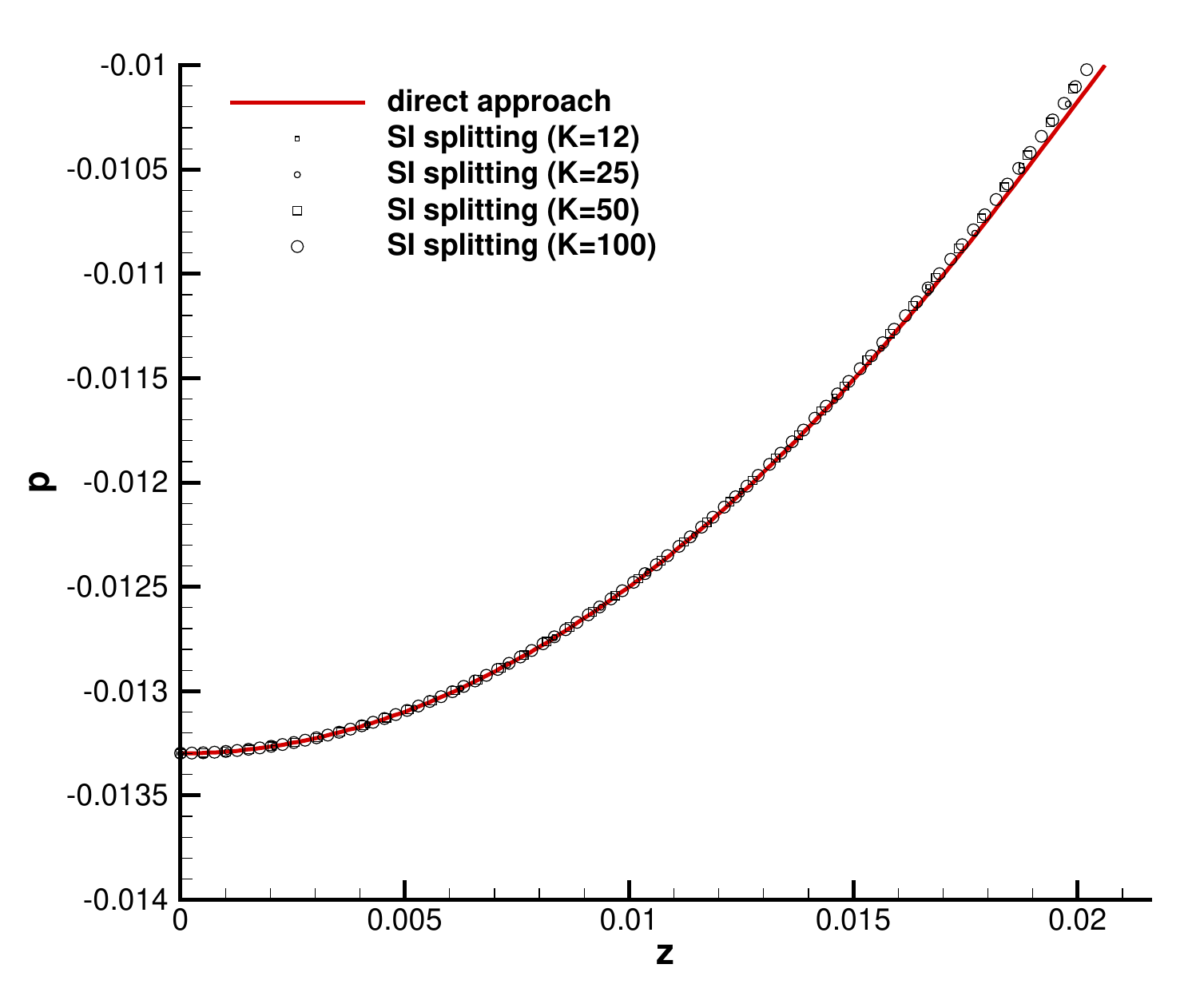} 
			\end{tabular} 
\caption{Radial profiles of axial velocity $u$ (left) and pressure $p$ (right) at $x=L/2$ computed with a direct approach and the present semi-implicit splitting method with various 
radial grid resolutions.}\label{fig.nhrad} 
\end{figure}

Figure \ref{fig:cap5SteadyQc} shows the distribution of the non-hydrostatic pressure $\mathit{q}$ along the axial section computed with the semi-implicit method proposed in this article. For comparison, we solve the same problem again with a \textit{direct approach} that solves the incompressible Navier-Stokes equations on a fine grid with 200 radial layers in the \textit{deformed geometry} given by 
(\ref{eqn.steady.uR}). The radial profiles for the axial velocity $u$ and the pressure distribution are shown in Fig. \ref{fig.nhrad},  
where one can note a very good agreement between the direct approach and the non-hydrostatic splitting method proposed in this paper.


\subsection{Womersley profiles}
\label{sec.wom} 

\begin{table}[b]
\begin{center}
\begin{tabular}{c c c c c c} 
\hline
$N_x$ 						 & $N_z$ 										  & $N_t$ & $N_\varphi$ &  $\epsilon^{u,N_t}_{L_2}$ & $\mathpzc{O}^u_{L_2}$ \\ 
\hline
$\phantom{0\,}100$ &    $\phantom{0}25$ 			  & $\phantom{0}25$ & $30$ & $9.8753E-05$ &       \\
$\phantom{0\,}200$ &    $\phantom{0}50$ 			  & $\phantom{0}50$ & $30$ & $2.7400E-05$ & $1.9$ \\ 
$\phantom{0\,}400$ &    $100$ 			  			    & $100$ 					& $30$ & $6.9837E-06$ & $2.0$ \\ 
$\phantom{0\,}800$ &    $200$         					& $200$ 					& $30$ & $1.7532E-06$ & $2.0$ \\
\hline
\end{tabular}
\end{center}
\caption{Numerical convergence results in $L_2$ error norm for $u(x,z,\varphi)$ at time $t_e=2.0$ (unsteady problem).}
\label{tab:cap5Unsteady}
\end{table}

Next, the performance of the proposed semi-implicit scheme is tested against an oscillating flow through a straight rigid tube. The flow is driven by a sinusoidal pressure gradient 
which is imposed at the ends of a tube of length $L$.  
According to Womersley \cite{Womersley:1955} the axial velocity profile is uniform in the axial $x$ direction and can be written in explicit form in terms of the zeroth order Bessel function  with reference to only one parameter, the Womersley number $\alpha = R \sqrt{\omega/\nu}$ (see \cite{Casulli:2012} for details).
For the present test, the chosen parameters are $L = 1$, $R=0.025$, $\hat{P}=1000$, $\rho = 1000$, $\omega = 2 \pi$, and $\beta = 10^{12}$ so that the tube wall is sufficiently rigid. Moreover, according to Womersley's approximations, the nonlinear advective terms are neglected. 

Depending on the viscosity, two cases are selected. Firstly, a low Reynolds number $R_e=50$ (based on the tube diameter $D =2R = 0.05$) is obtained by choosing $\nu = 10^{-3}$. In this regime the resulting viscous effects dominate the entire cross section, and the velocity profile resembles the parabolic Hagen-Poiseuille flow. In the second case, by choosing $\nu = 10^{-5}$, the flow is characterized by a Reynolds number $R_e = 5\,000$, and the resulting viscous effects are basically confined to a boundary layer close to the walls. In this regime, the velocity profile is essentially flat (inviscid) at the center of the cross section, whereas a sharp boundary layer with high velocity gradients develops near the walls.
The initial conditions are taken to be $u(x,z,\varphi,0)$ from the Womersley solution, and $R(x,0)=R_0$. Then, the time dependent boundary conditions are specified at the two ends of the tube by applying, for times $t > 0$, an oscillating pressure gradient.

\begin{figure} 
\centering 
			\begin{subfigure}[b]{0.45\textwidth}
			\centering
			\includegraphics[width=\textwidth]{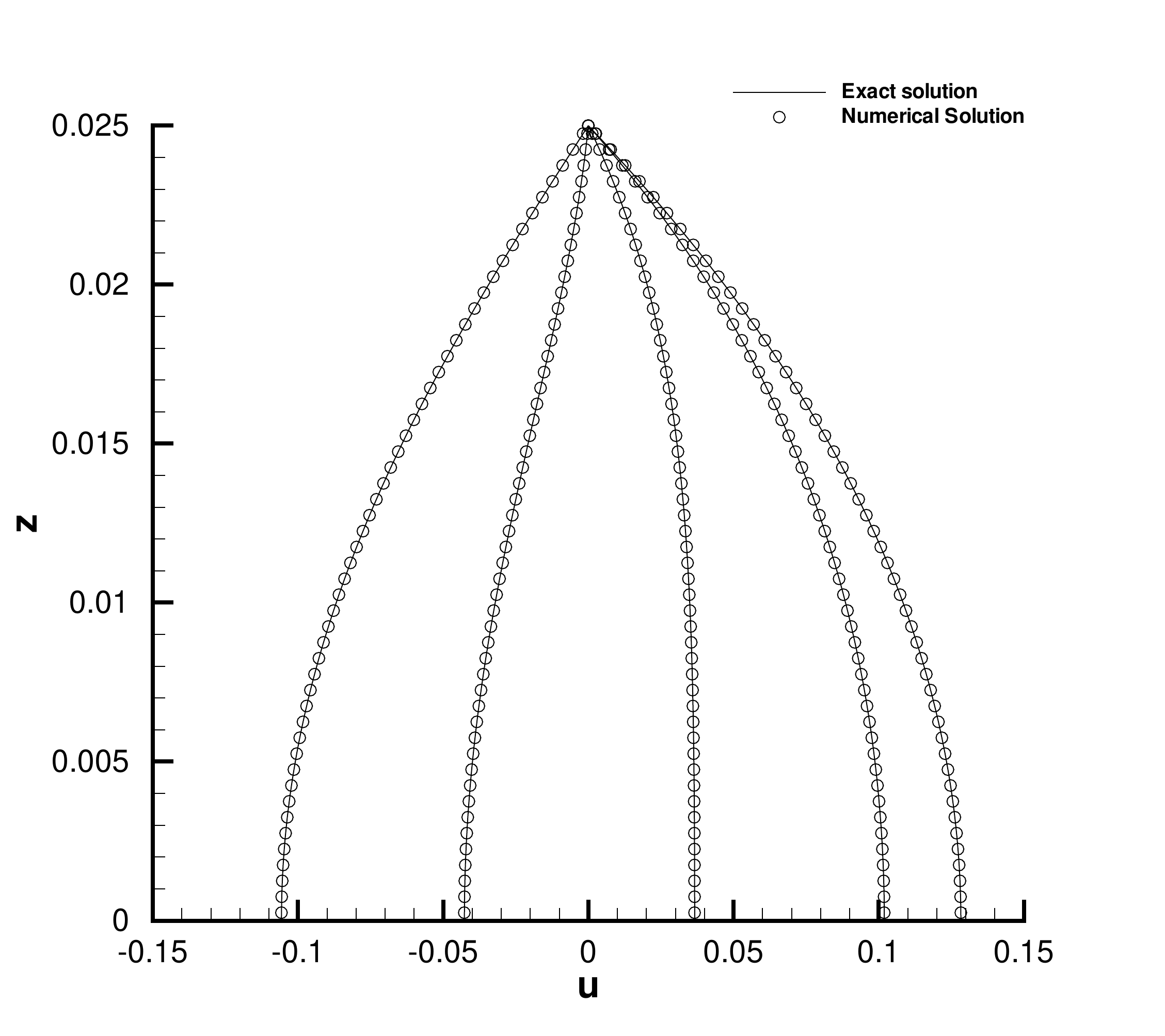}
			\end{subfigure}\;
			\begin{subfigure}[b]{0.45\textwidth}
			\centering
			\includegraphics[width=\textwidth]{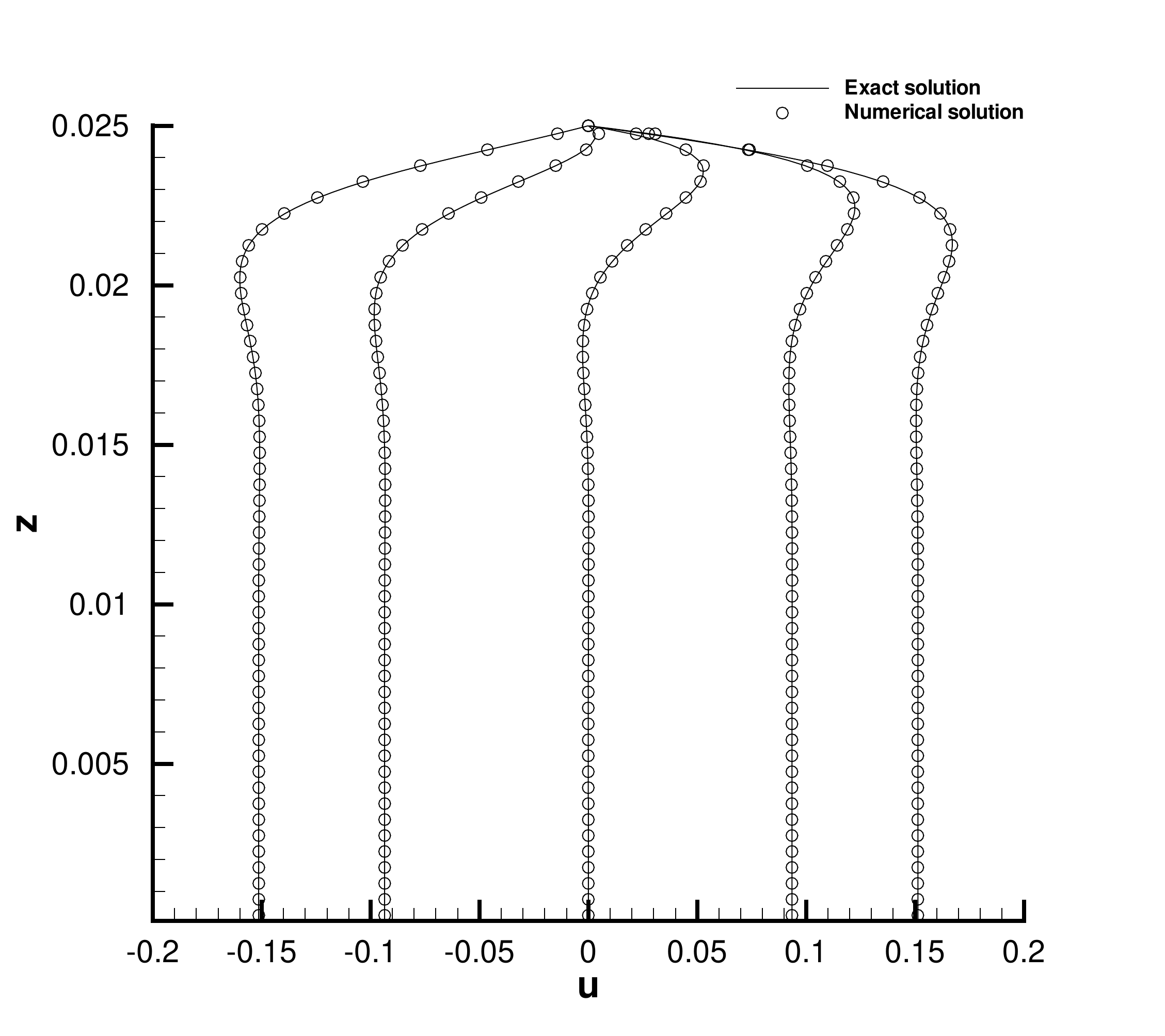}
			\end{subfigure}
\caption{Comparison of the exact solution of Womersley \cite{Womersley:1955} with the numerical results at different times. Left ($R_e = 50$): the different graphs left to right correspond to the times $t=1.7$, $t=1.8$, $t=1.9$, $t=2.0$, and $t=2.1$, respectively. Right ($R_e = 5000$): the different graphs left to right correspond to the times $t=1.8$, $t=1.9$, $t=2.0$, $t=2.1$ and $t=2.2$, respectively.}\label{fig:Womersley}
\end{figure}

The computational domain is discretized with $N_x=100$ segments, $N_z = 50$ rings and $N_\varphi=30$ angular slices. By using $\theta=0.5$ and $\theta'=1$, the simulation is advanced with a time-step size $\Delta t = 0.01$ until a final time $t_e=3$ so that three cycles of oscillation are covered. The computed results are illustrated in Figure \ref{fig:Womersley} next to the analytical solution for both the low and the high Reynolds number regime. An excellent agreement between the numerical and the exact solution is clearly shown. The order of accuracy is confirmed numerically by measuring the error against the exact solution for the high Reynolds number flow ($R_e = 5000$) on a sequence of meshes that are successively refined in both space and time. The angular discretization number is kept constant $N_\varphi=30$, because the solution is still independent of the angular  coordinate. For this test, $\theta = 0.5$ and $\theta'=1$. The convergence results are listed in Table \ref{tab:cap5Unsteady}, confirming that second order of accuracy is achieved in space and time for \emph{unsteady} flow problems. 

\subsection{Oscillating flow in a straight vessel of elliptical cross section}
\label{sec.wom.ell} 

Once the numerical method is shown to be accurate in simulating fluid flow in circular vessels, the case of elliptical cross sections is investigated. In the present section, the oscillatory flow in a straight rigid tube of elliptical cross section is considered. 

Similar to the Womersley problem before, the tube walls are rigid and a sinusoidal pressure gradient is applied at the ends. Analytical studies of this family of pulsatile flow problems were first given in \cite{Khamrui:1957,Verma:1960}, and more recently in \cite{Haslam:1998}, where a more detailed analysis has been provided. The radial and angular elliptical coordinates $\xi$ and $\eta$ are introduced as 
\begin{equation*}
z\cos{\varphi} = d \cosh{\xi} \cos{\eta}, \qquad   z\sin{\varphi} = d \sinh{\xi} \sin{\eta},
\end{equation*}
where $2d$ is the interfocal distance of the ellipse. %
The radial coordinate $\xi$ varies in $[0,\xi_0]$ from the interfocal line to the tube walls, while the angular coordinate $\eta$ varies in $\left[\right.0,2\pi\left.\right)$. 
By neglecting advective and axial viscous terms and by assuming the pressure to be hydrostatic, according to \cite{Haslam:1998}, the axial velocity profile is uniform in the axial $x$ direction and is given by the real part of the expression
\begin{gather}
u(\xi,\eta,t) = \frac{\hat{P}}{\rho} \frac{1}{i\omega} \left[ 1- 2\pi \sum_{n=0}^{\infty} \frac{A_0^{(2n)}}{Ce_{2n}(\xi_0,-q)I_{2n}}Ce_{2n}(\xi,-q)ce_{2n}(\eta,-q) \right] e^{i\omega t}, \label{eq:Khamrui}
\end{gather}
where $Ce_{2n}$ and $ce_{2n}$ are respectively the ordinary and modified Mathieu functions of order $2n$; $A^{(2n)}_{2r}$ are constant coefficients that satisfy a recurrence relation (see \cite{Haslam:1998b} for details); 
\begin{gather*}
I_{2n} = \int_0^{2\pi} ce^2_{2n}(\eta,-q) \,d\eta = 2\pi \left[ A_0^{(2n)}\right]^2 + \pi \sum_{r=2}^{\infty} \left[A_{2r}^{(2n)}\right]^2;
\end{gather*}
$q = \left. i \lambda d^2 \middle/ 4\sigma^2 \right. $, where $\sigma=\sqrt{ \left. 2 \alpha_1^2 \alpha_2^2 \middle/(\alpha_1^2+\alpha_2^2) \right.}$ is the characteristic length parameter of an ellipse, by defining the major and minor axis of the elliptical cross section area of the tube as $\alpha_1 = R(x,0)$ and $\alpha_2 = R(x,\pi/2)$ respectively,  and $\lambda = \omega \sigma^2 / \nu$ is the frequency parameter that reduces to the square of the Womersley number $\alpha_W$ for circular cross sections. For the present test, the chosen parameters are $L = 1$, $R=0.025$, $\hat{P}=1000$, $\rho = 1000$, $\nu = 10^{-4}$,  $\beta = 10^{12}$ so that the tube wall is sufficiently rigid, with aspect radio $\alpha_1/\alpha_2 = \tan \xi_0$ fixed equal to $0.3$. 

\begin{figure} 
\centering 
			\begin{subfigure}[b]{\textwidth}
			\centering
			\includegraphics[width=0.45\textwidth]{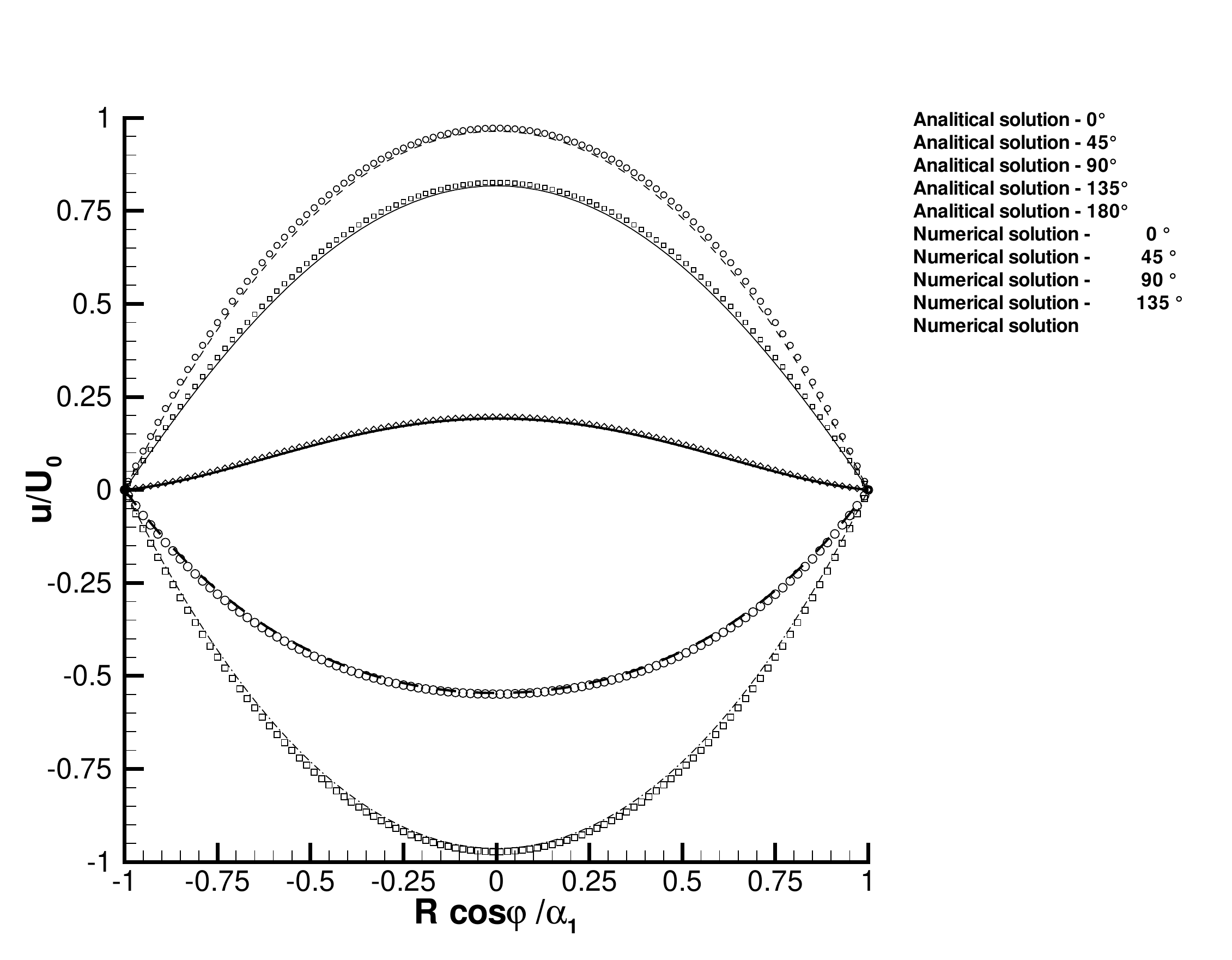}\;
			\includegraphics[width=0.45\textwidth]{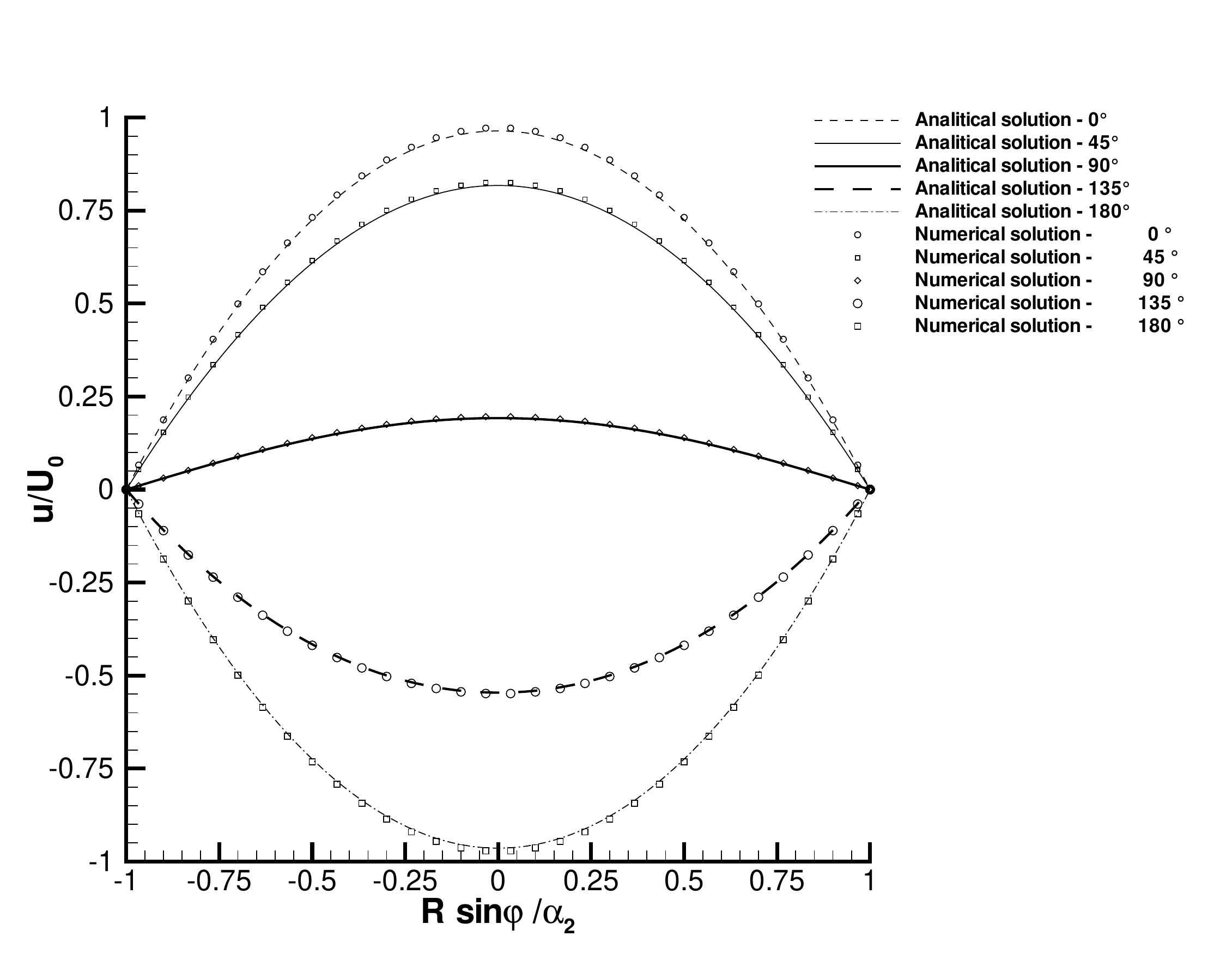}
			\caption{$\lambda=1$.}
			 \label{fig:EllUnsteadyLambda1}
			\end{subfigure}\\
			\begin{subfigure}[b]{\textwidth}
			\centering
			\includegraphics[width=0.45\textwidth]{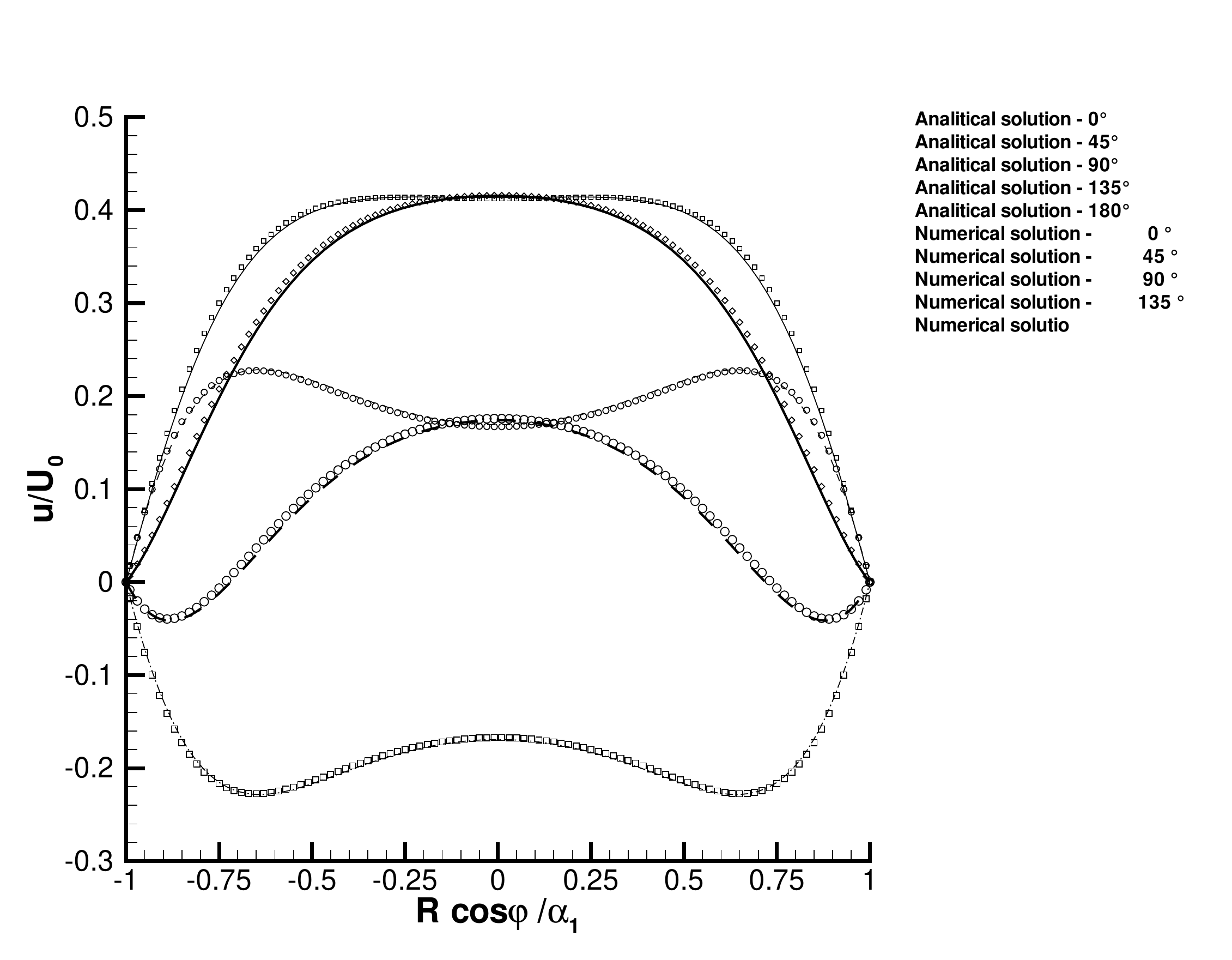}\;
			\includegraphics[width=0.45\textwidth]{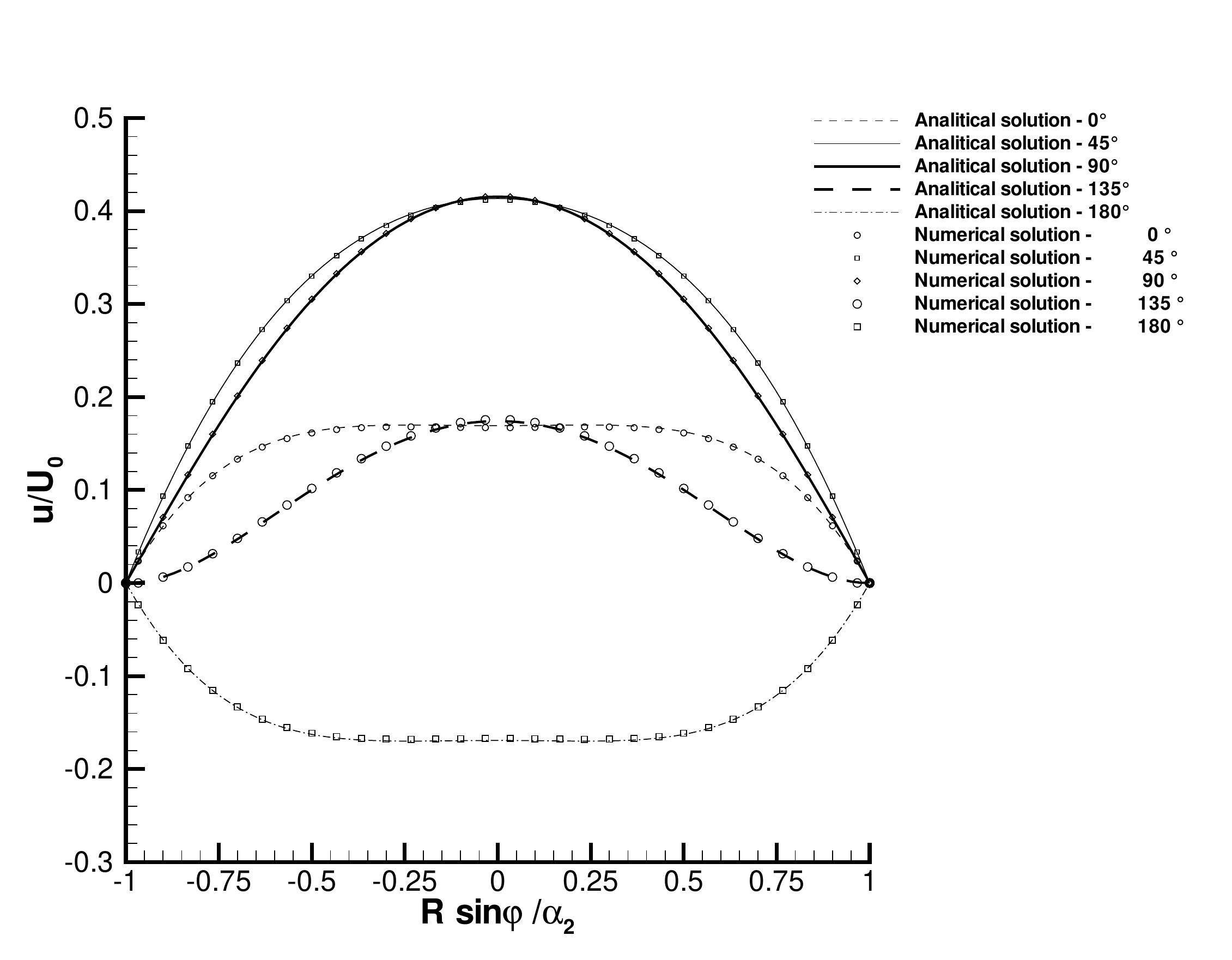}
			\caption{$\lambda=10$.}
			 \label{fig:EllUnsteadyLambda10}
			\end{subfigure}\\
			\begin{subfigure}[b]{\textwidth}
			\centering
			\includegraphics[width=0.45\textwidth]{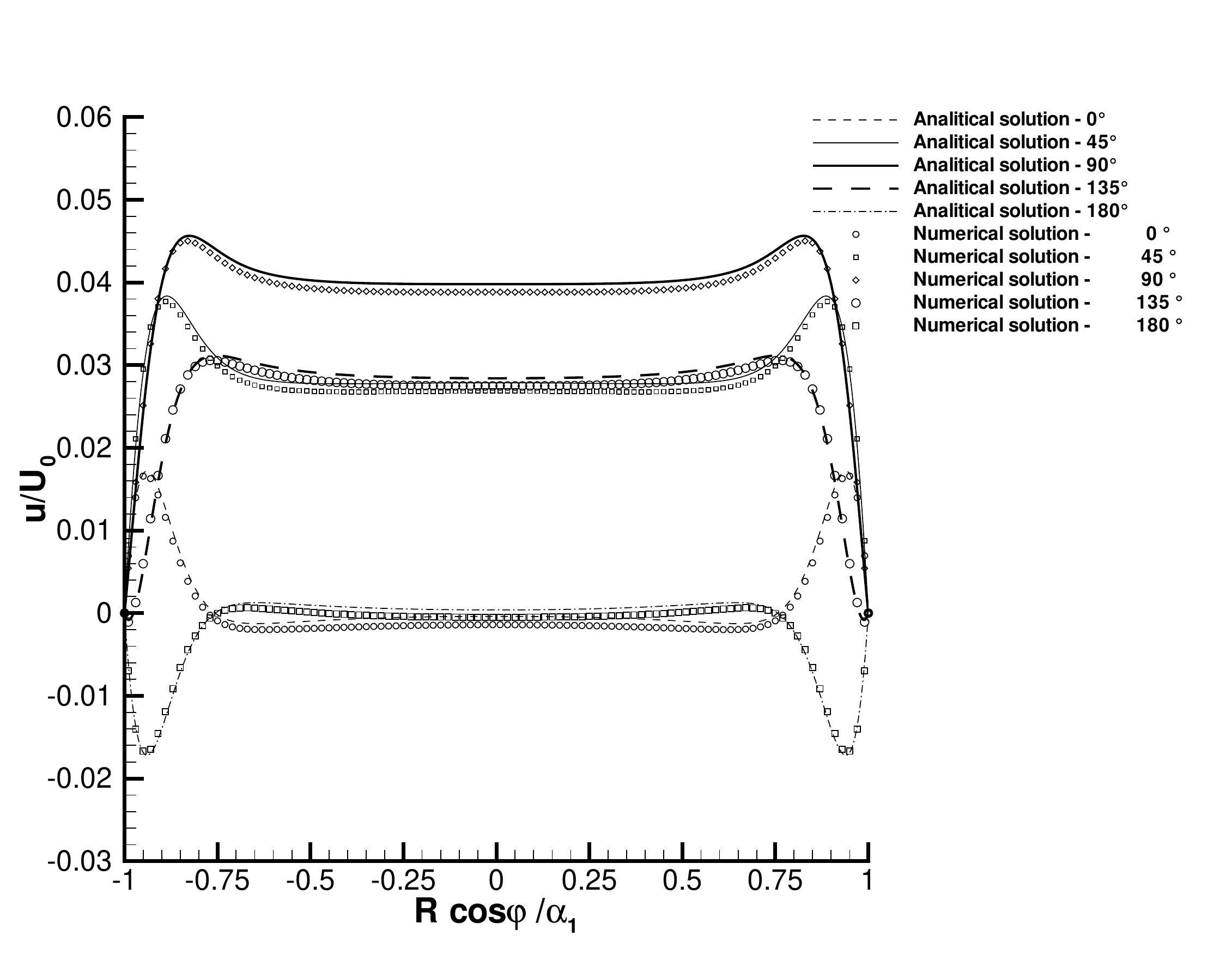}\;
			\includegraphics[width=0.45\textwidth]{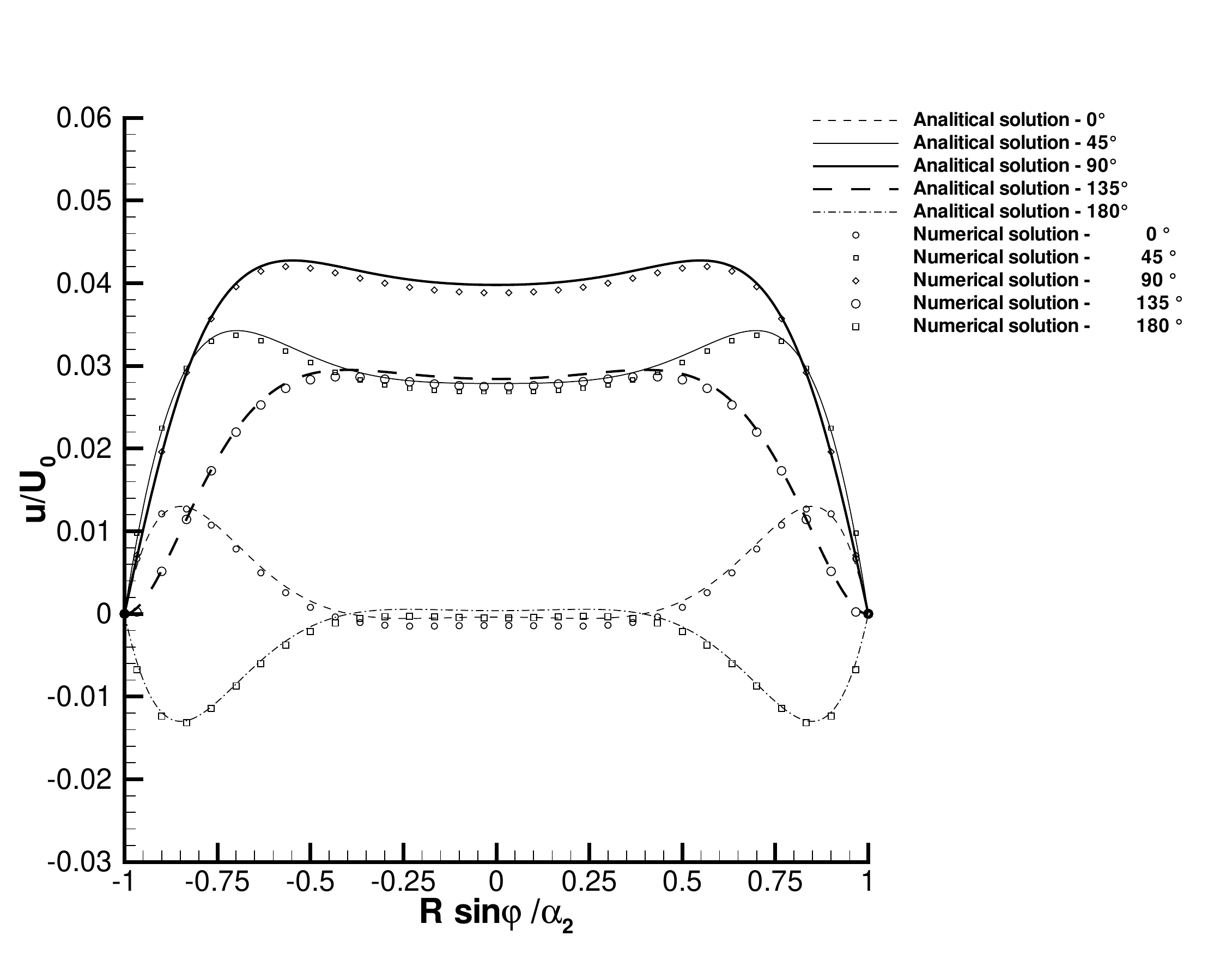}
			\caption{$\lambda=100$.}
			 \label{fig:EllUnsteadyLambda100}
			\end{subfigure}
\caption{Comparison of the exact solution \cite{Haslam:1998} with the numerical results at different times during the first half of a period. The axial velocity field is interpolated along the major axis (left) and the minor axis (right).}\label{fig:EllUnsteady}
\end{figure}

In order to verify the accuracy of the numerical results, three different values of the frequency parameter are considered. First, a low frequency number $\lambda=1$, corresponding to a low Reynolds number flow, is chosen. In this  regime the entire cross section is dominated by the viscous effect, and the velocity profile resembles the elliptic paraboloid profile, typical of steady flows (Figure \ref{fig:EllUnsteadyLambda1}). In the second case, by choosing $\lambda=100$, the resulting velocity profiles show an essentially inviscid (flat) core and a sharp boundary layer close to the walls in which viscous effect dominates (Figure \ref{fig:EllUnsteadyLambda100}). Then, an intermediate frequency number $\lambda=10$ is further chosen in order to emphasize the development of the double peaked velocity profile, peculiar of elliptical sections, that arise smoothly by increasing the frequency parameter (Figure \ref{fig:EllUnsteadyLambda10}).

The initial conditions are taken to be $u(x,z,\varphi,0)=0$, $R(x,\varphi,0)=R_0$ with $p(x,\varphi,0) = -x\hat{P}/\rho$.
Then, the time dependent boundary conditions are specified at the two ends of the tube by applying, for times $t > 0$, the same oscillating pressure gradient of the Womersley problem. %

The computational domain is discretized with $N_x=50$ segments, $N_z = 50$ rings with reference to the cylinder of radius $\alpha_1$ that contains the whole elliptic tube, and $N_\varphi=160$ angular slices in order to build a sufficiently fine mesh to resolve all details of the flow field well. By using $\theta=0.5$ and $\theta'=0.5$, the simulation is advanced with a time-step size $ \omega \Delta t = 5^\circ$ until three cycles of oscillation are  covered. The computed results are illustrated in Figure \ref{fig:EllUnsteady} next to the analytical solution for the three chosen regimes. A very good agreement between the numerical and the exact solution is shown, similar to the circular Womersley problem. 
The test problems described in this section and the previous one show how the present method can accurately reproduce complex hydrostatic flows. 

  %

\subsection{{$\mathit{90^{\circ}}$} bended tube}
The geometry of the grid of a uniformly curved section is defined with the reference to the \emph{curvature radius} $R_c$, and the \emph{curvature angle} $\phi_c$. Toroidal coordinates are introduced in the equation of motion without loss of stability. 
Notice that, with regard to this, toroidal advective and viscous terms are correctly computed because of the Eulerian-Lagrangian method.

With the aim of validating the non-hydrostatic approach it is important to test the model against a physical problem in which non-hydrostatic components play a leading role in the fluid dynamics. In the following, the results of the simulation of steady and pulsatile flow in an uniformly curved rigid tube are presented trying to reproduce the \emph{experimental data} and \emph{numerical results} given by \cite{Timite:2010}. Experimental measurements of axial velocity components are performed by a Laser Doppler Velocimeter.

The geometry consists of a uniformly curved tube of circular cross section with rigid walls ($\beta = 10^{12} \text{Pa/m}$), angle of curvature $\phi_c=90^\circ$, radius of curvature $R_c=0.22$m  and cross sectional diameter $D_0 = 2 R_0 = 0.04$m. For the simulations, the curvature of the numerical model is characterized by $\phi_c = 1.05\times 90^{\circ}$ in order to avoid numerical boundary effects. The properties of the fluid are approximately the ones of water at $300$K: $\nu = 10^{-6}\text{m}^2/\text{s}$, $\rho = 1000 \text{Kg}/\text{m}^3$. We use the same Laplace law \eqref{eq:cap5laplace} as for 
the straight tube, which is a valid approximation only for $R_0/R_c \ll 1$. 

\subparagraph{Steady flow.}
Starting from rest, the numerical flow is driven by a fixed entry condition of Poiseuille type and a fixed hydrostatic pressure at the exit. Every simulation is described by a fixed Reynolds number $R_e$ defined as
\begin{gather*}
R_e = 2 R_0 \frac{Q}{\nu \pi R_0^2} = 2 R_0 \frac{ U_0}{\nu}  %
\end{gather*}
where $Q$ is the mean flux and $U_0=Q/(\pi R_0^2)$ is the respective mean axial velocity. The respective value of the Dean number, based on the hydraulic diameter, are reported as
$D = 4  R_e\sqrt{2 R_0/R_c} $.

Nonlinear advective terms cannot be neglected because to the continuous change of direction of streamlines in curved sections. Then, an Eulerian-Lagrangian approach is consequently used. Numerical results show clearly what happens  if non-hydrostatic contributions are excluded, leading to unphysical solutions. 
The $\theta$-method is run with $\theta=1$ and $\theta'=1$.
For the present test, the discretization numbers are $N_x=63$, $N_z=40$ and $N_\varphi=48$. The simulation is advanced for $N_t = 300$ time steps with a time-step size $\Delta t = t_e / N_t$ until a final time $t_e=60s$.

\begin{figure}[!htbp]
\centering 
			\begin{subfigure}{0.48\textwidth}
			\centering
			\includegraphics[width=\textwidth]{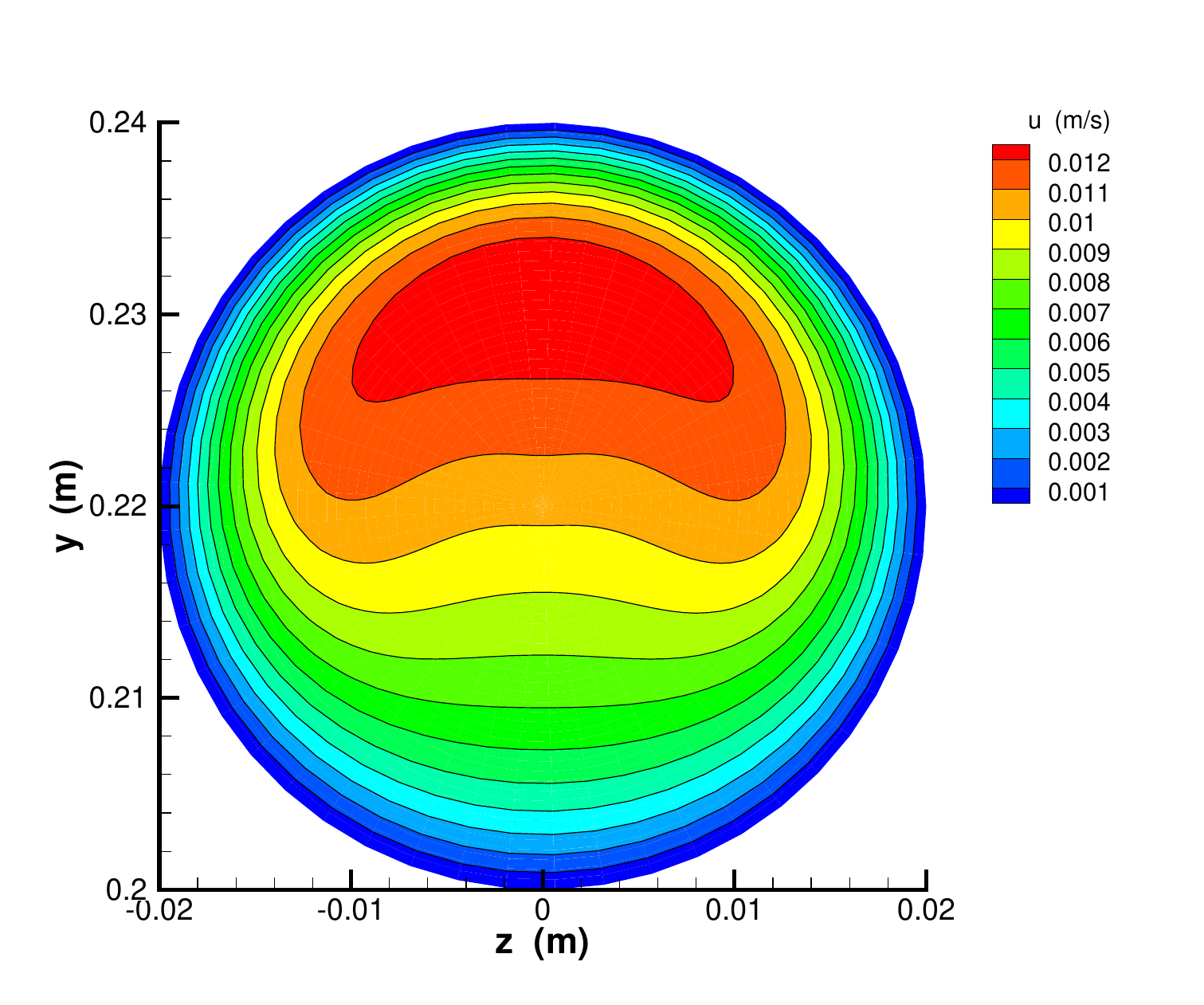}
			\end{subfigure} \;
			\begin{subfigure}{0.48\textwidth}
			\centering
			\includegraphics[width=\textwidth]{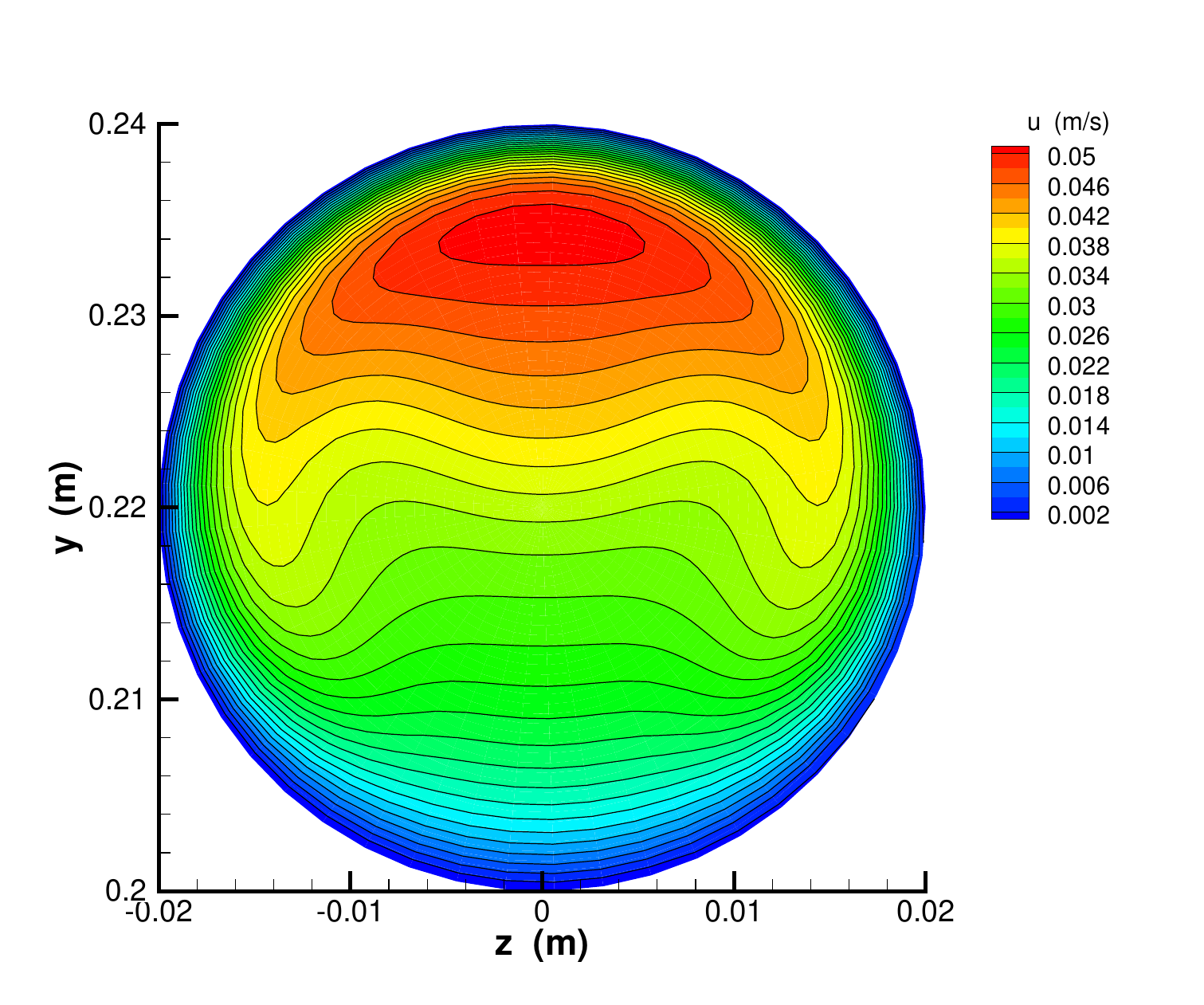}
			\end{subfigure}\\
			\begin{subfigure}{0.48\textwidth}
			\centering
			\includegraphics[width=\textwidth]{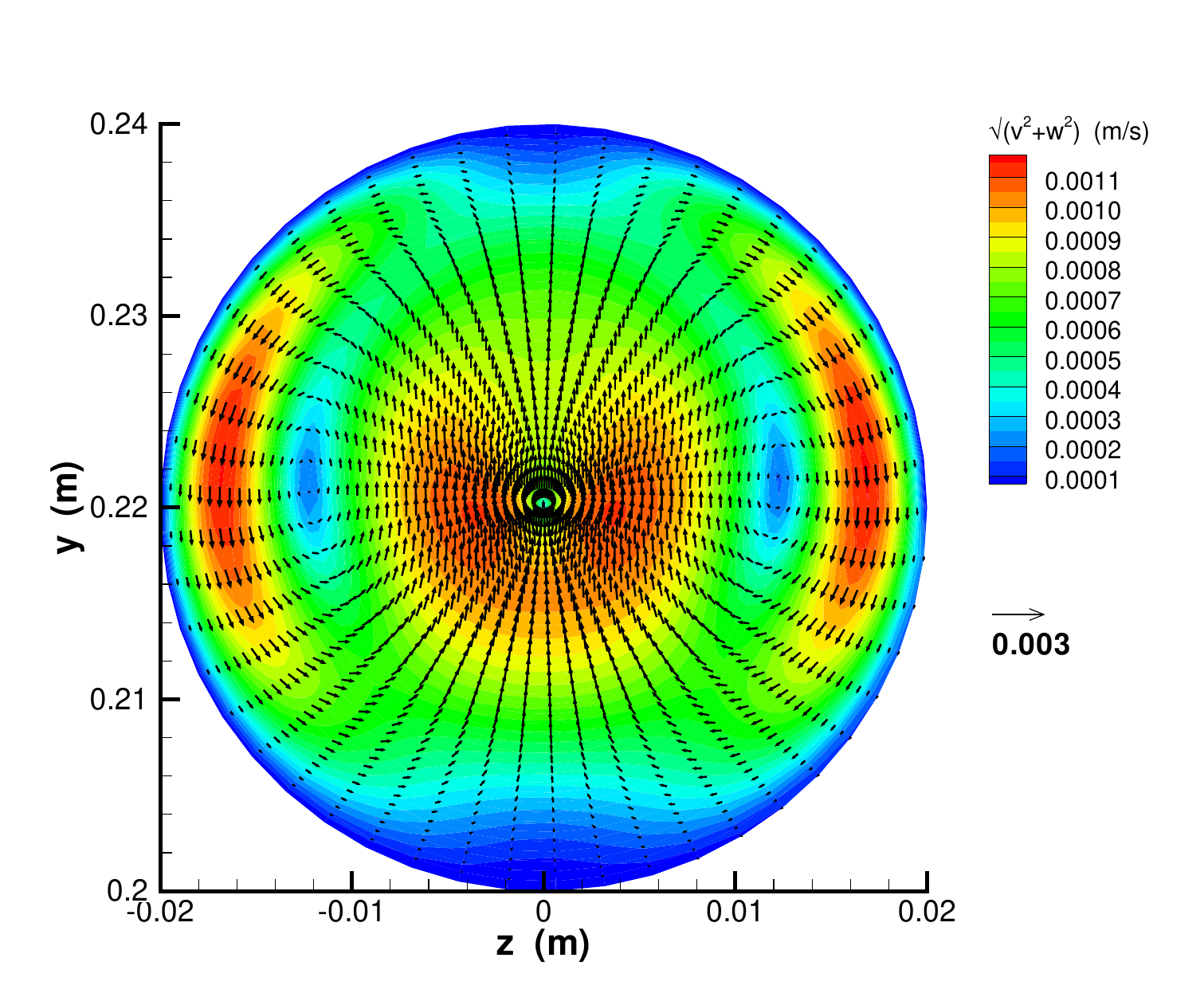}
			\end{subfigure} \;
			\begin{subfigure}{0.48\textwidth}
			\centering
			\includegraphics[width=\textwidth]{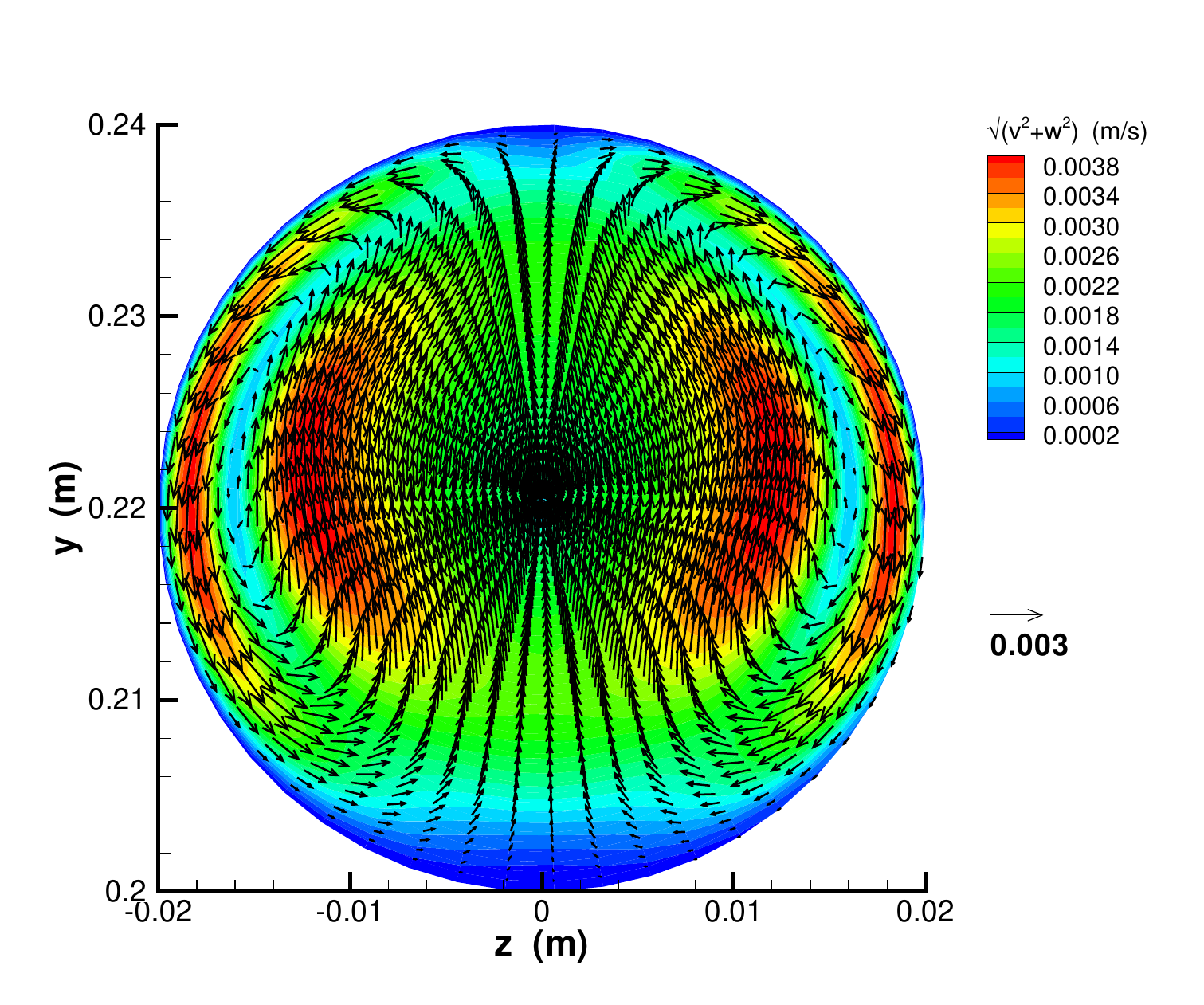}
			\end{subfigure}\\
			\begin{subfigure}{0.48\textwidth}
			\centering
			\includegraphics[width=\textwidth]{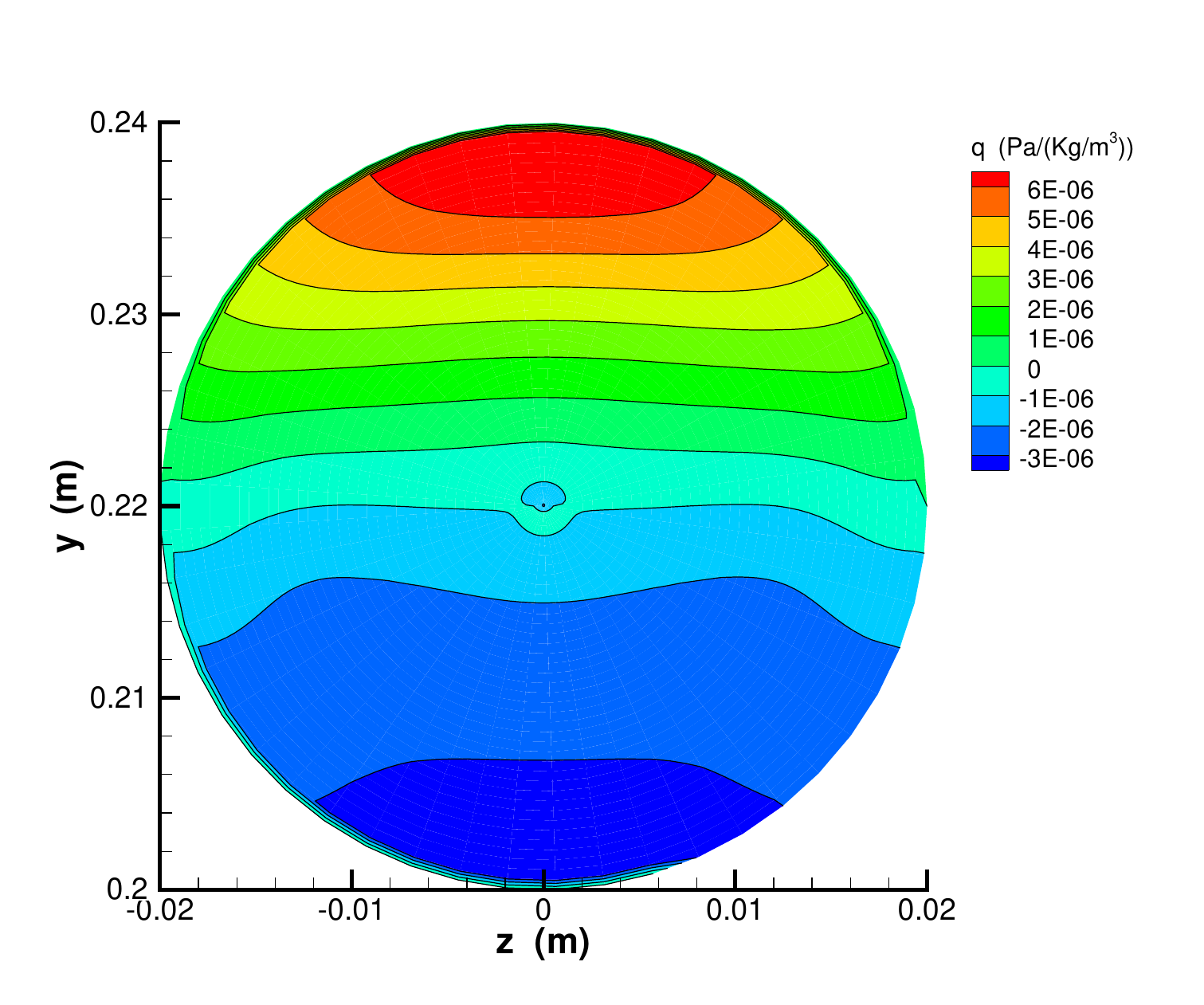}
			\caption{$R_e = 300$ ($D \sim 512$).}
			 \label{fig:cap5CrossSecRe300new}
			\end{subfigure} \;
			\begin{subfigure}{0.48\textwidth}
			\centering
			\includegraphics[width=\textwidth]{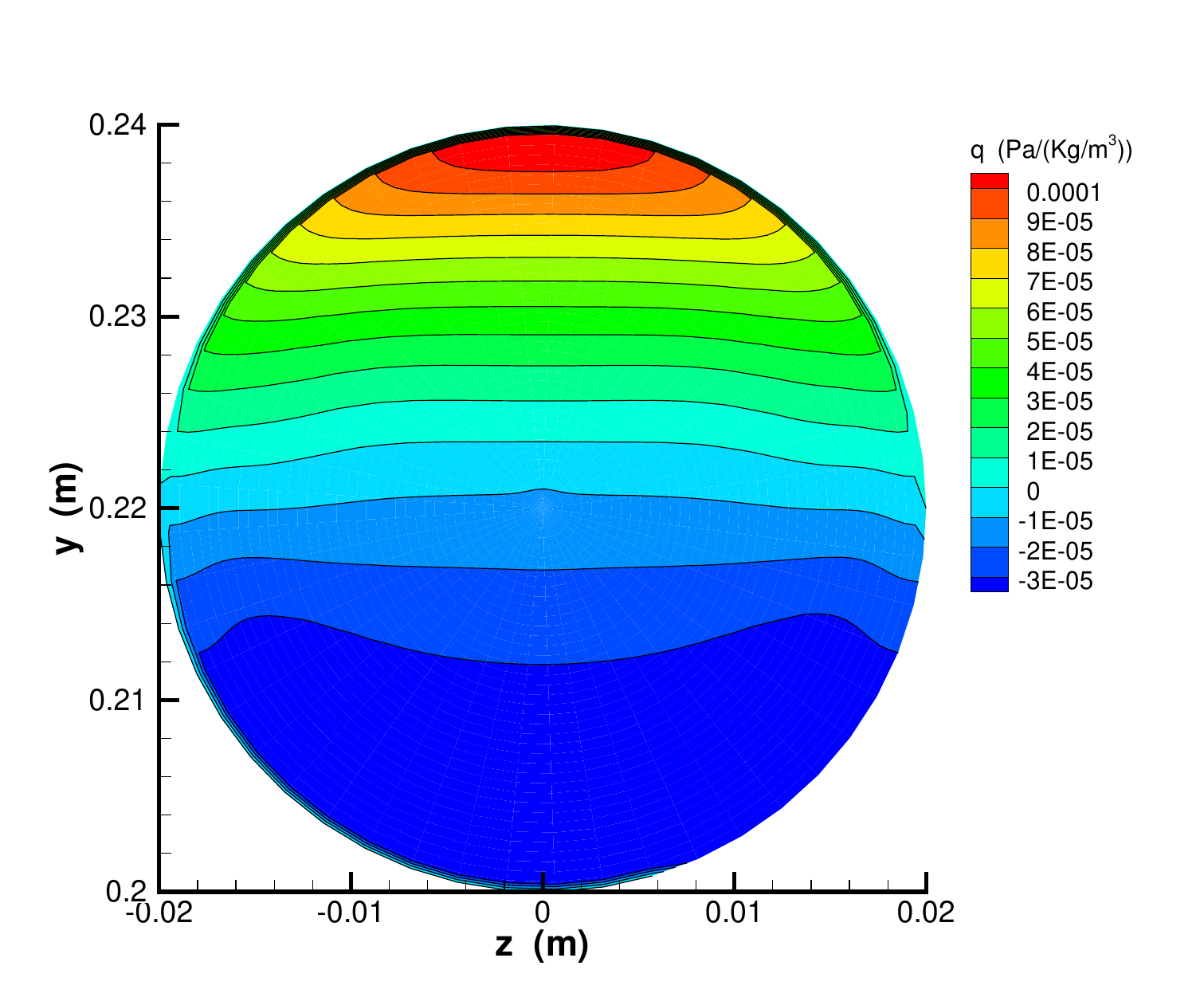}
			\caption{$R_e = 1200$ ($D \sim 1706$).}
			\label{fig:cap5CrossSecRe1200new}
			\end{subfigure}\\
\caption{Physical quantities interpolated along the $90^\circ$ cross section, for $R_e = 300$ and $R_e = 1200$. From the top to the bottom: axial velocity contours; velocity vector field, tangential to the plane; non-hydrostatic correction $\mathit{q}$.}\label{fig:cap5CrSec3001200new}
\end{figure}

Figure \ref{fig:cap5CrSec3001200new} shows the velocity and pressure fields interpolated along the exit cross sections at different Reynolds numbers, from $R_e=300$ to $R_e =1200$. Curvature gives rise to a centripetal pressure gradient that, combined with the centrifugal effect, generate the peculiar cross sectional circulation of curved flows: secondary motions are well established as two pairs of symmetric, counter-rotating vortices. Results confirm that curvature effects becomes stronger with the increasing of the Reynolds number: the centrifugal effect shifts the peak of the axial velocity components to the outer side of the curve; the centripetal pressure gradient and the intensity of the secondary motion increase.

\begin{figure}[!htbp]
\centering {
			\begin{subfigure}{0.4\textwidth}
			\centering
			\includegraphics[width=\textwidth]{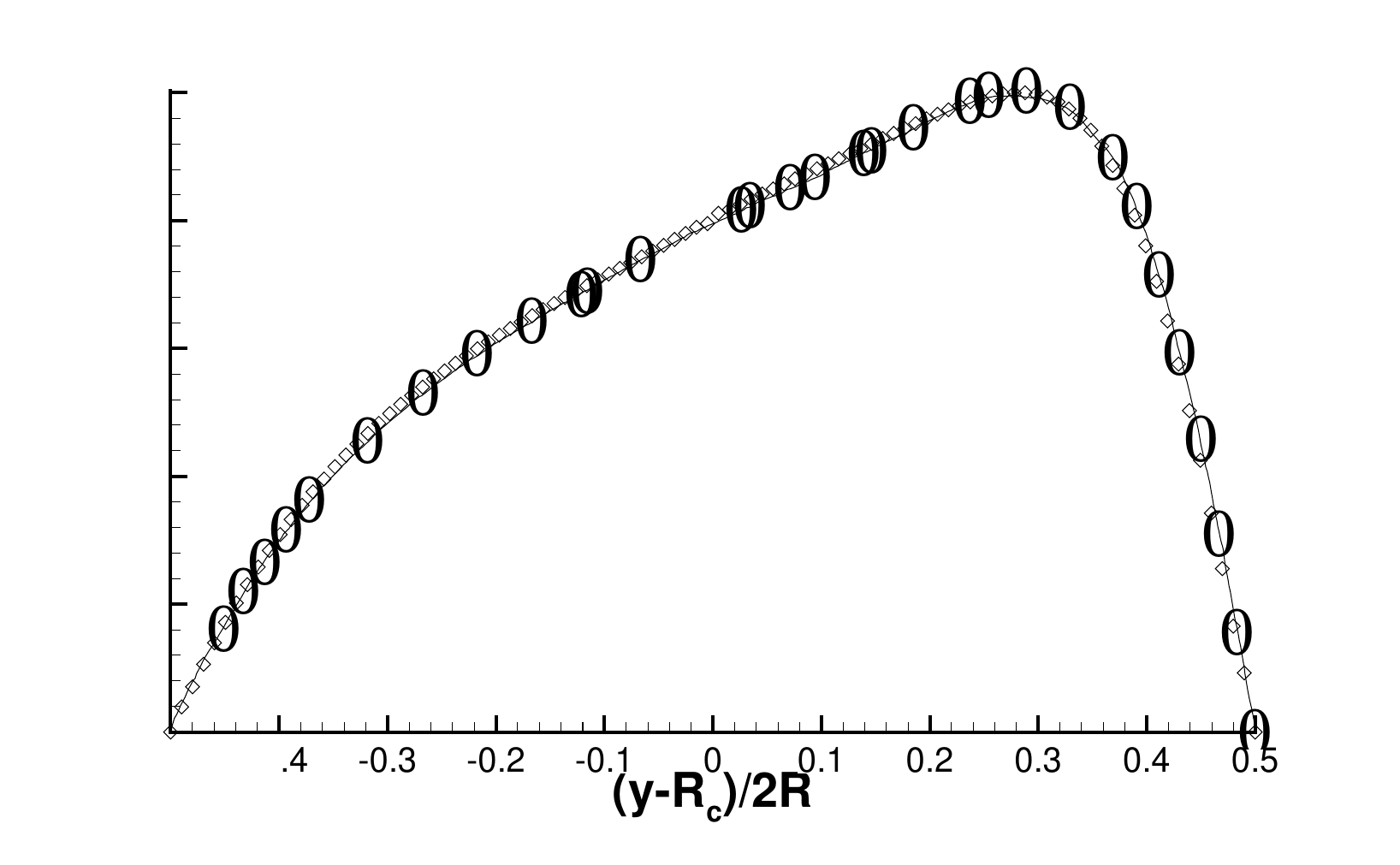} \nonumber 
			\end{subfigure} \hspace{5mm}
			\begin{subfigure}{0.4\textwidth}
			\centering
			\includegraphics[width=\textwidth]{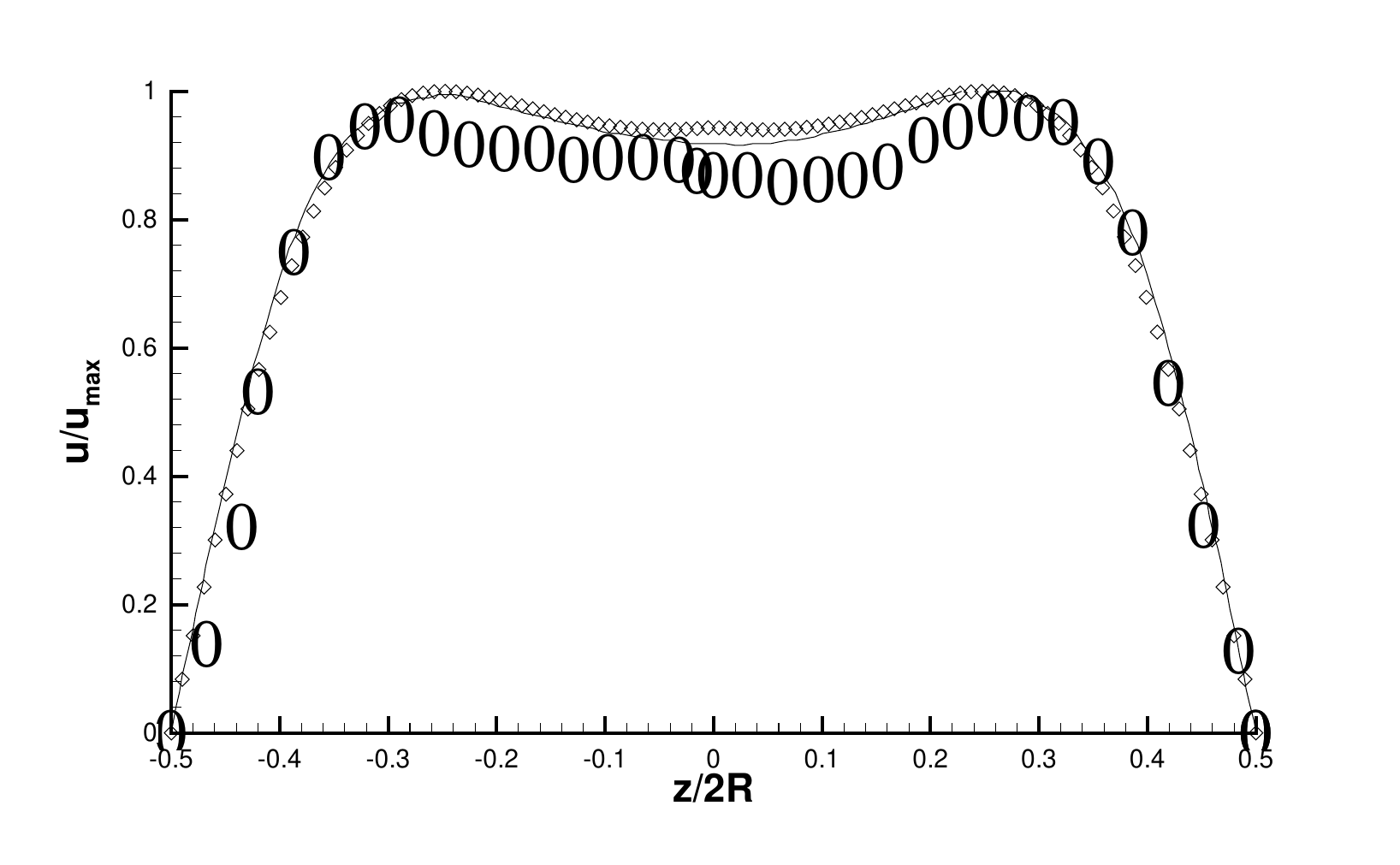} \nonumber
			\end{subfigure}
\subcaption{$R_e = 300$.} \label{fig:cap5Pro300new} 
			\begin{subfigure}{0.4\textwidth}
			\centering
			\includegraphics[width=\textwidth]{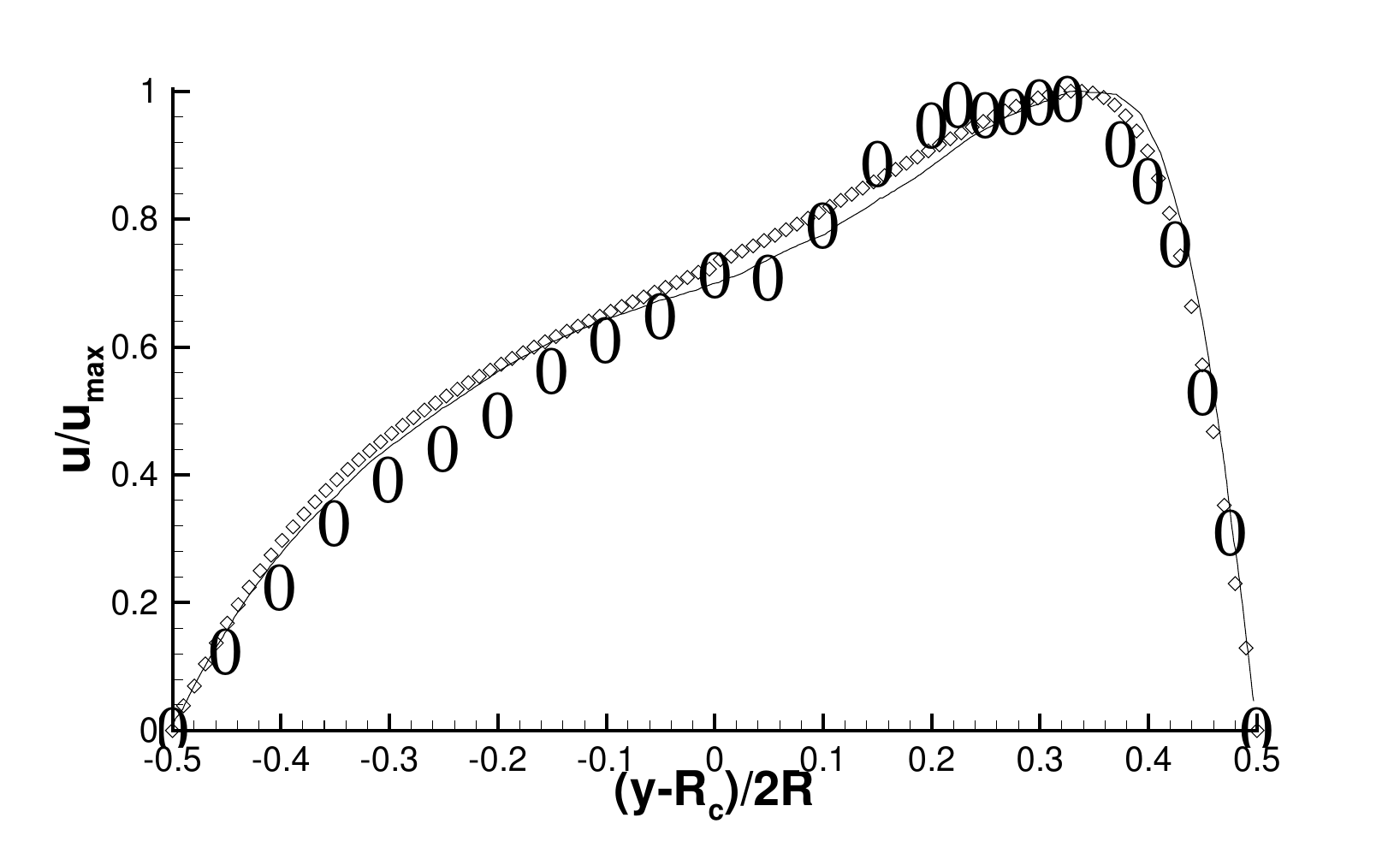} \nonumber  
			\end{subfigure} \hspace{5mm}
			\begin{subfigure}{0.4\textwidth}
			\centering
			\includegraphics[width=\textwidth]{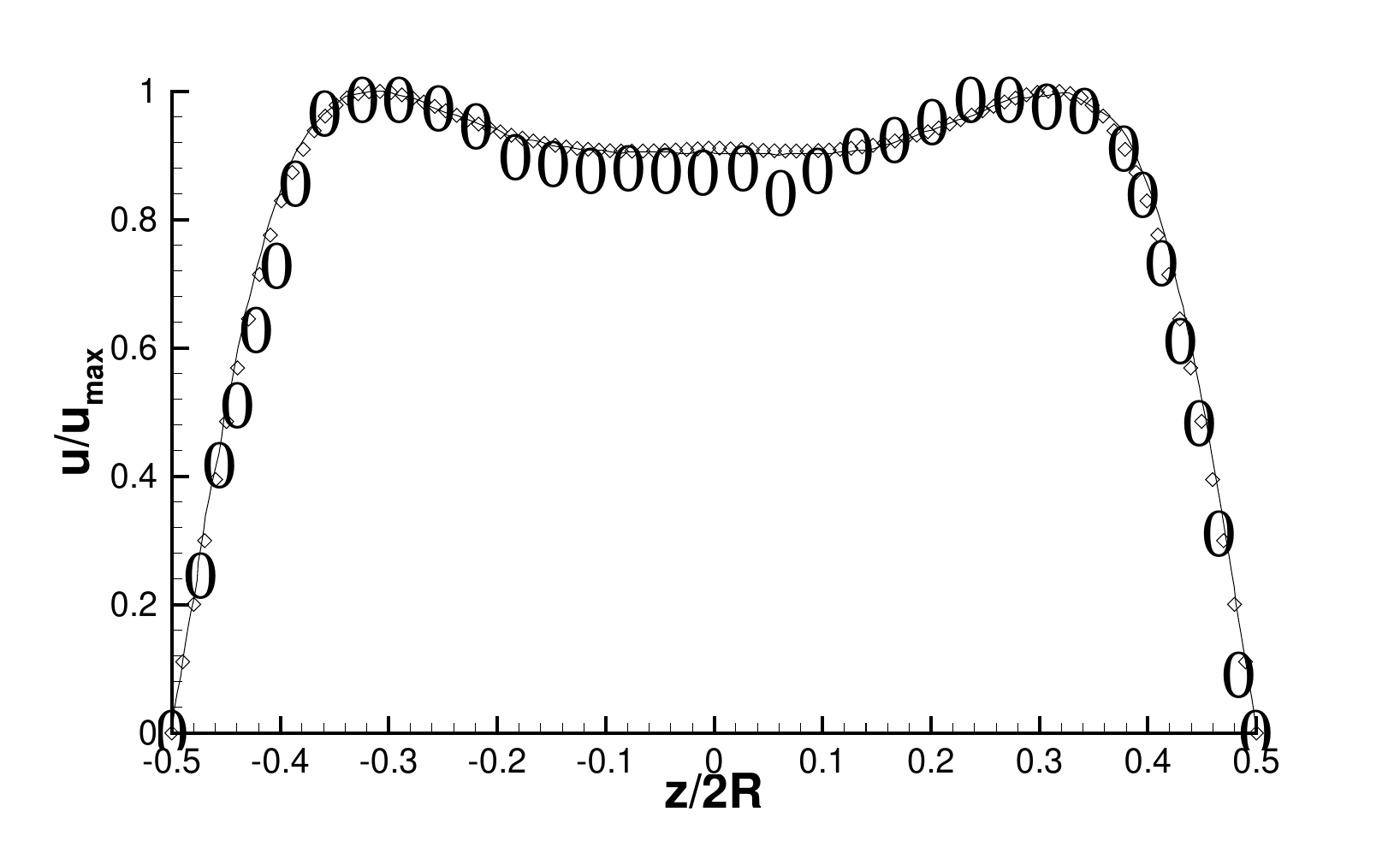} \nonumber
			\end{subfigure}
\subcaption{$R_e = 600$.} \label{fig:cap5Pro600new} 
			\begin{subfigure}{0.4\textwidth}
			\centering
			\includegraphics[width=\textwidth]{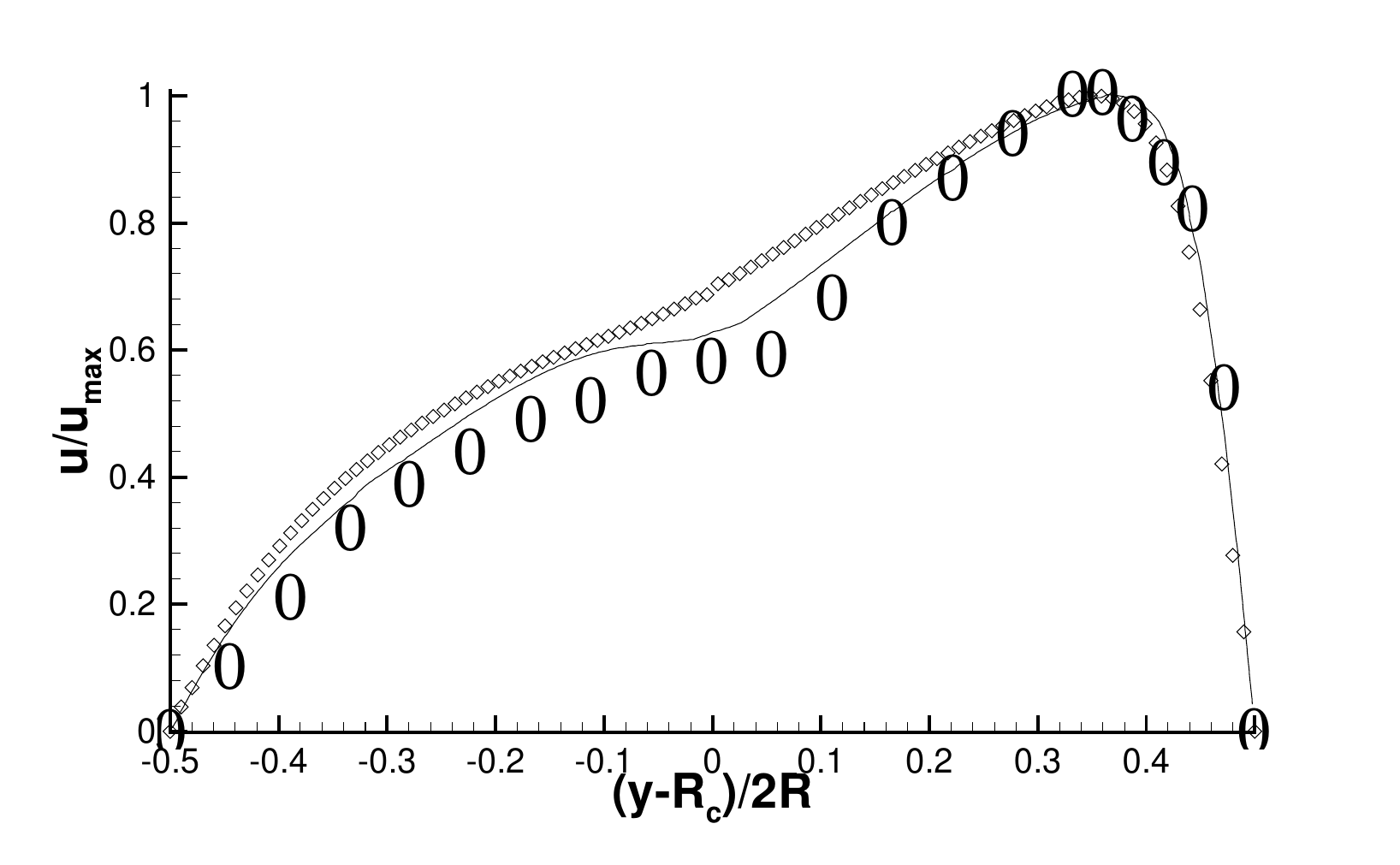} \nonumber 
			\end{subfigure} \hspace{5mm}
			\begin{subfigure}{0.4\textwidth}
			\centering
			\includegraphics[width=\textwidth]{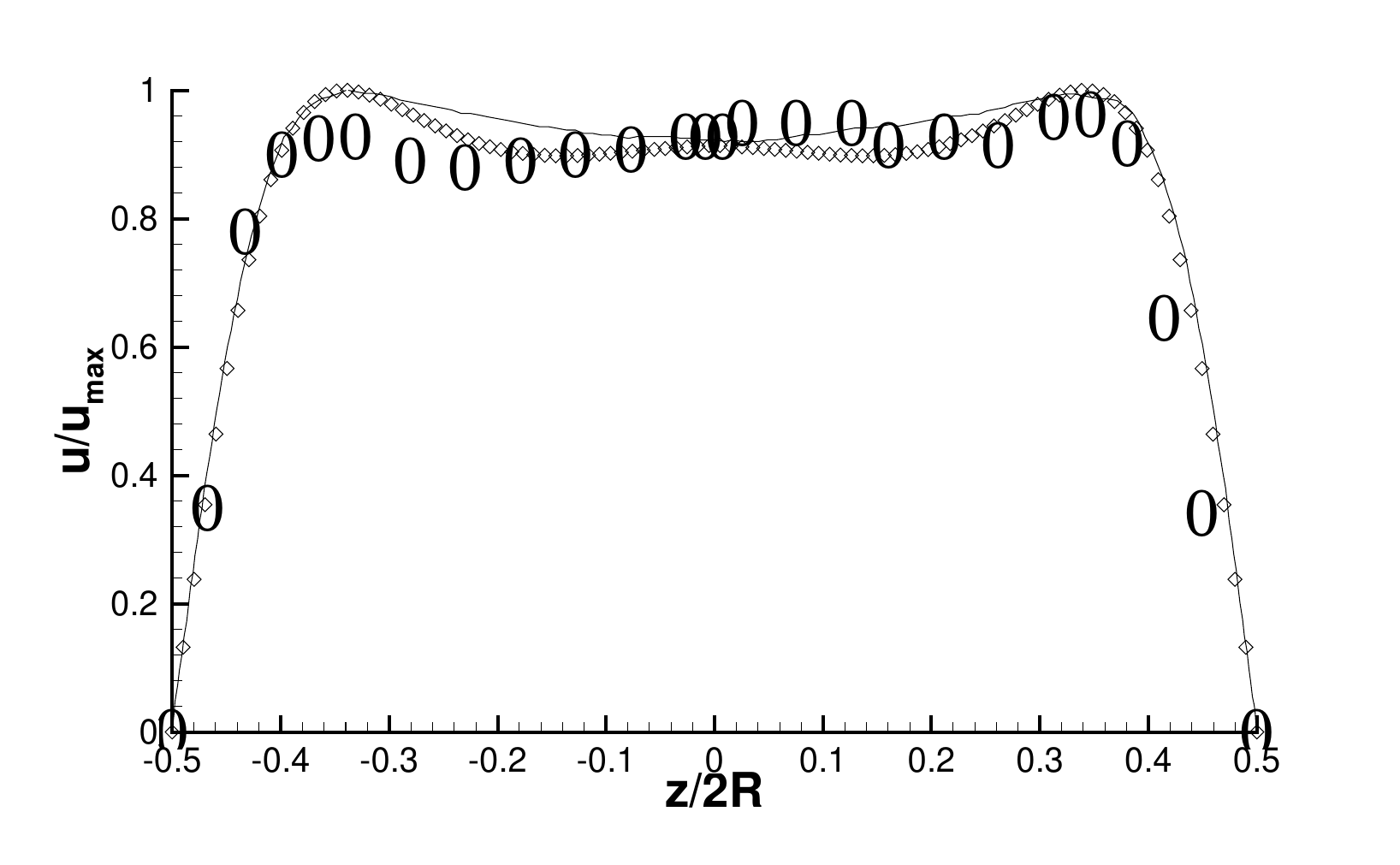} \nonumber
			\end{subfigure}
\subcaption{$R_e = 900$.} \label{fig:cap5Pro900new} 
			\begin{subfigure}{0.4\textwidth}
			\centering
			\includegraphics[width=\textwidth]{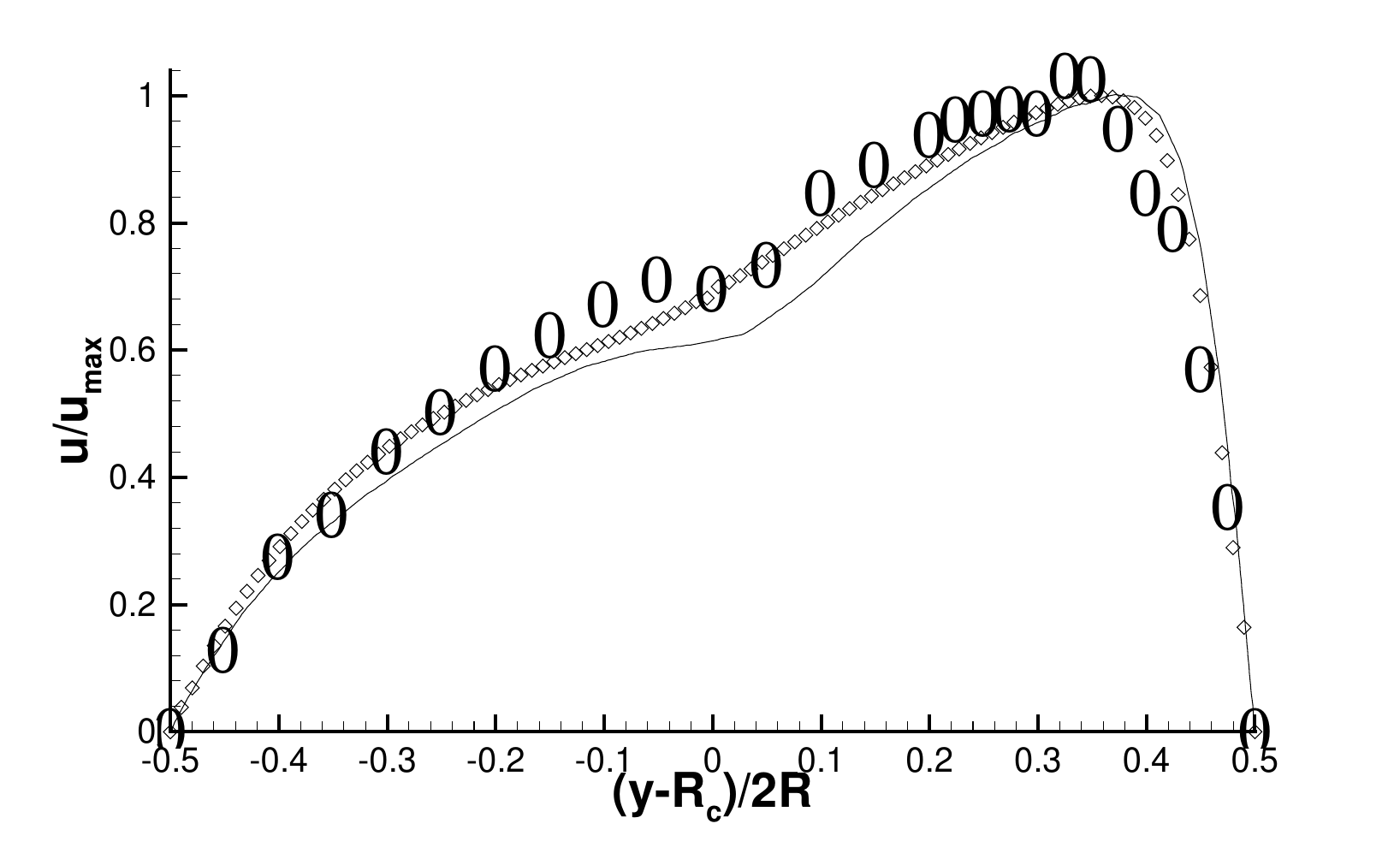} \nonumber 
			\end{subfigure} \hspace{5mm}
			\begin{subfigure}{0.4\textwidth}
			\centering
			\includegraphics[width=\textwidth]{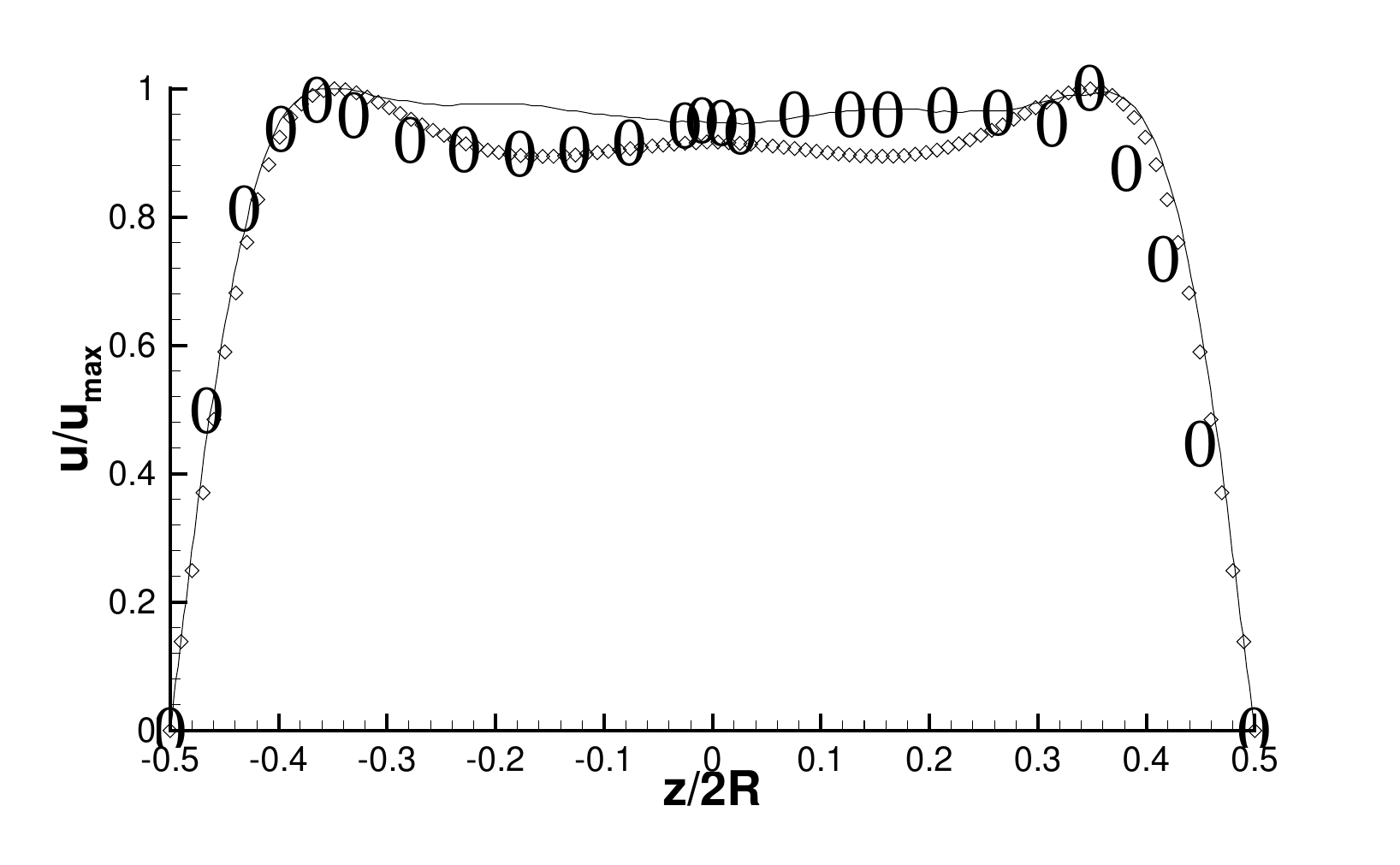} \nonumber
			\end{subfigure}
\subcaption{$R_e = 1000$.} \label{fig:cap5Pro1000new} 
			\begin{subfigure}{0.4\textwidth}
			\centering
			\includegraphics[width=\textwidth]{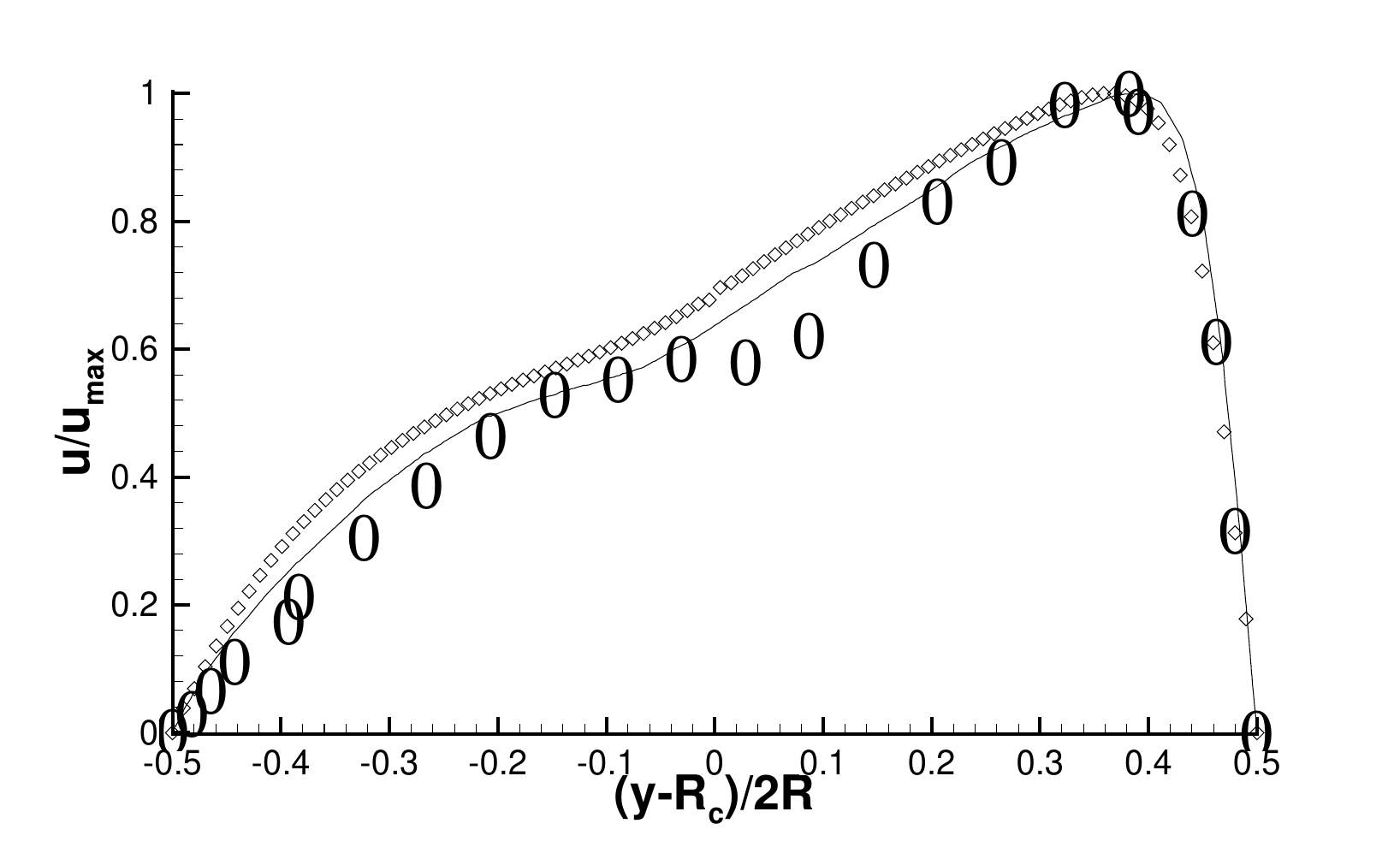} \nonumber 
			\end{subfigure} \hspace{5mm}
			\begin{subfigure}{0.4\textwidth}
			\centering
			\includegraphics[width=\textwidth]{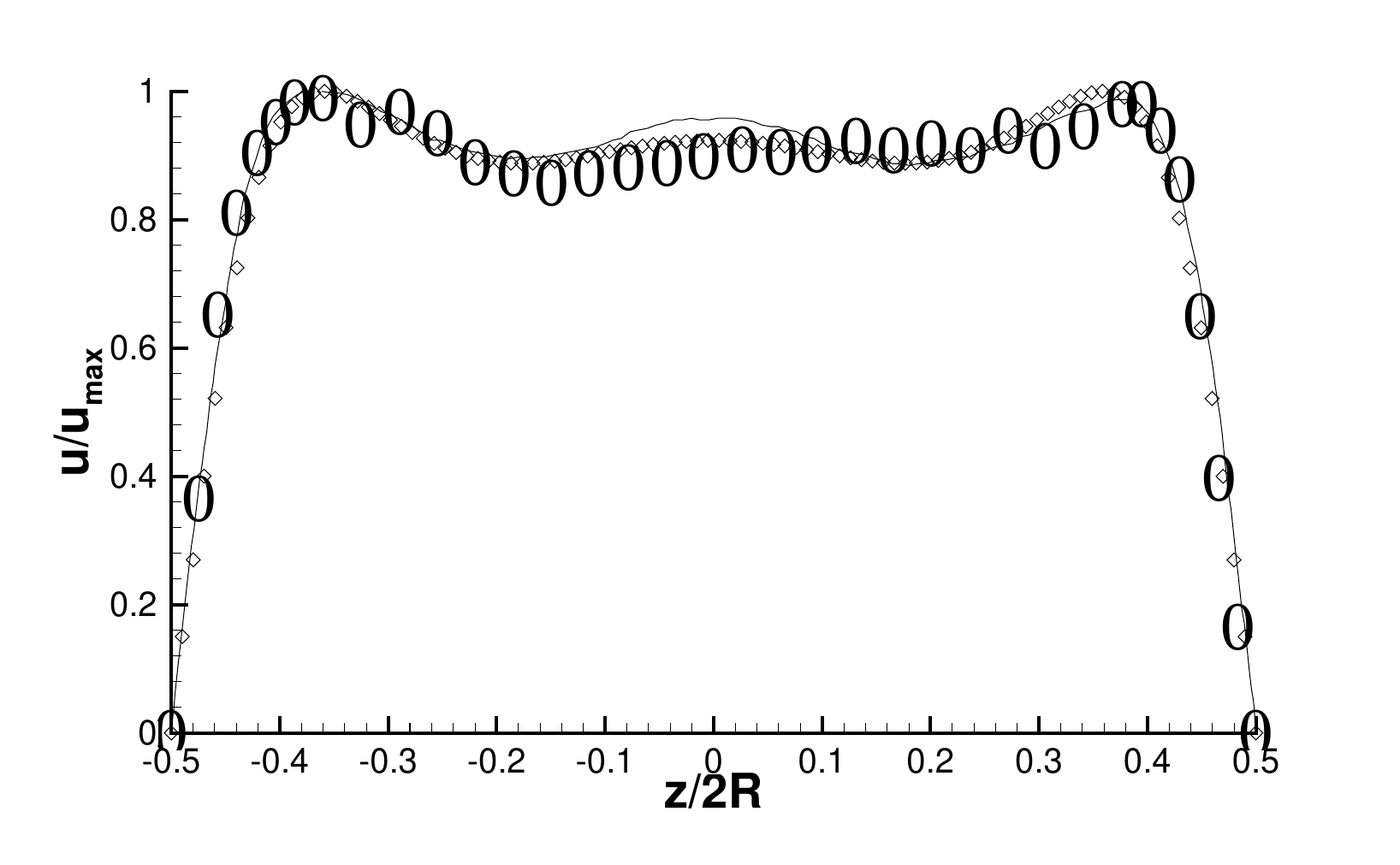} \nonumber
			\end{subfigure}
\subcaption{$R_e = 1200$.} \label{fig:cap5Pro1200new} }
\caption{Axial velocity profiles interpolated along the diameters of the $90^\circ$ cross section in the $\mathit{y}$- and $\mathit{z}$ direction, for different Reynolds numbers. $\diamond$ numerical results; $-$ numerical results of \cite{Timite:2010}; $\mathbf{0}$ experimental data of \cite{Timite:2010}.}
\end{figure}

Axial velocity components of the $90^\circ$ cross section have been interpolated along the orthogonal diameters for the purpose of comparing the digitized data of literature. Figures \ref{fig:cap5Pro300new}-\ref{fig:cap5Pro1200new} show our numerical results next to the numerical and experimental data of \cite{Timite:2010}. The numerical results obtained with the present semi-implicit non-hydrostatic approach are very similar to the numerical and experimental results presented in \cite{Timite:2010}.

\subparagraph{Pulsatile flow.}
The flow is now driven by a pulsatile entry condition, consisting in a steady Poiseuille-type component plus an oscillating Womersley-type component, and a fixed hydrostatic pressure at the exit. Each simulation is described 
by  the \emph{steady} Reynolds number $R_e$ (based on the axial velocity of the steady component) 
and by the Womersley number $\alpha_W = R \sqrt{\omega / \nu}$. The oscillatory component is defined as the real part of the expression 
\begin{gather*}
u_W(x,z,t) = \Re\frac{\hat{P}}{\rho} \frac{1}{i\omega} \left[ 1- \frac{\mathpzc{J}_0\left(\alpha_W y i^\frac{3}{2}\right)}{\mathpzc{J}_0 \left( \alpha_W i^\frac{3}{2}\right)}\right] e^{i\omega t},
\end{gather*}
with $y = z/R$. The amplitude of the oscillation is chosen to be equal to the mean velocity $U_0$, hence $\hat{P}= U_0\rho \omega$.
The simulation starts with 
\begin{gather*}
u(x,z,\varphi,0) = 2 Q\frac{ \left( R_0^2 - z^2 \right)}{R_0^4} + u_W(x,z,0)
\end{gather*}

The chosen parameters of the simulations are $R_e = 600$, $\alpha_W=17.17$, $U_0 = R_e \nu/2 R_0$ and $\omega =  \nu\left(\alpha_W/R_0\right)^2$, by definition. The $\theta$-method is run with $\theta=0.5$ and $\theta'=1.0$.
For the present test, the discretization numbers are $N_x=63$, $N_z=40$ and $N_\varphi=48$. Simulation data are collected after five periods $\tau = 2\pi/\omega$ in order to ensure periodicity with a time-step size so that $\omega \Delta t = 5^\circ$.

\begin{figure}[!htbp]
\centering 
		\begin{subfigure}{\textwidth}
		\centering
					\begin{subfigure}{0.40\textwidth}
					\centering
					\includegraphics[width=\textwidth]{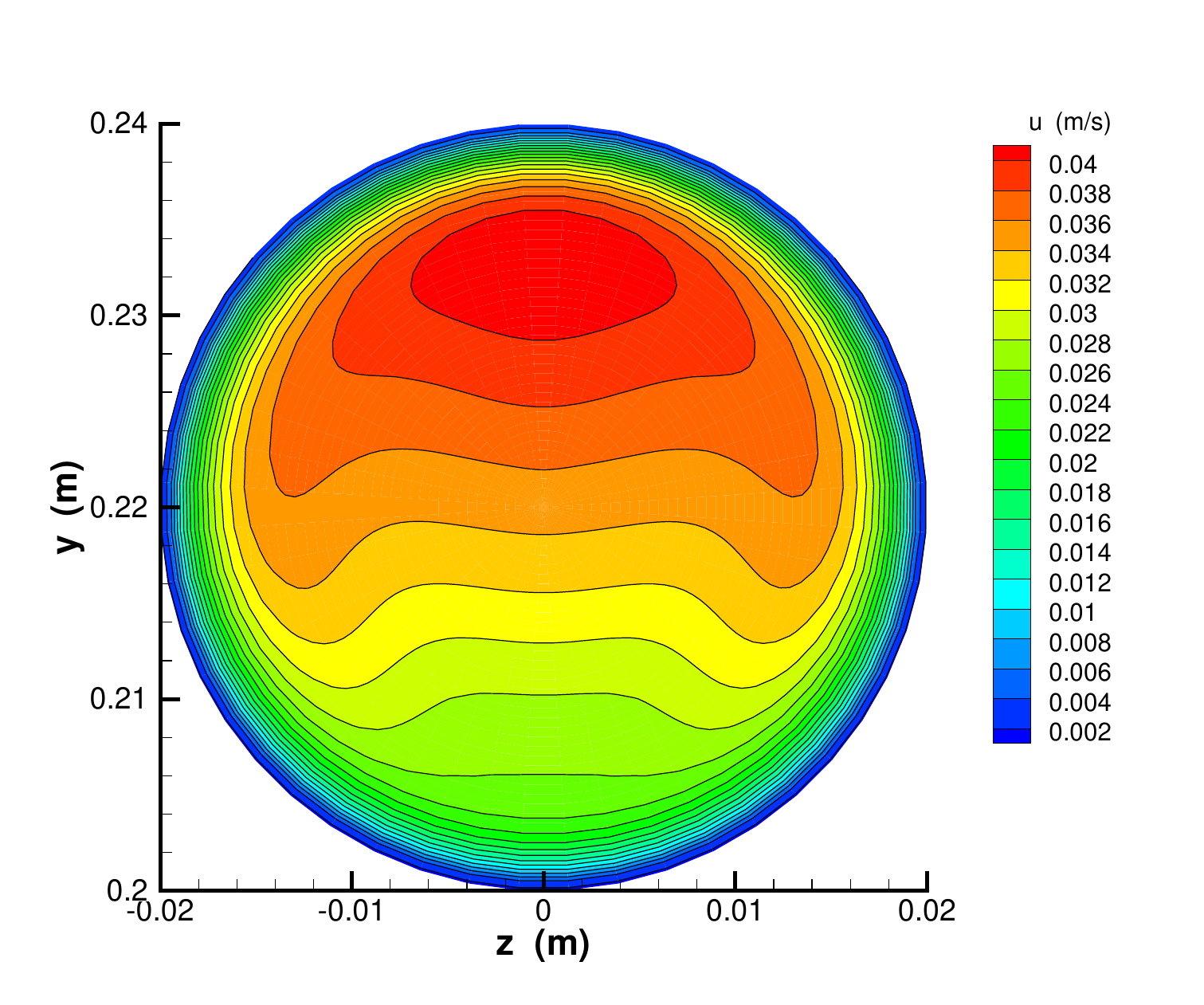}
					\end{subfigure} 
					\begin{subfigure}{0.40\textwidth}
					\centering
					\includegraphics[width=\textwidth]{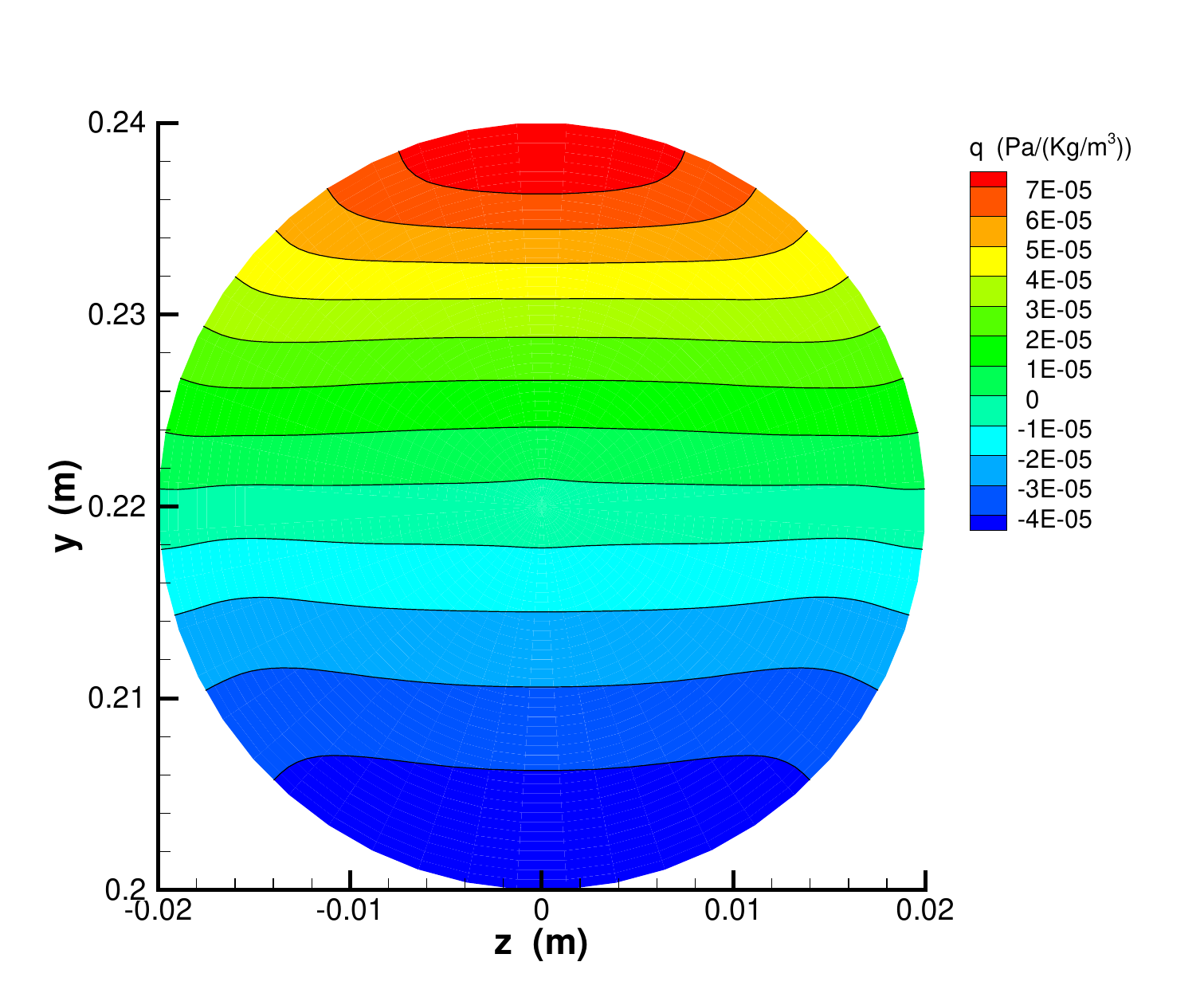}
					\end{subfigure}
					\caption{$\omega t = 90^\circ$}
		\end{subfigure}
		\begin{subfigure}{\textwidth}
		\centering
					\begin{subfigure}{0.40\textwidth}
					\centering
					\includegraphics[width=\textwidth]{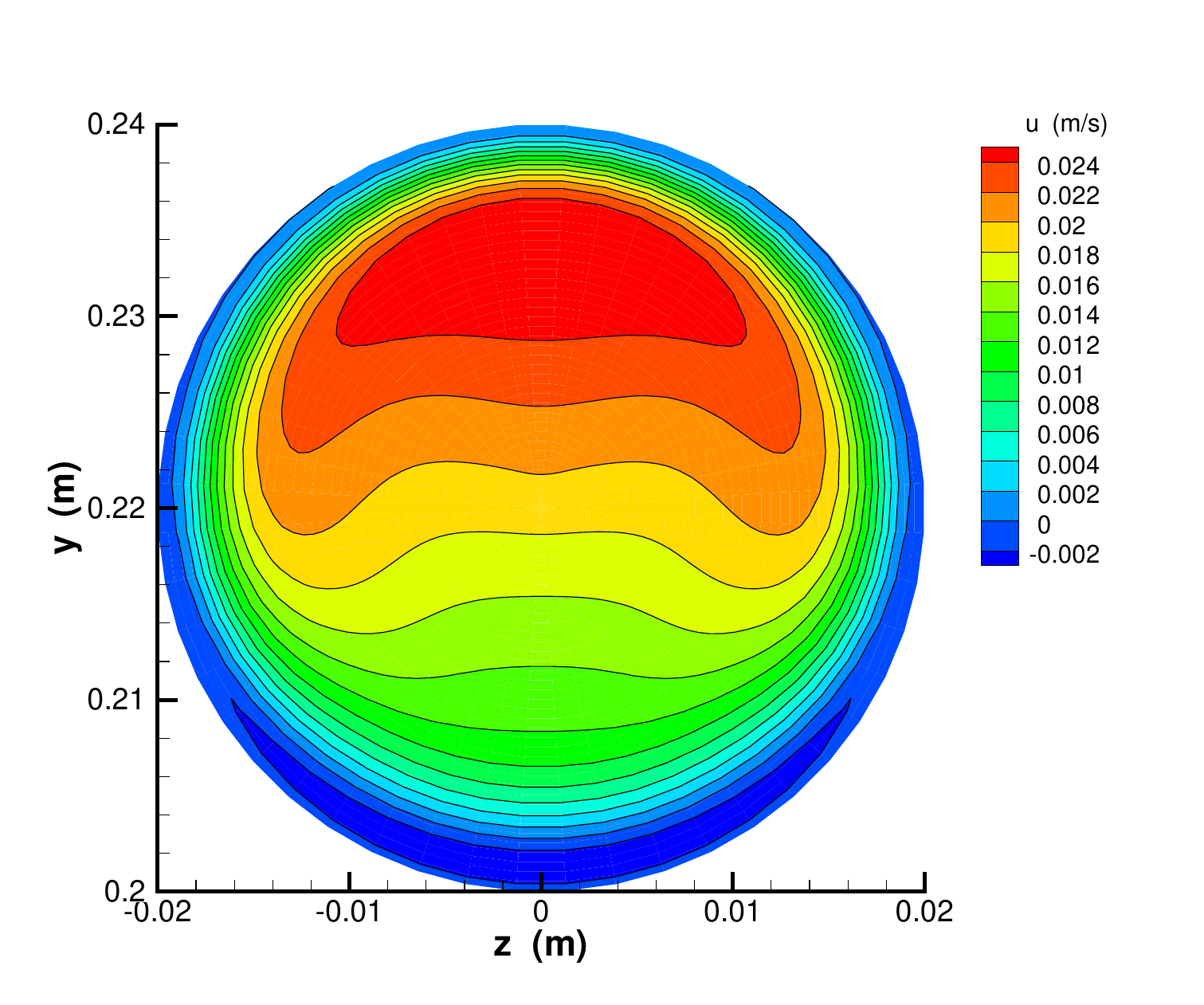}
					\end{subfigure} 
					\begin{subfigure}{0.40\textwidth}
					\centering
					\includegraphics[width=\textwidth]{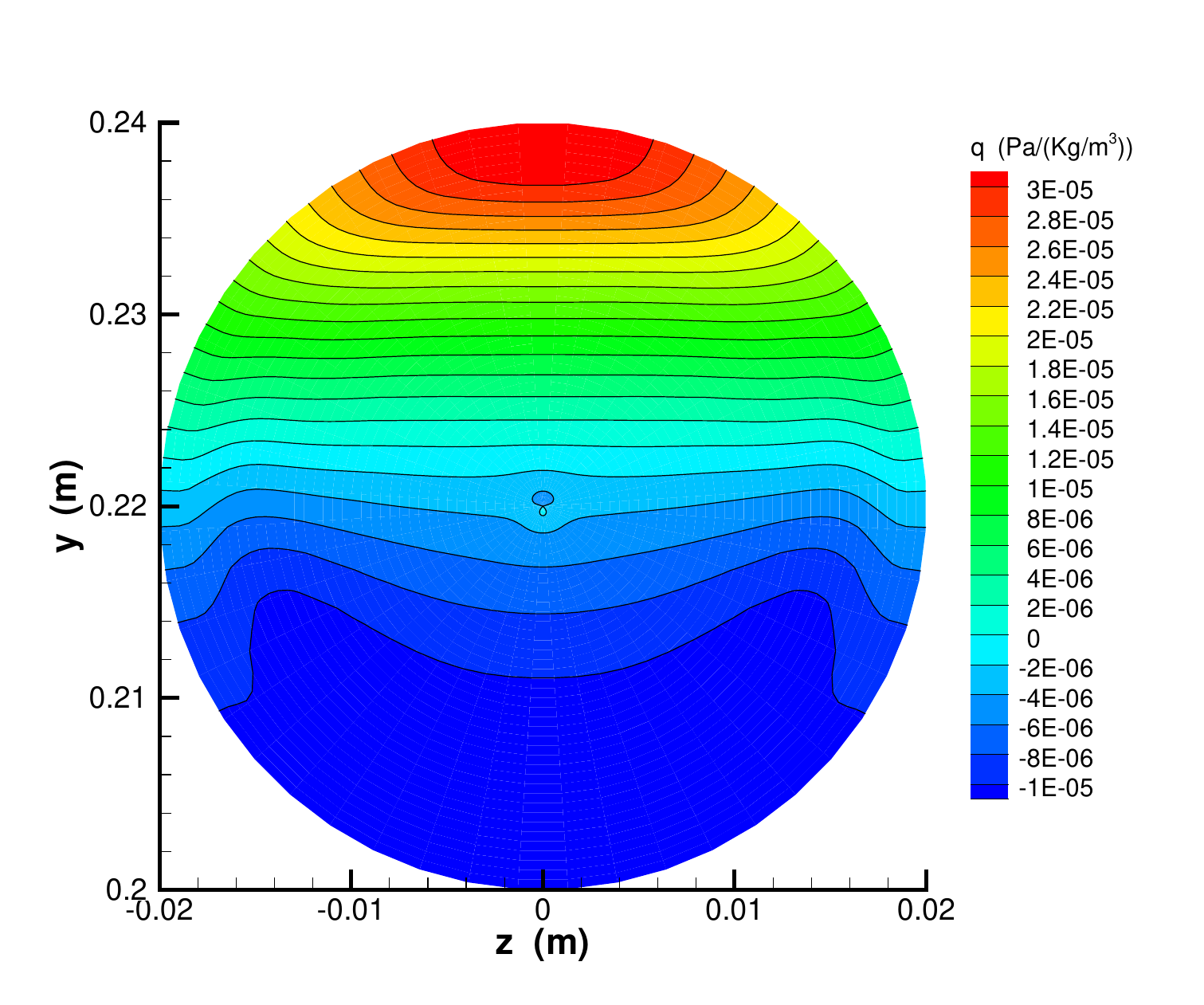}
					\end{subfigure}
					\caption{$\omega t = 180^\circ$}
		\end{subfigure}
		\begin{subfigure}{\textwidth}
		\centering
					\begin{subfigure}{0.40\textwidth}
					\centering
					\includegraphics[width=\textwidth]{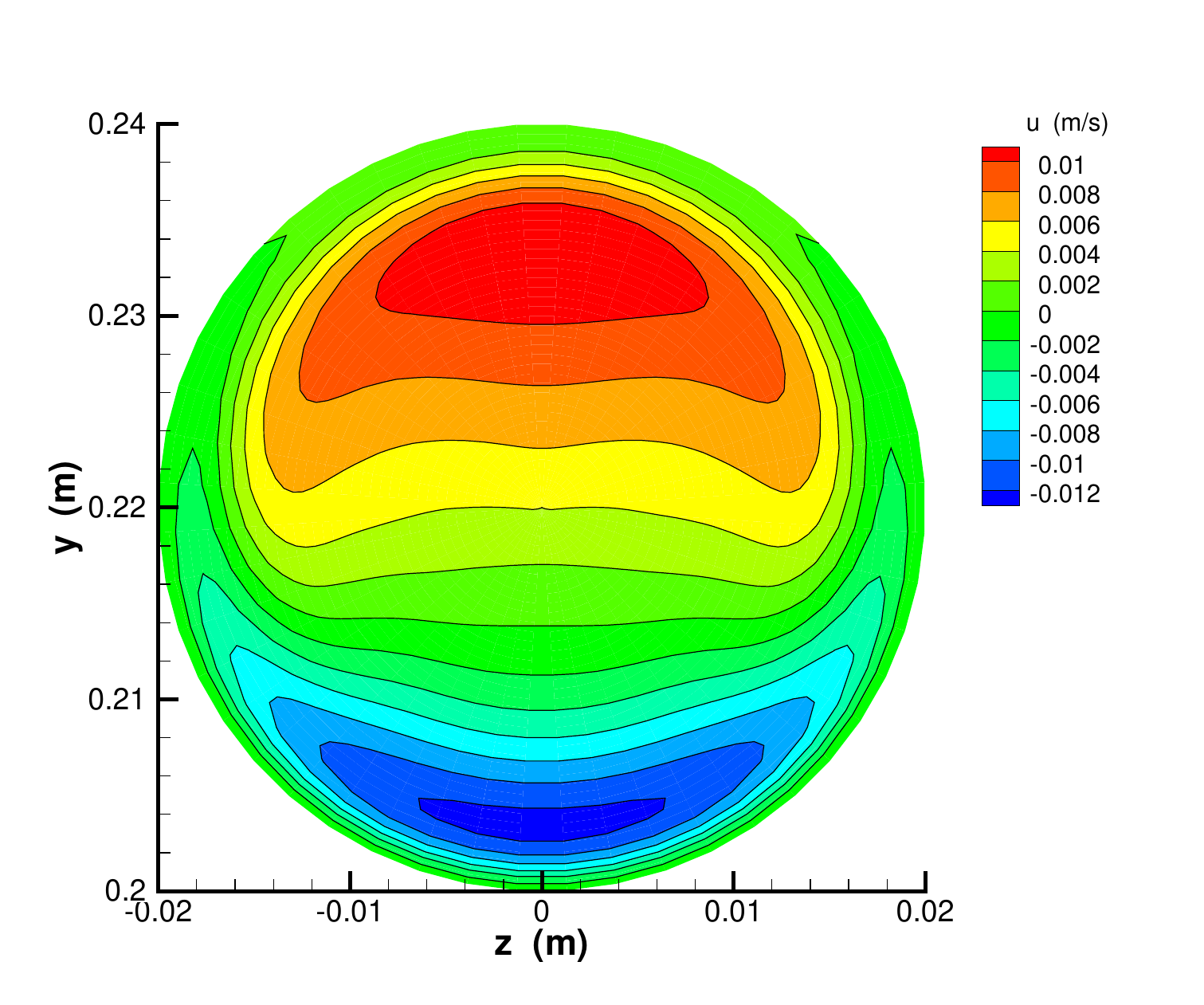}
					\end{subfigure} 
					\begin{subfigure}{0.40\textwidth}
					\centering
					\includegraphics[width=\textwidth]{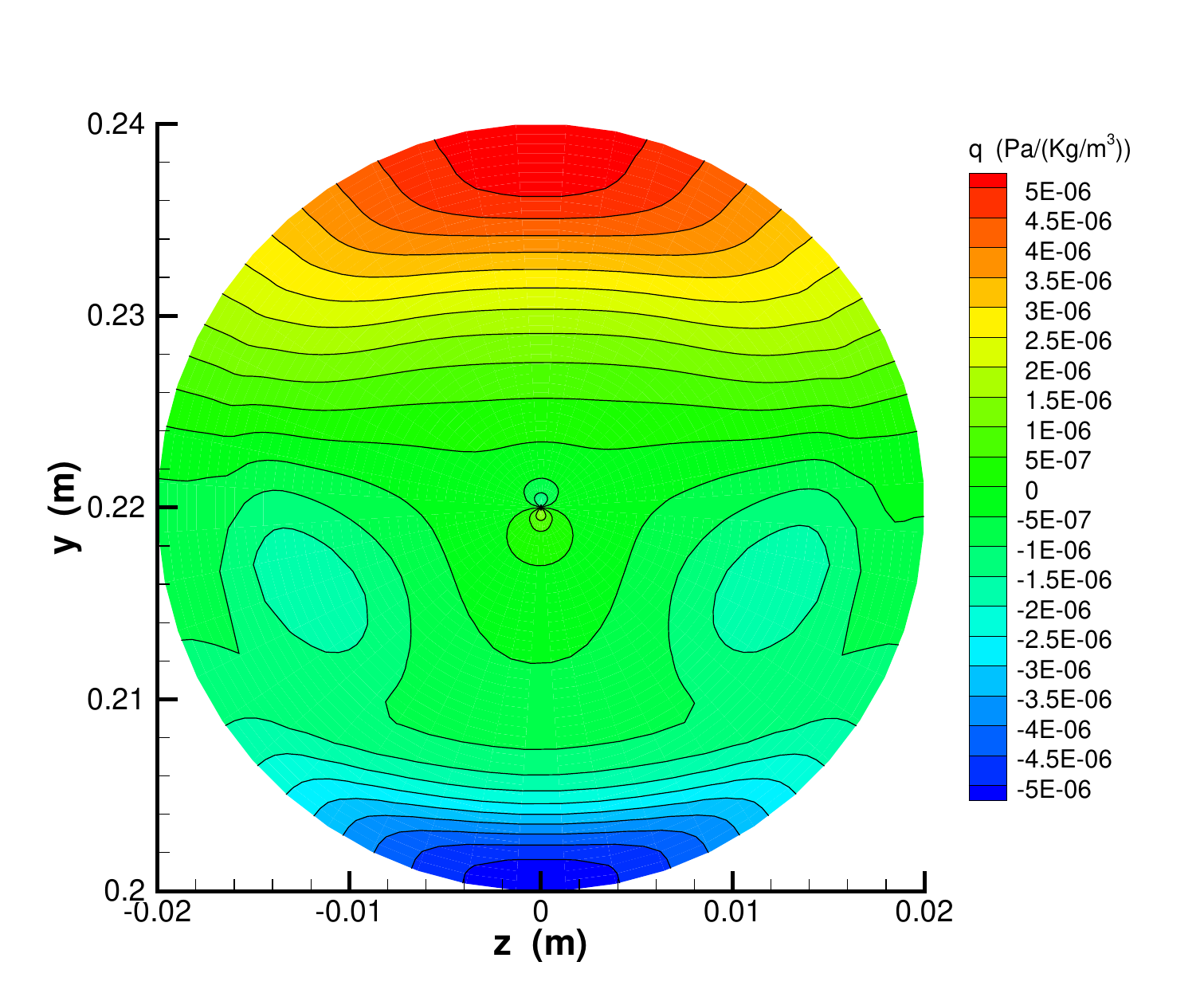}
					\end{subfigure}
					\caption{$\omega t = 270^\circ$}
		\end{subfigure}
		\begin{subfigure}{\textwidth}
		\centering
					\begin{subfigure}{0.40\textwidth}
					\centering
					\includegraphics[width=\textwidth]{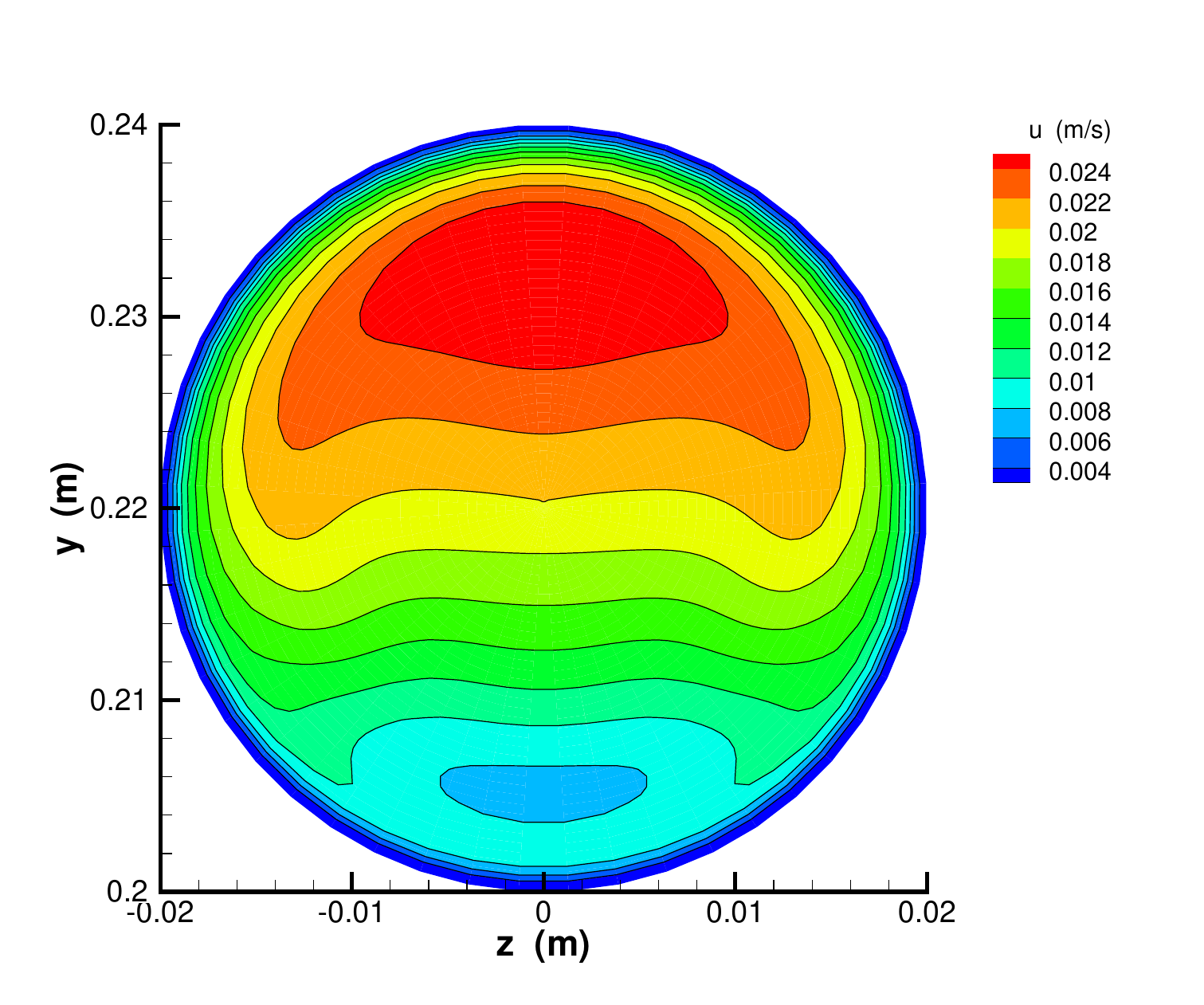}
					\end{subfigure}
					\begin{subfigure}{0.40\textwidth}
					\centering
					\includegraphics[width=\textwidth]{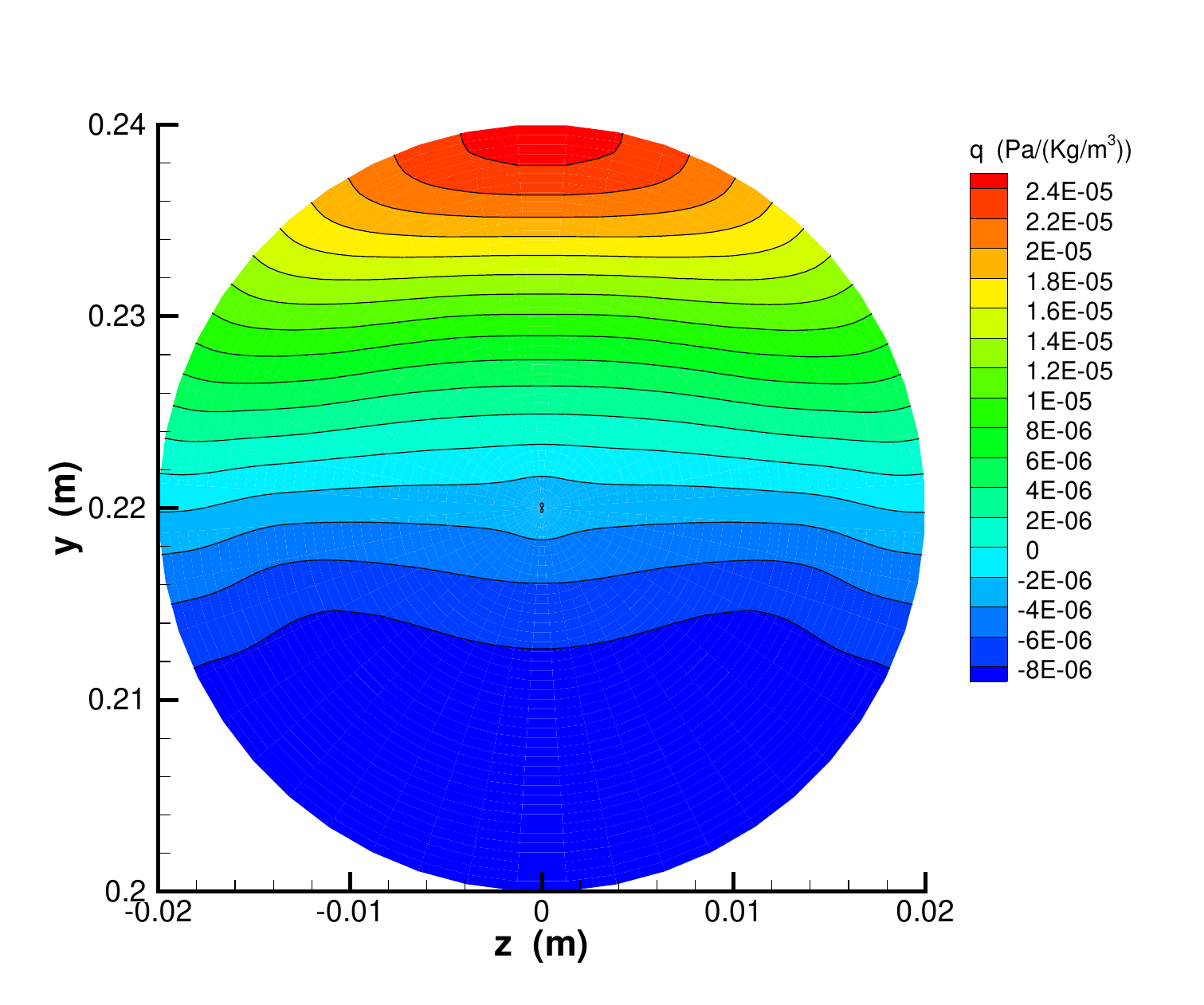}
					\end{subfigure}
					\caption{$\omega t = 360^\circ$}
		\end{subfigure}
\caption{Physical quantities interpolated along the $90^\circ$ cross section 
at different times, for $R_e = 600$ and $\alpha_W = 17.17$.  From left to 
right: axial velocity contours; non-hydrostatic correction $\mathit{q}$.}\label{fig:PulsatileUP} 
\end{figure}

\begin{figure}[!htbp]
\centering 
					\begin{subfigure}{0.4\textwidth}
						\centering
						\includegraphics[width=\textwidth]{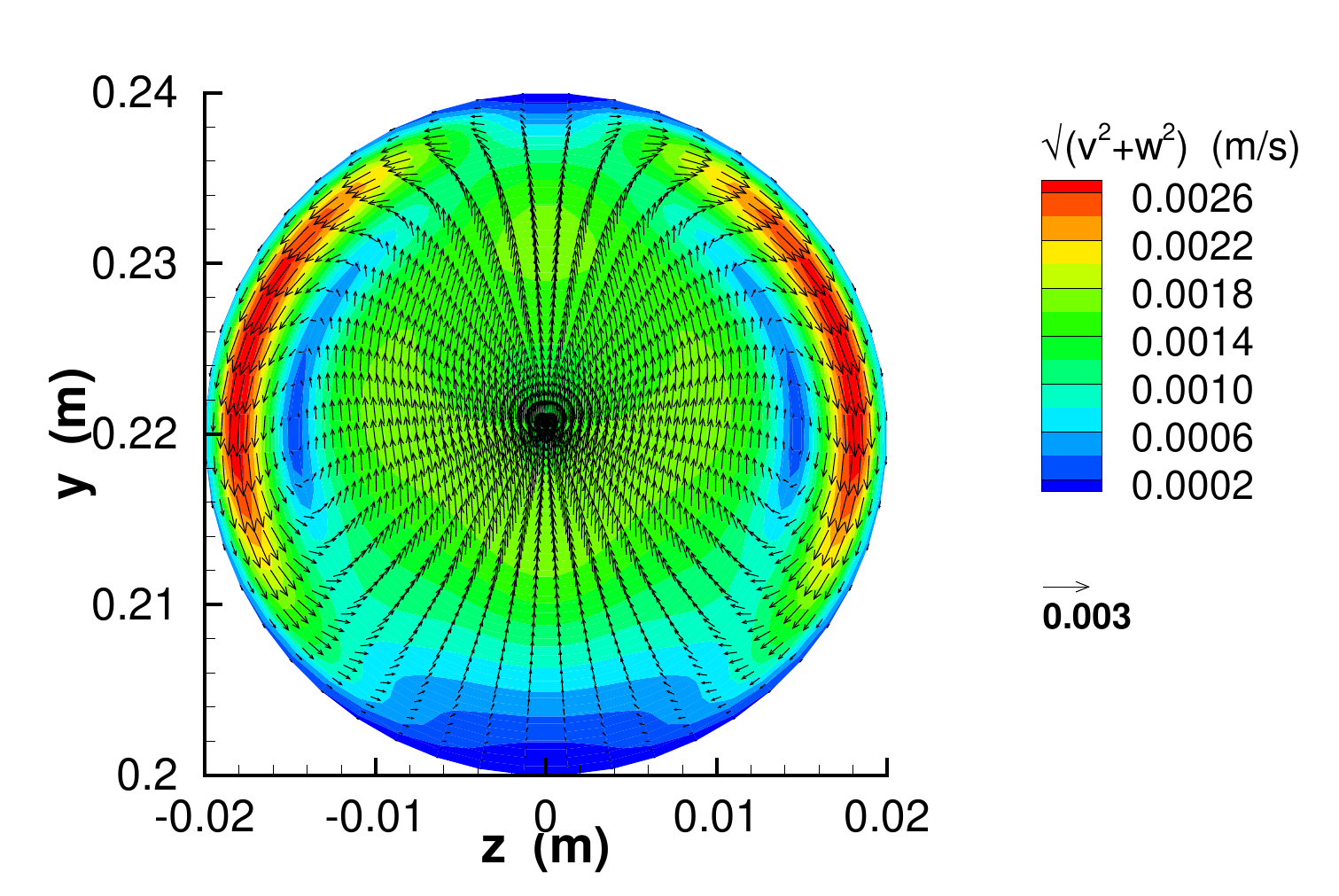}
						\caption{$\omega t = 90^\circ$}
					\end{subfigure} 
					\begin{subfigure}{0.4\textwidth}
						\centering
						\includegraphics[width=\textwidth]{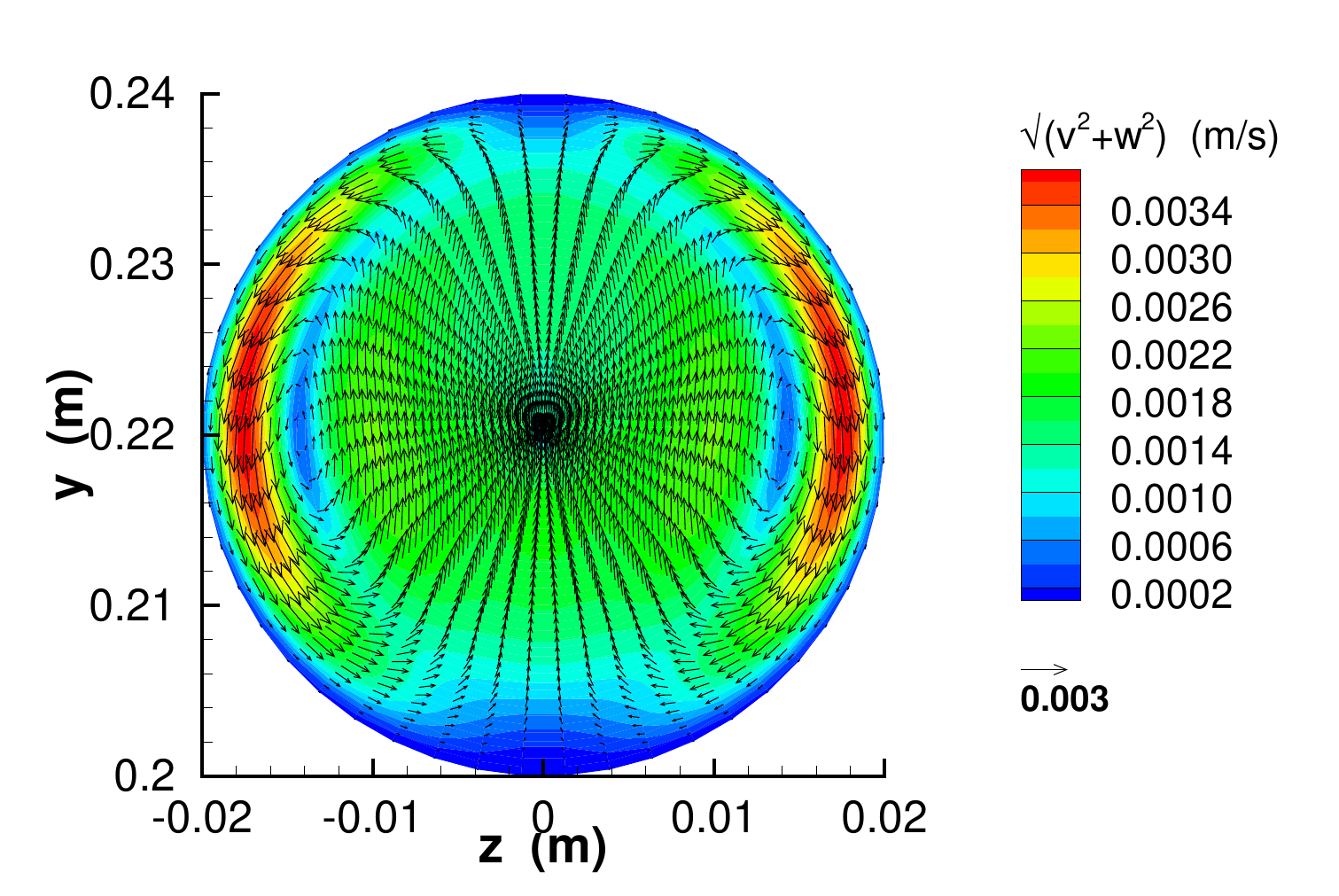}
						\caption{$\omega t = 180^\circ$}
					\end{subfigure}\\
					\begin{subfigure}{0.4\textwidth}
						\centering
						\includegraphics[width=\textwidth]{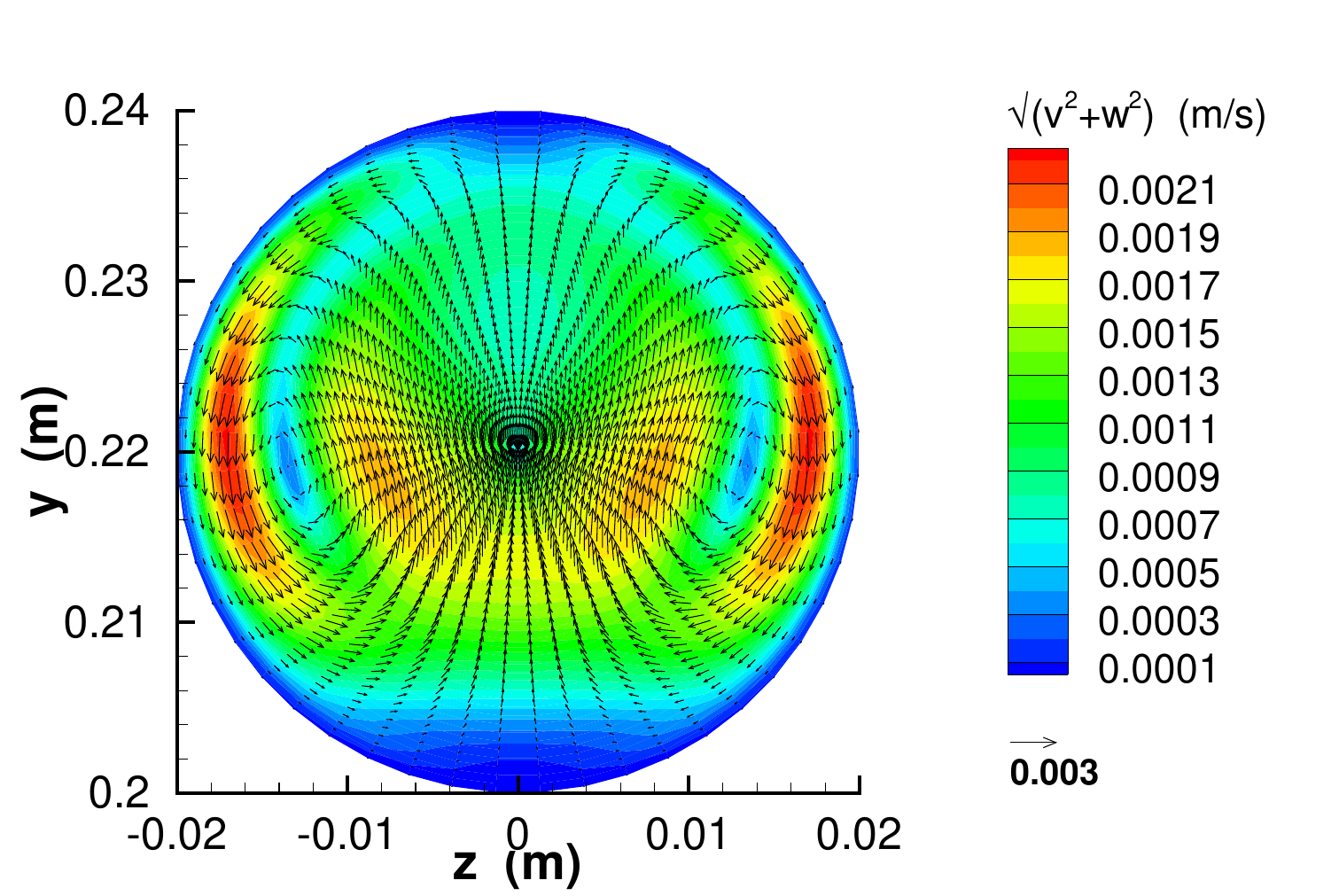}
						\caption{$\omega t = 270^\circ$}
					\end{subfigure} 
					\begin{subfigure}{0.4\textwidth}
						\centering
						\includegraphics[width=\textwidth]{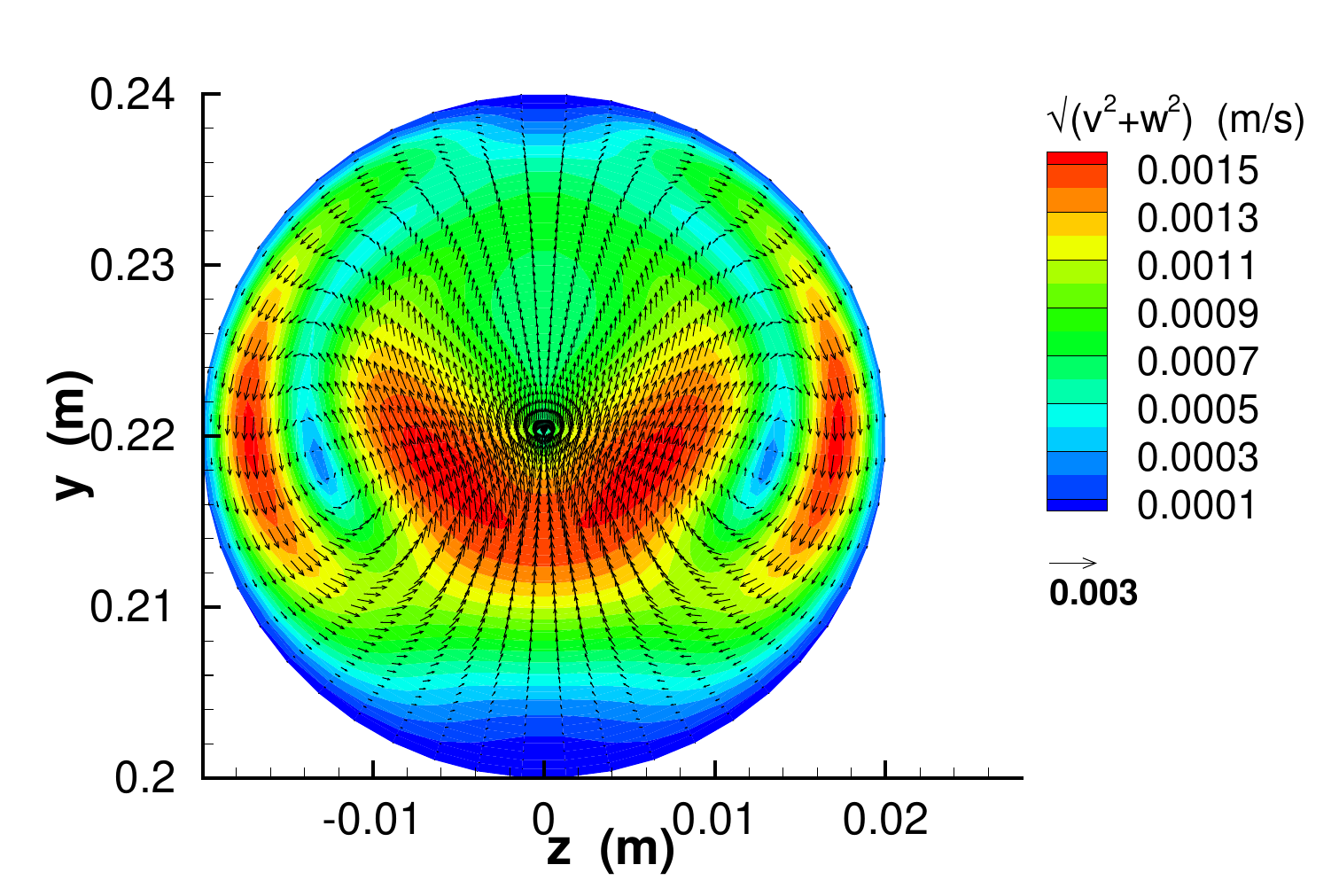}
						\caption{$\omega t = 360^\circ$}
					\end{subfigure}
\caption{Velocity vector field in the $90^\circ$ cross section for $R_e = 600$ and $\alpha_W = 17.17$.}\label{fig:PulsatileVw2} 
\end{figure}

Figures \ref{fig:PulsatileUP}-\ref{fig:PulsatileVw2} show the velocity and pressure fields interpolated along the exit cross section at different times, throughout one oscillation. Secondary motions are well established as two  pairs of symmetric, counter-rotating vortices, the intensity of which is oscillating with the Womersley frequency. 

\begin{figure}[!htbp]
\centering 
		\begin{subfigure}{\textwidth}
		\centering
					\begin{subfigure}{0.45\textwidth}
					\centering
					\includegraphics[width=\textwidth]{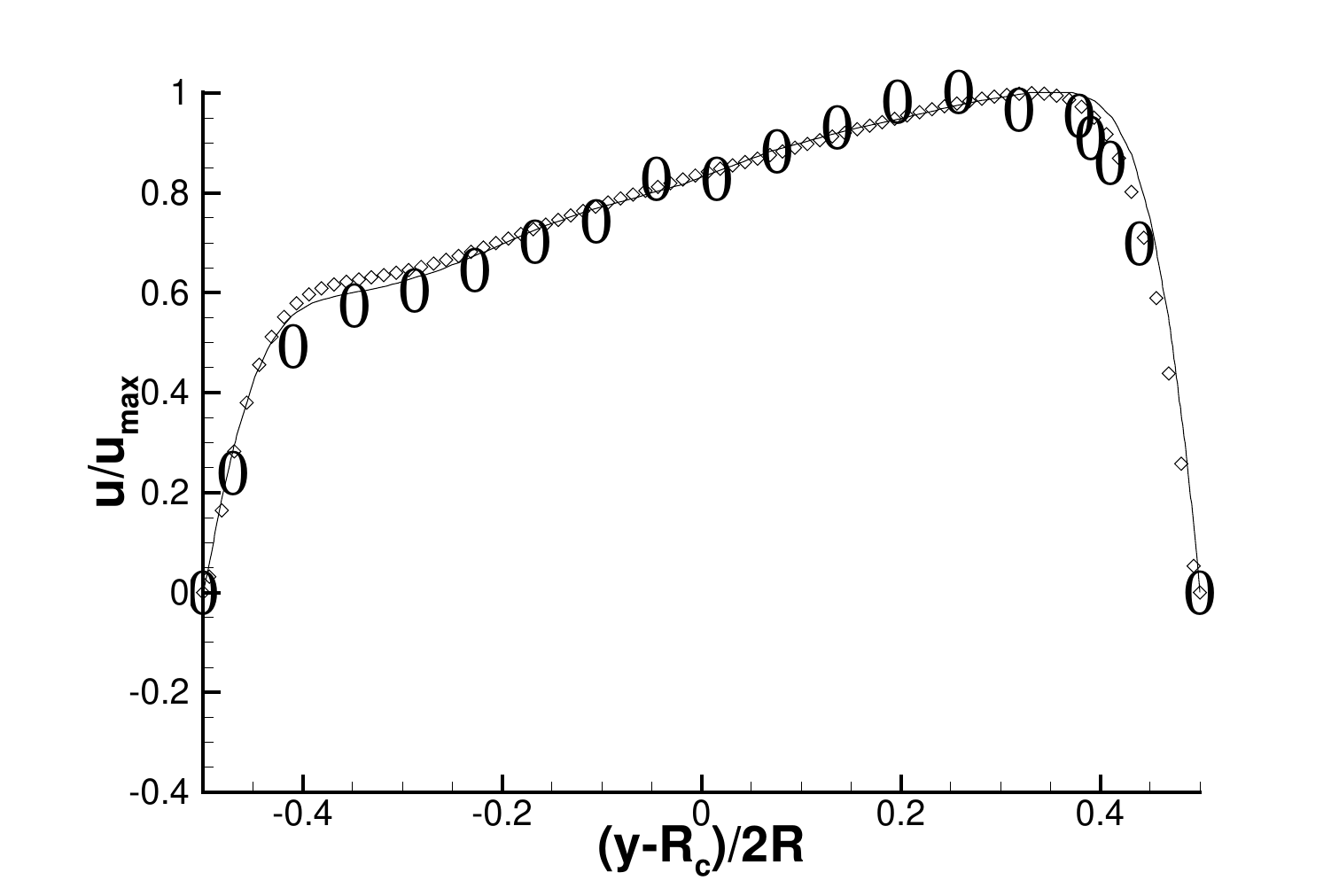}
					\end{subfigure} 
					\begin{subfigure}{0.45\textwidth}
					\centering
					\includegraphics[width=\textwidth]{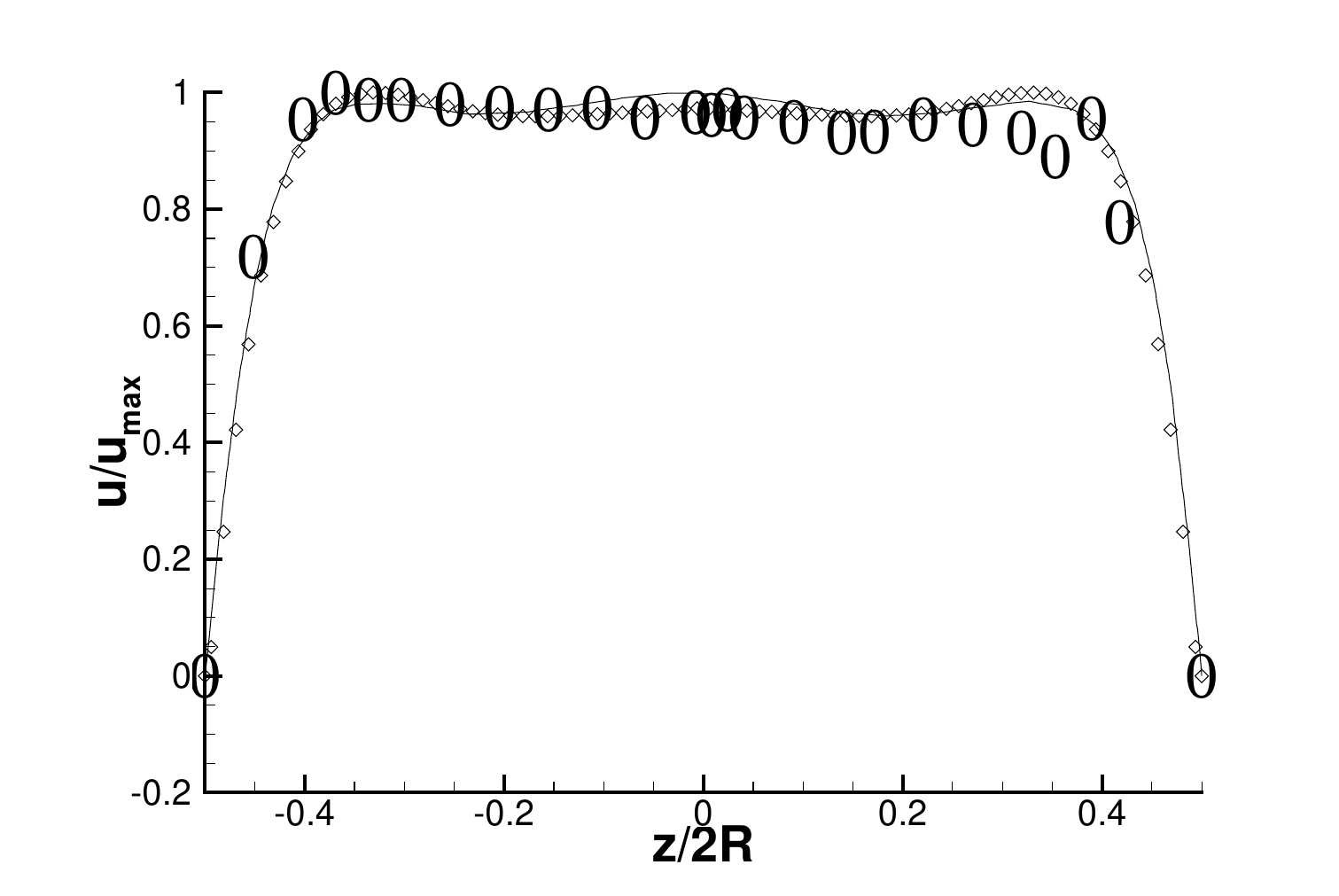}
					\end{subfigure}
					\caption{$\omega t = 90^\circ$}
		\end{subfigure}
		\begin{subfigure}{\textwidth}
		\centering
					\begin{subfigure}{0.45\textwidth}
					\centering
					\includegraphics[width=\textwidth]{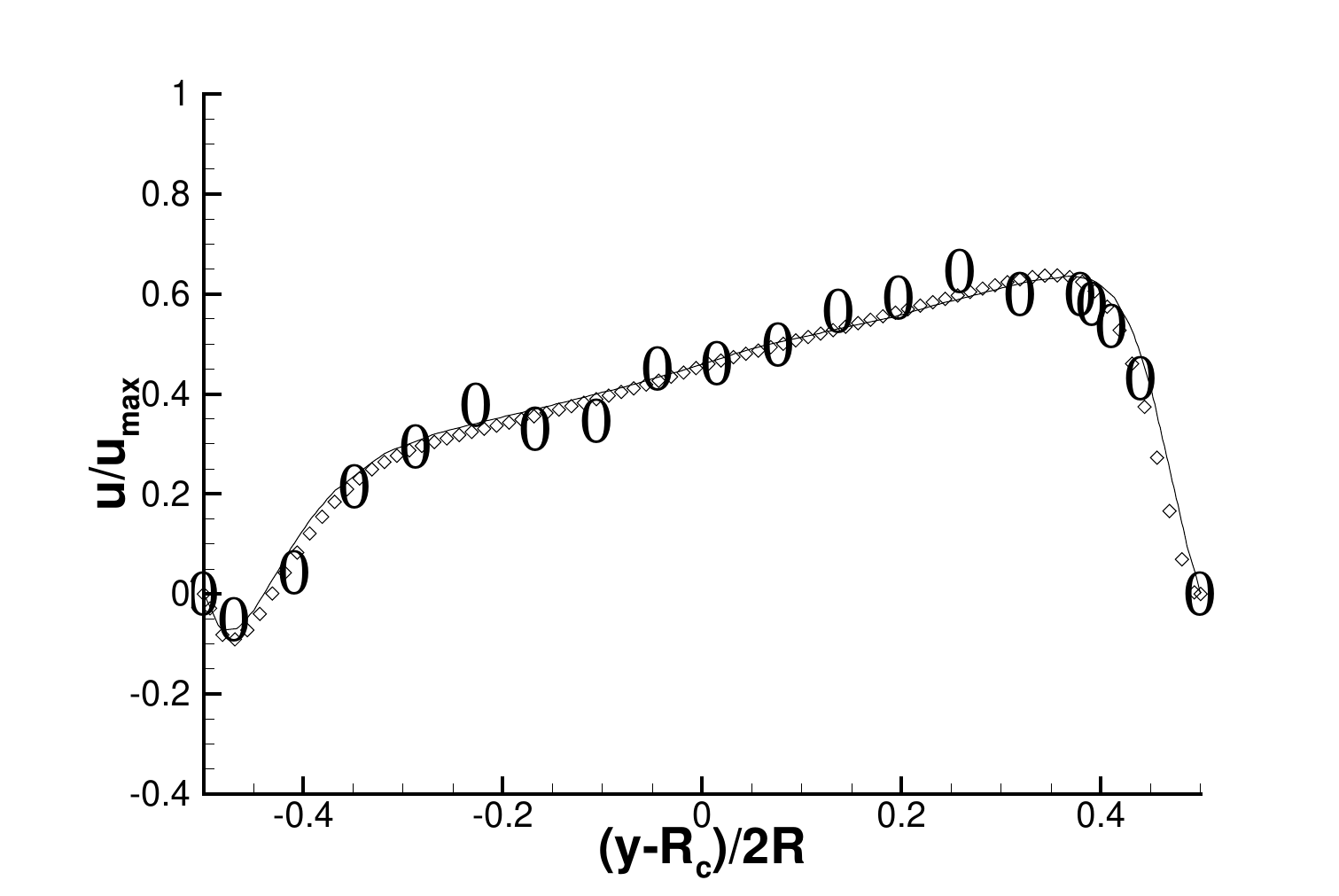}
					\end{subfigure} 
					\begin{subfigure}{0.45\textwidth}
					\centering
					\includegraphics[width=\textwidth]{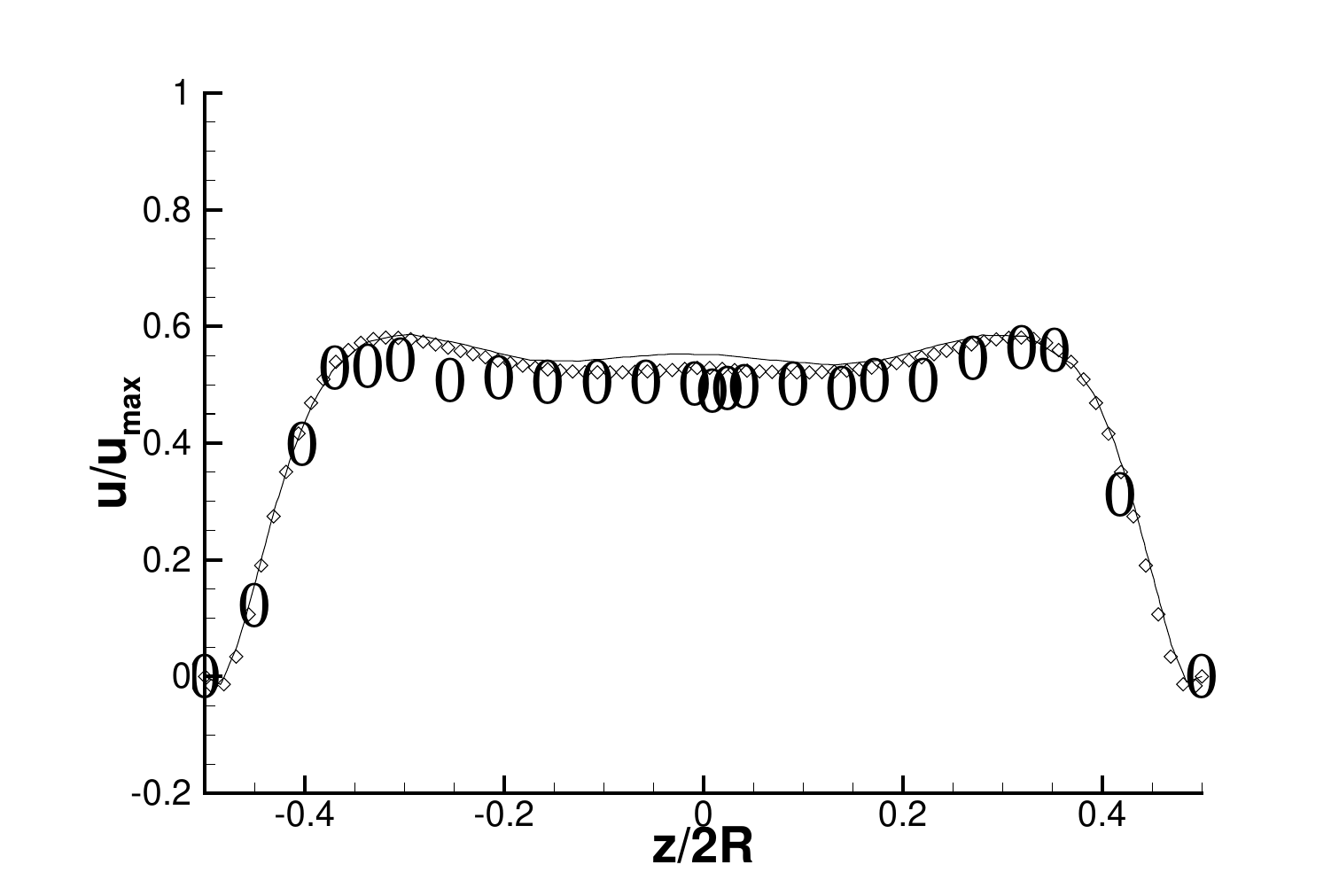}
					\end{subfigure}
					\caption{$\omega t = 180^\circ$}
		\end{subfigure}
		\begin{subfigure}{\textwidth}
		\centering
					\begin{subfigure}{0.45\textwidth}
					\centering
					\includegraphics[width=\textwidth]{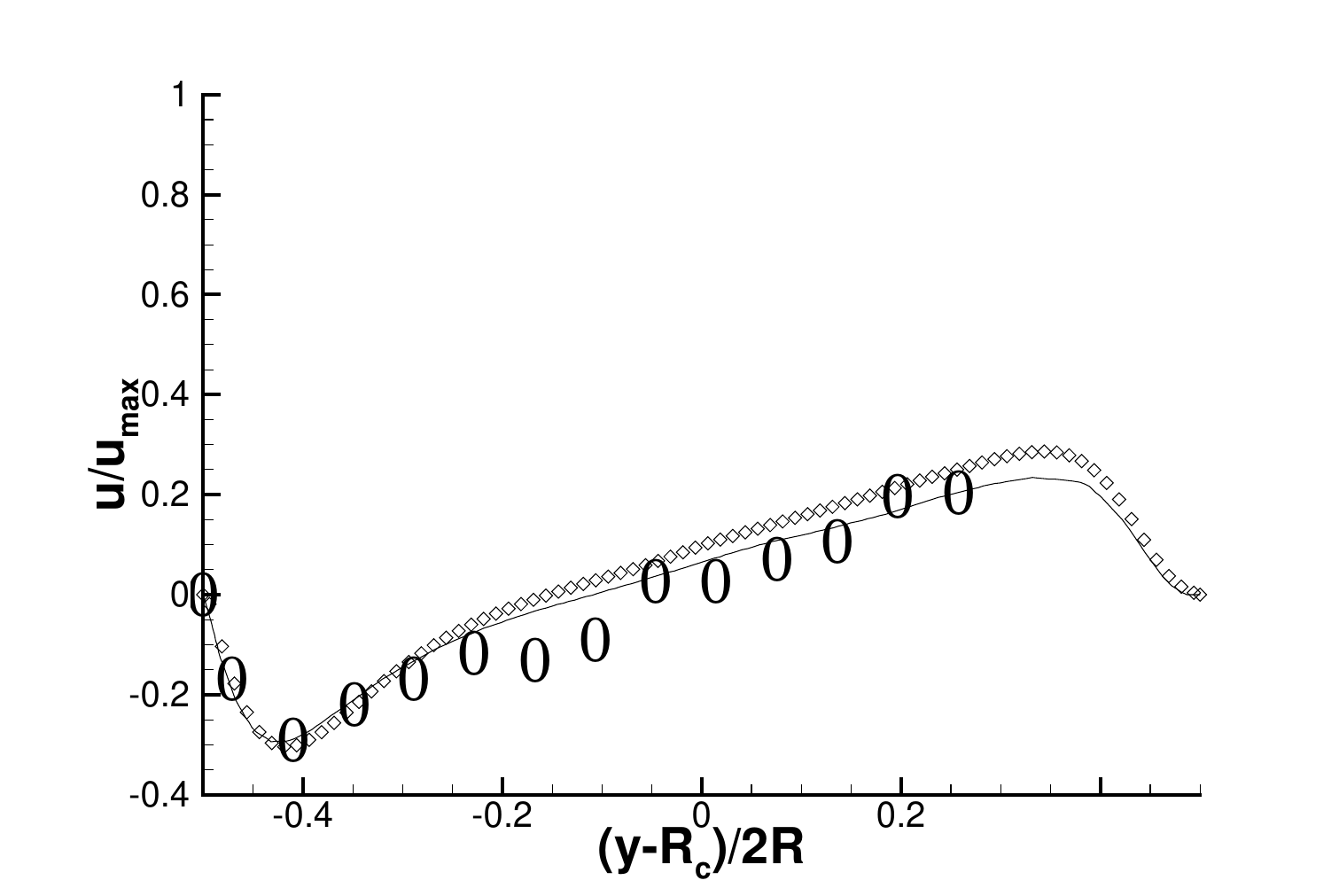}
					\end{subfigure} 
					\begin{subfigure}{0.45\textwidth}
					\centering
					\includegraphics[width=\textwidth]{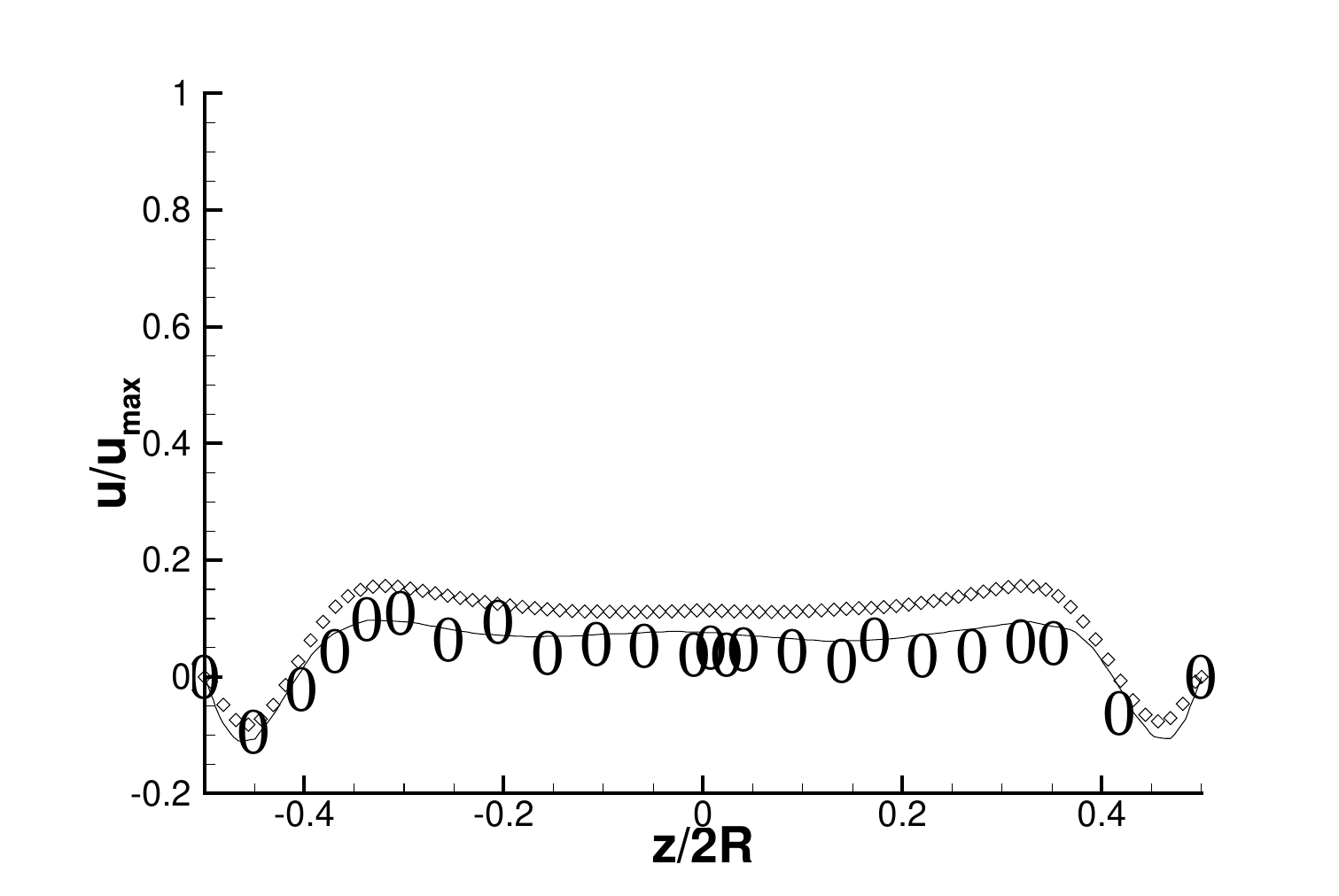}
					\end{subfigure}
					\caption{$\omega t = 270^\circ$}
		\end{subfigure}
\caption{Axial velocity profiles interpolated along the diameters of the $90^\circ$ cross section in the $\mathit{y}$- and $\mathit{z}$-direction at different times, for $R_e = 600$ and $\alpha_W = 17.17$. $\diamond$ numerical results; $-$ numerical results of \cite{Timite:2010}; $0$ experimental data of \cite{Timite:2010}.}\label{fig:cap7UYZ} 
\end{figure}

Figure \ref{fig:cap7UYZ} shows  the axial velocity components interpolated at different times along the orthogonal diameters of the $90^\circ$ cross section, next to the numerical and experimental data of \cite{Timite:2010}. Likewise, the numerical results appear to be very similar to the numerical simulations and the experimental data of \cite{Timite:2010}, also for the pulsatile flow. In particular, the reverse flow in the inner part of the 
curve is well established during the deceleration phase.

\subparagraph{Pulsatile flow in compliant tubes.} 
In order to verify the influence of the moving boundaries on the velocity and pressure fields, a pulsatile flow in a compliant tube is simulated.
The initial and the entry conditions are the same as before, while a sinusoidal pressure $p_{\text{out}}$ is applied at the exit as
\begin{gather*}
p_{\text{out}}(t) = - \frac{\hat{P}L}{\rho} \cos{\omega t},
\end{gather*}
where $L$ is the length of the tube. The parameters of the simulations are chosen to be $R_0 = 0.025 \text{m}$, $R_e = 600$, $\alpha_W=17.17$, with a rigidity coefficient $\beta=1$. The discretization numbers are again $N_x=63$, $N_z=40$ and $N_\varphi=48$; $\theta=0.5$ and $\theta'=1.0$. Numerical data are collected after five periods $\tau = 2\pi/\omega$, with a time-step size so that $\omega \Delta t = 5^\circ$.

\begin{figure}[p]
\centering 
					\begin{subfigure}[b]{0.4\textwidth}
						\centering
						\includegraphics[width=\textwidth]{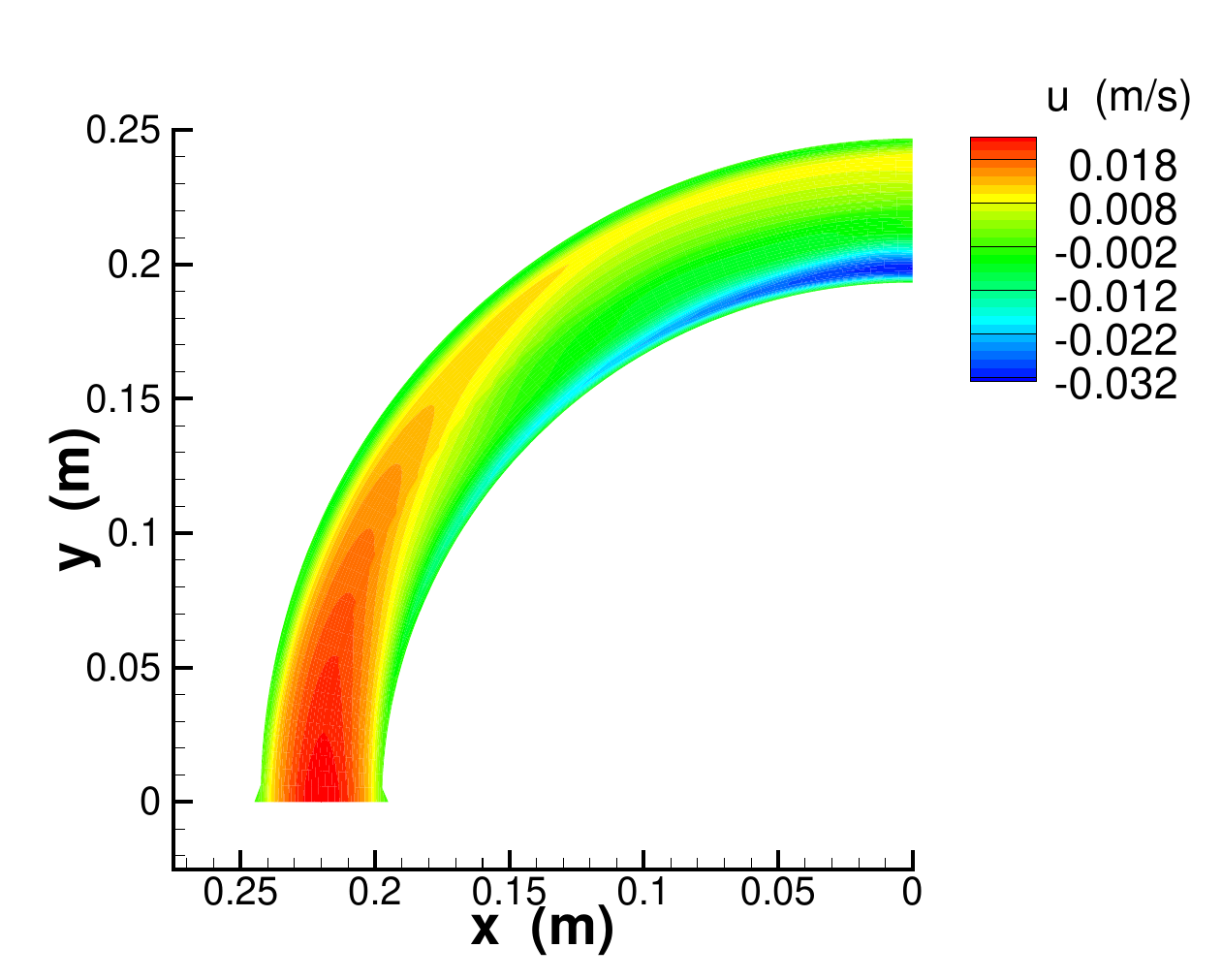}
						\caption{$\omega t = 180^\circ$}
					\end{subfigure} 
					\begin{subfigure}[b]{0.4\textwidth}
						\centering
						\includegraphics[width=\textwidth]{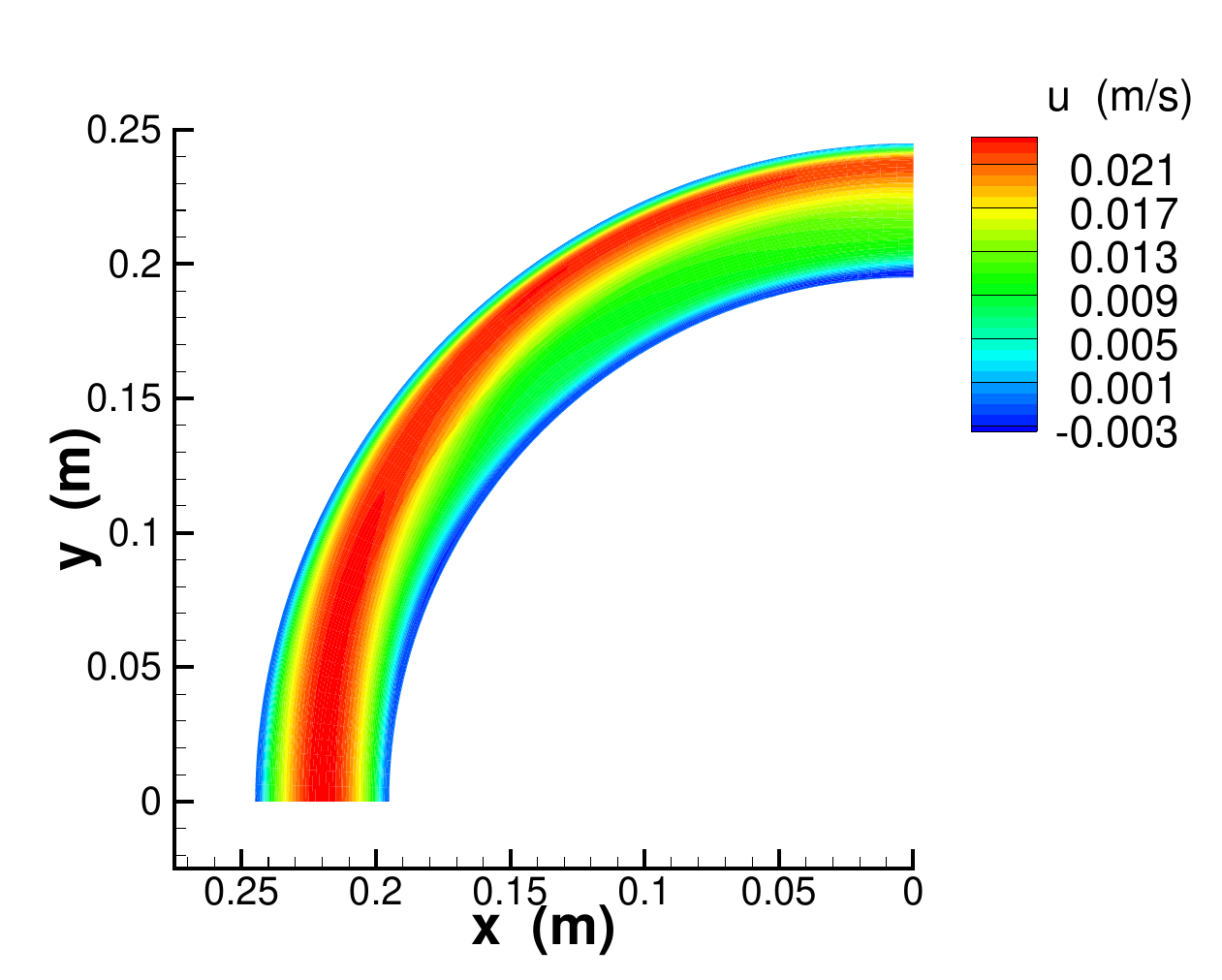}
						\caption{$\omega t = 180^\circ$}
					\end{subfigure}\\
					\begin{subfigure}[b]{0.4\textwidth}
						\centering
						\includegraphics[width=\textwidth]{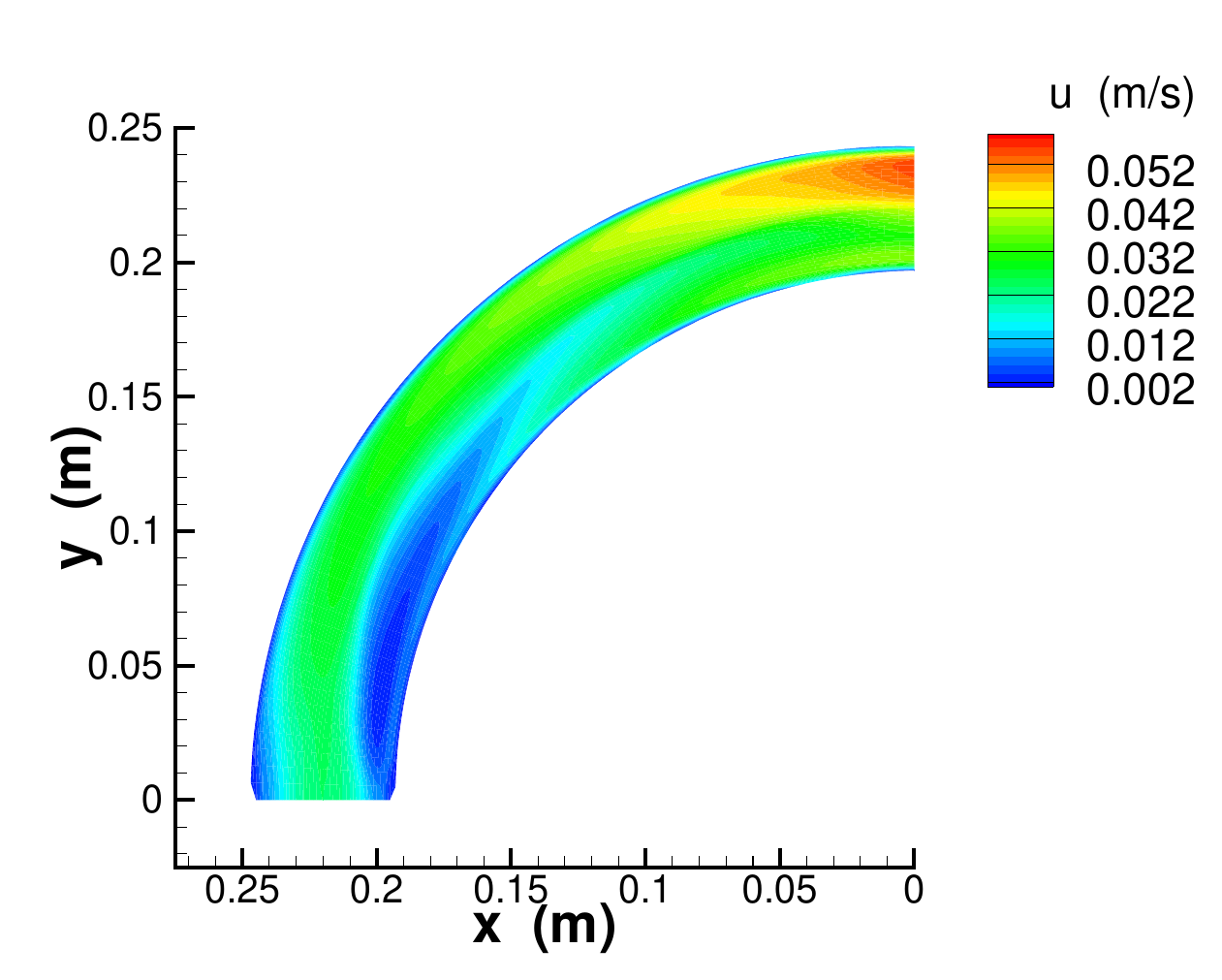}
						\caption{$\omega t = 360^\circ$}
					\end{subfigure} 
					\begin{subfigure}[b]{0.4\textwidth}
						\centering
						\includegraphics[width=\textwidth]{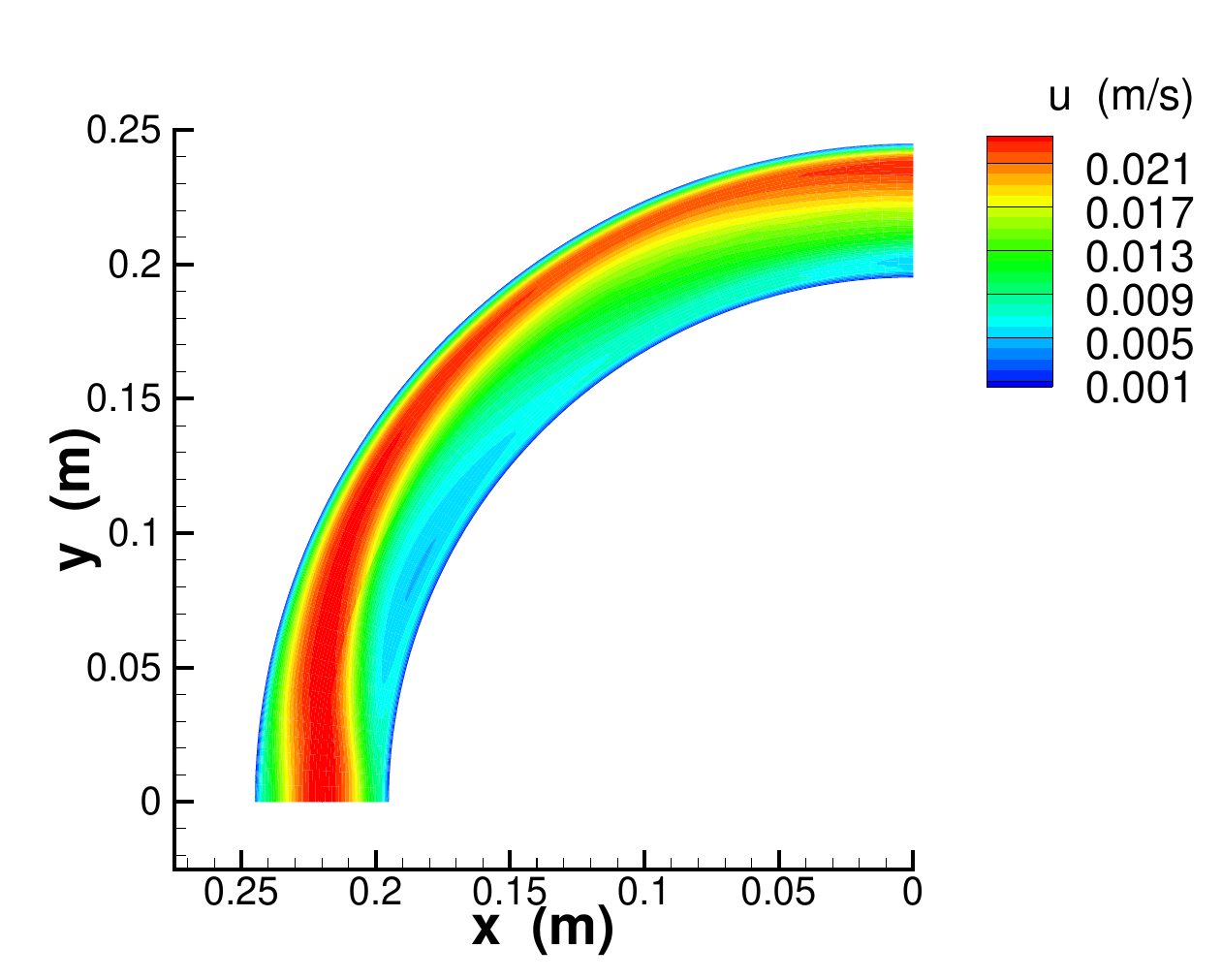}
						\caption{$\omega t = 360^\circ$}
					\end{subfigure}
\caption{Axial velocity field interpolated along the symmetry plane at different times and different rigidity coefficient, for $R_e = 600$, $\alpha_W = 17.17$, $\beta=1$ (left) and $\beta=10^{12}$ (right).}\label{fig:ComplianceProfile} 
\end{figure}

\begin{figure}[p]
\centering 
					\begin{subfigure}[]{0.40\textwidth}
						\centering
						\includegraphics[width=\textwidth]{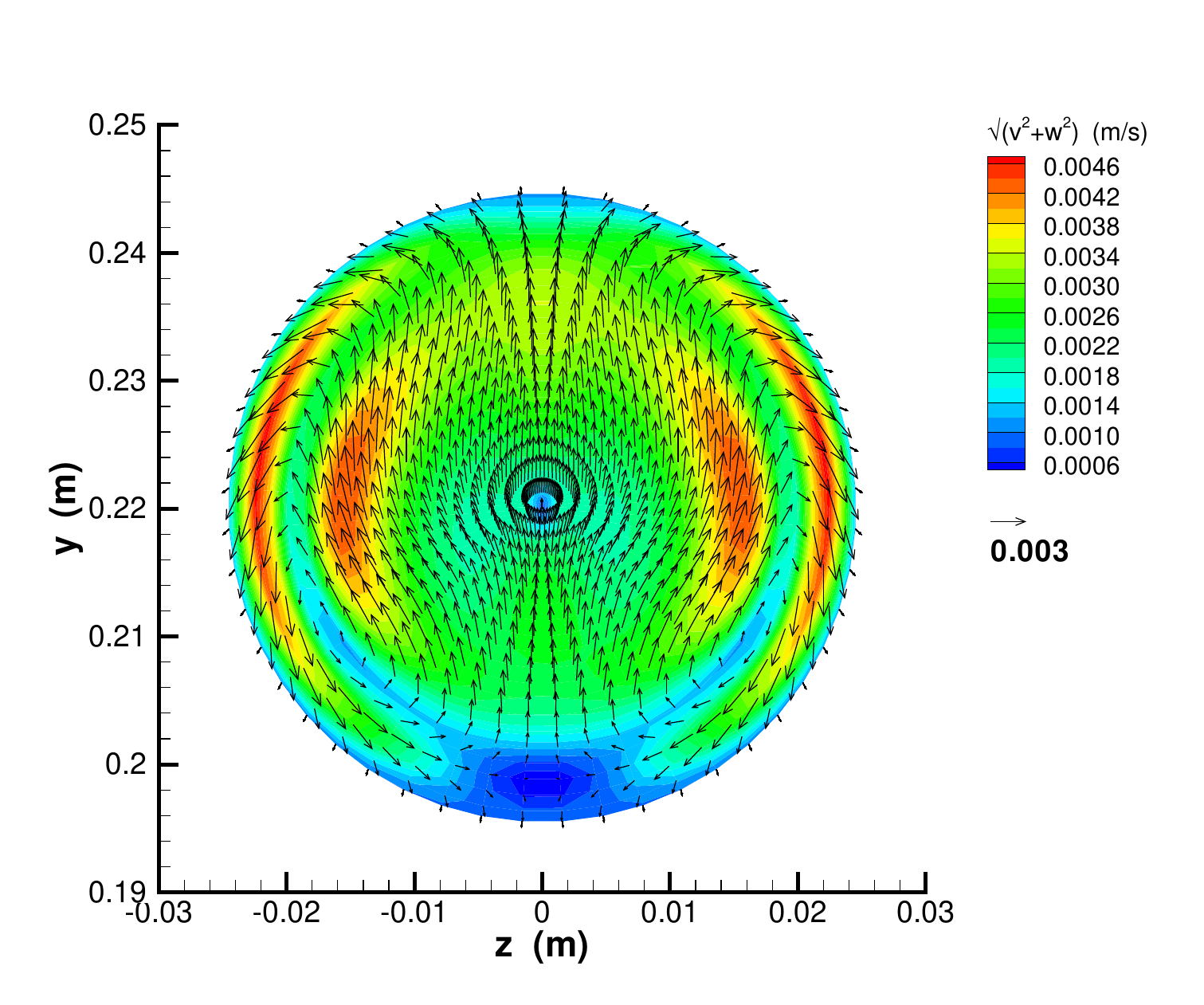}
						\caption{$\omega t = 90^\circ$}\label{fig:Compliance90B1} 
					\end{subfigure} 
					\begin{subfigure}[]{0.40\textwidth}
						\centering
						\includegraphics[width=\textwidth]{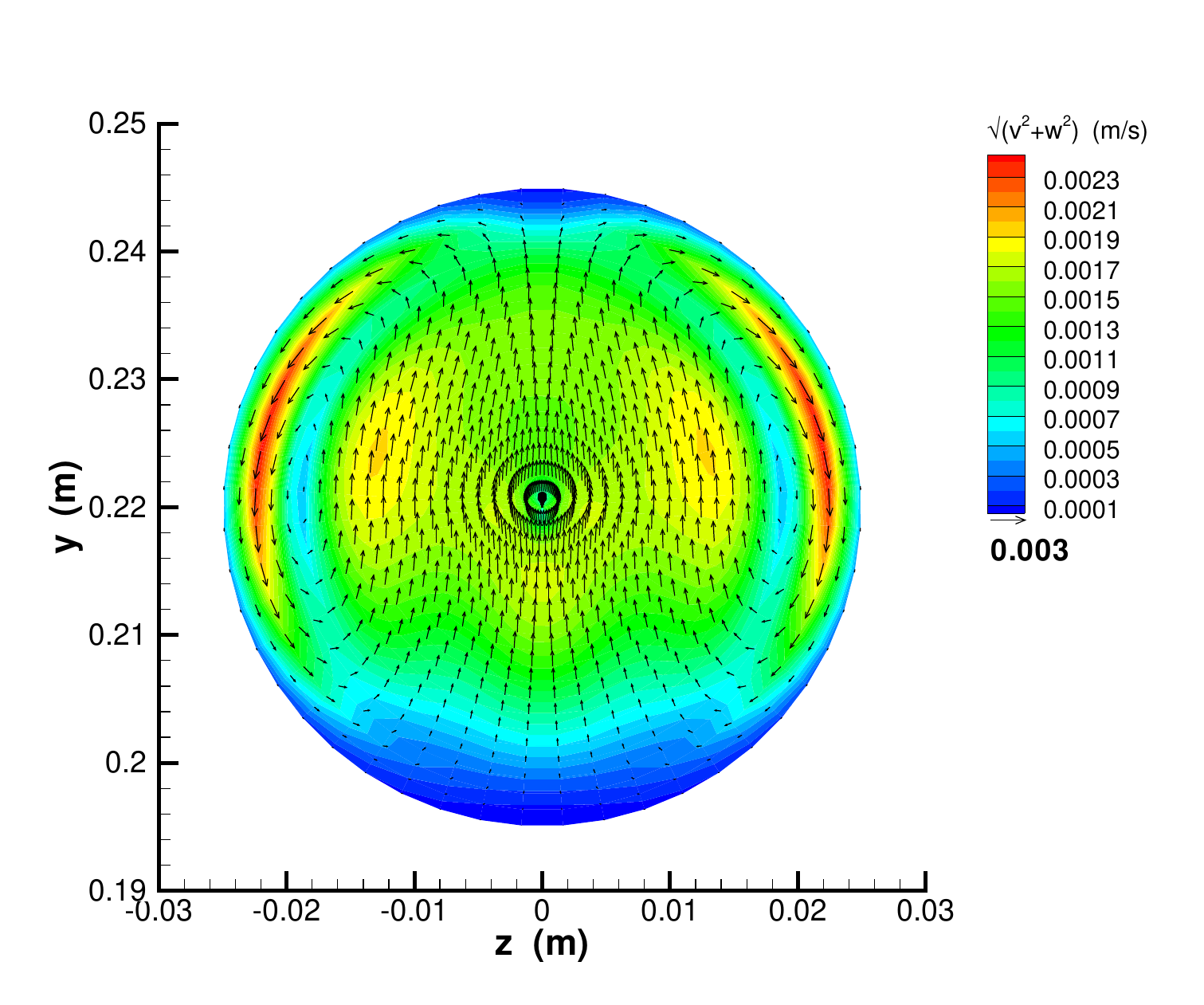}
						\caption{$\omega t = 90^\circ$}
					\end{subfigure}\\
					\begin{subfigure}[]{0.40\textwidth}
						\centering
						\includegraphics[width=\textwidth]{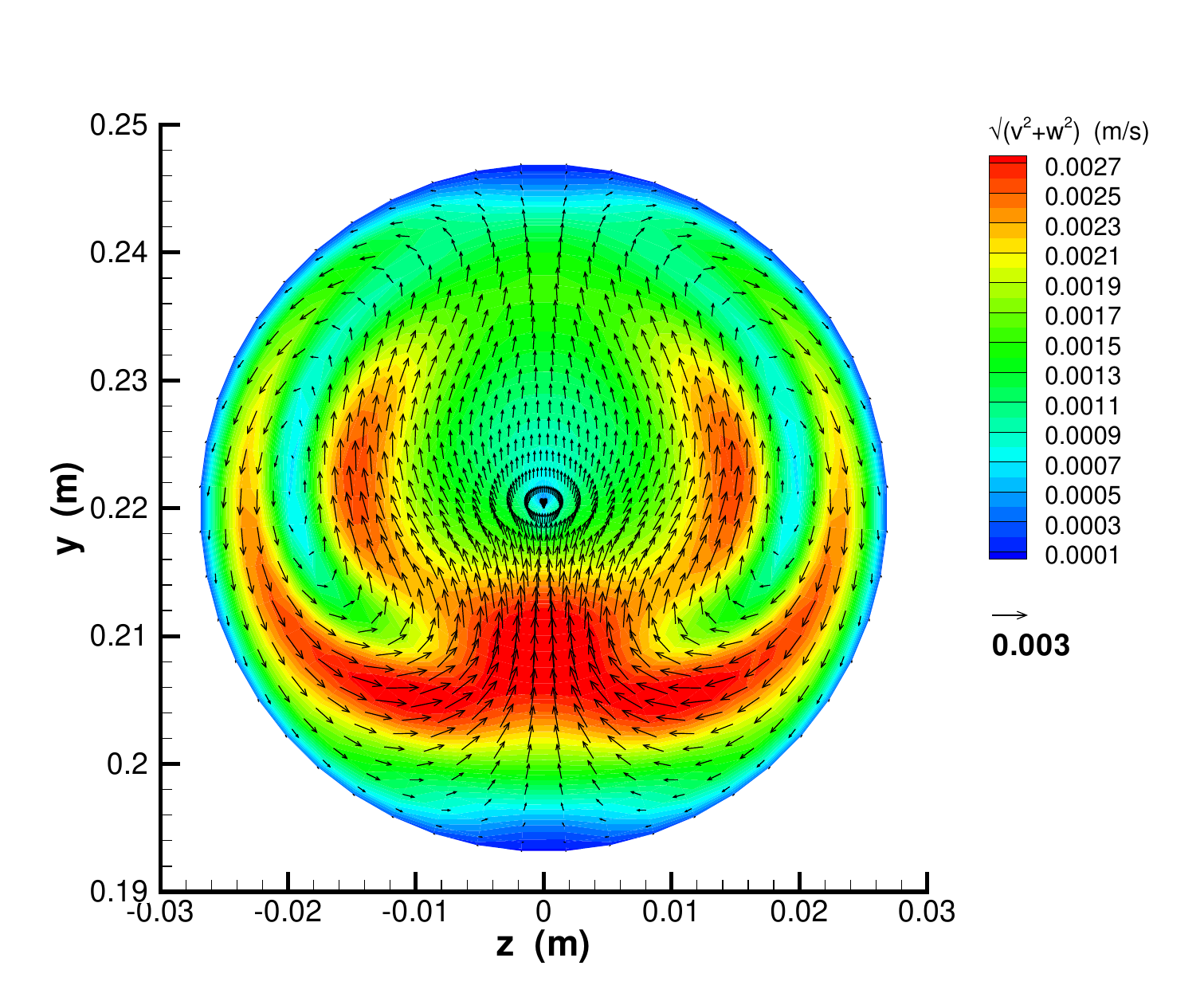}
						\caption{$\omega t = 180^\circ$}\label{fig:Compliance180B1} 
					\end{subfigure} 
					\begin{subfigure}[]{0.40\textwidth}
						\centering
						\includegraphics[width=\textwidth]{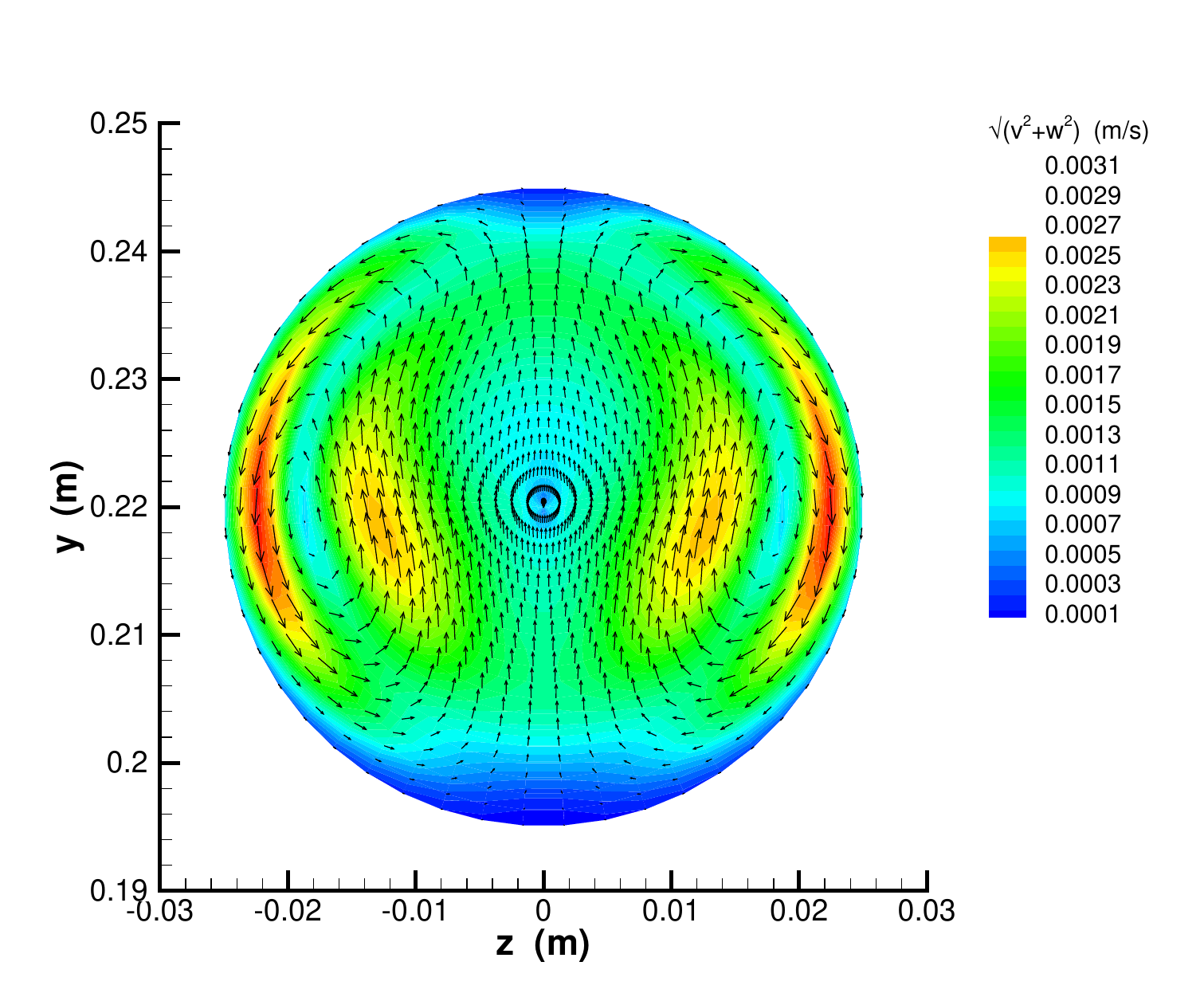}
						\caption{$\omega t = 180^\circ$}
					\end{subfigure}\\
					\begin{subfigure}[]{0.40\textwidth}
						\centering
						\includegraphics[width=\textwidth]{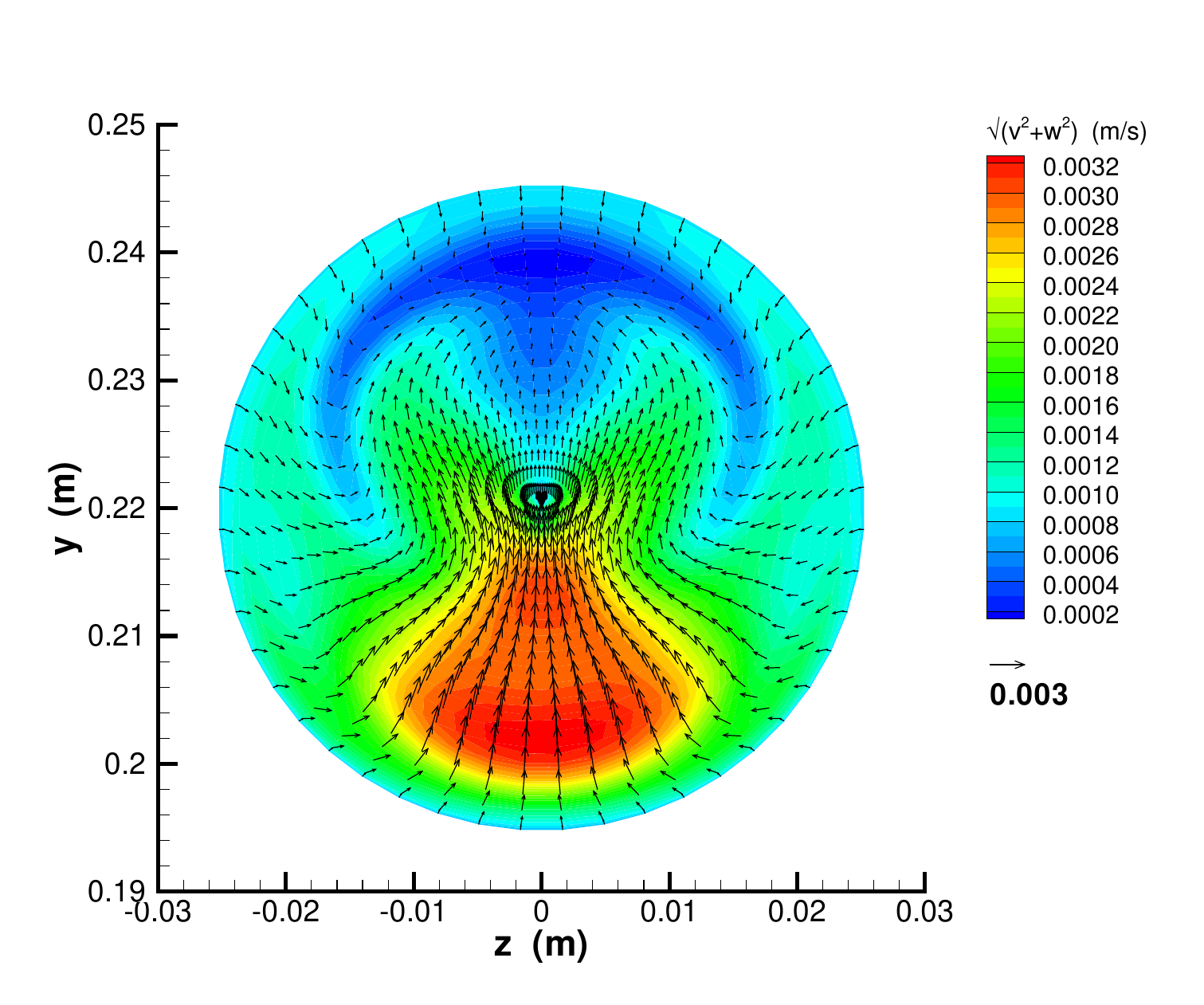}
						\caption{$\omega t = 270^\circ$}\label{fig:Compliance270B1} 
					\end{subfigure} 
					\begin{subfigure}[]{0.40\textwidth}
						\centering
						\includegraphics[width=\textwidth]{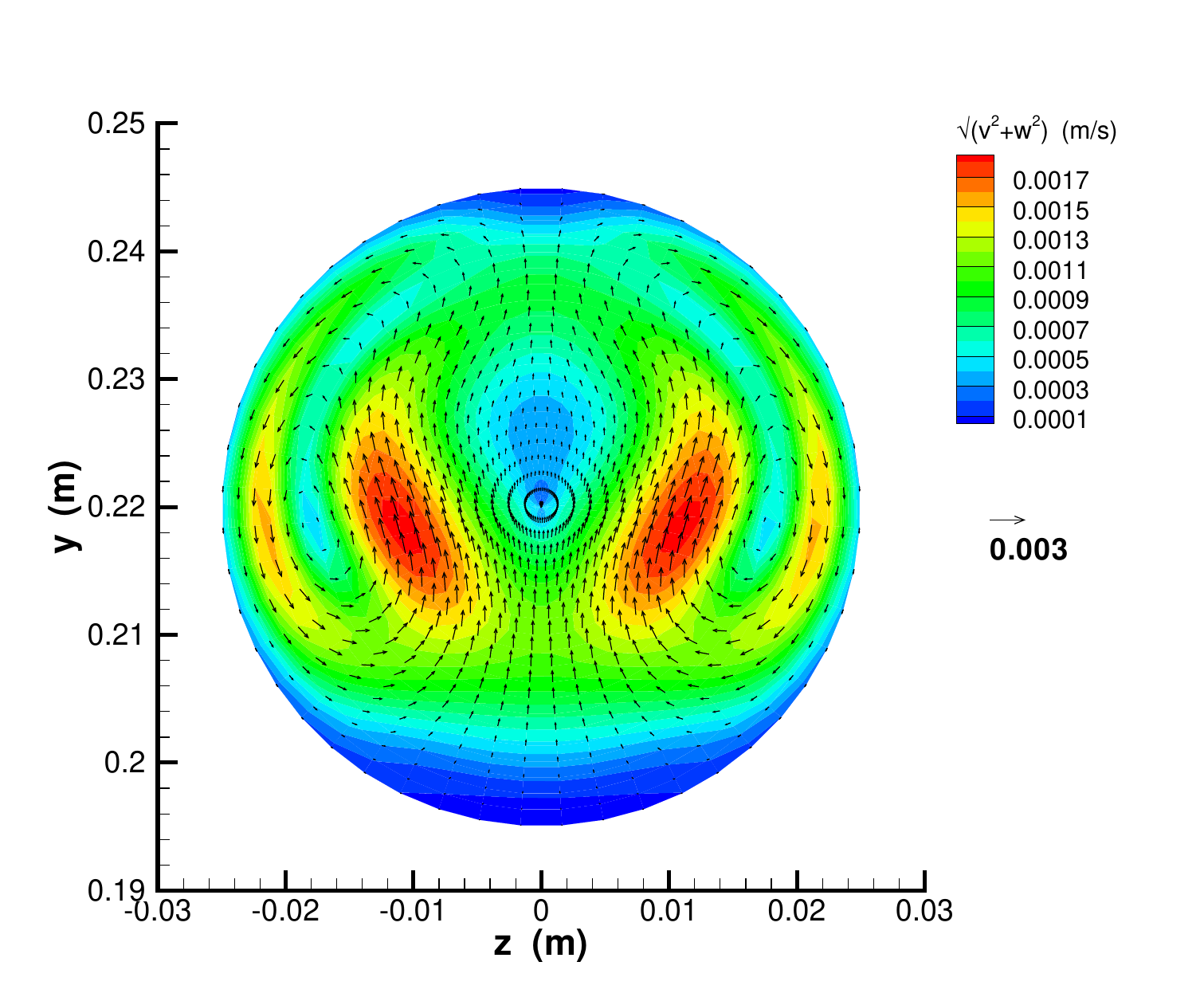}
						\caption{$\omega t = 270^\circ$}
					\end{subfigure}\\
					\begin{subfigure}[]{0.40\textwidth}
						\centering
						\includegraphics[width=\textwidth]{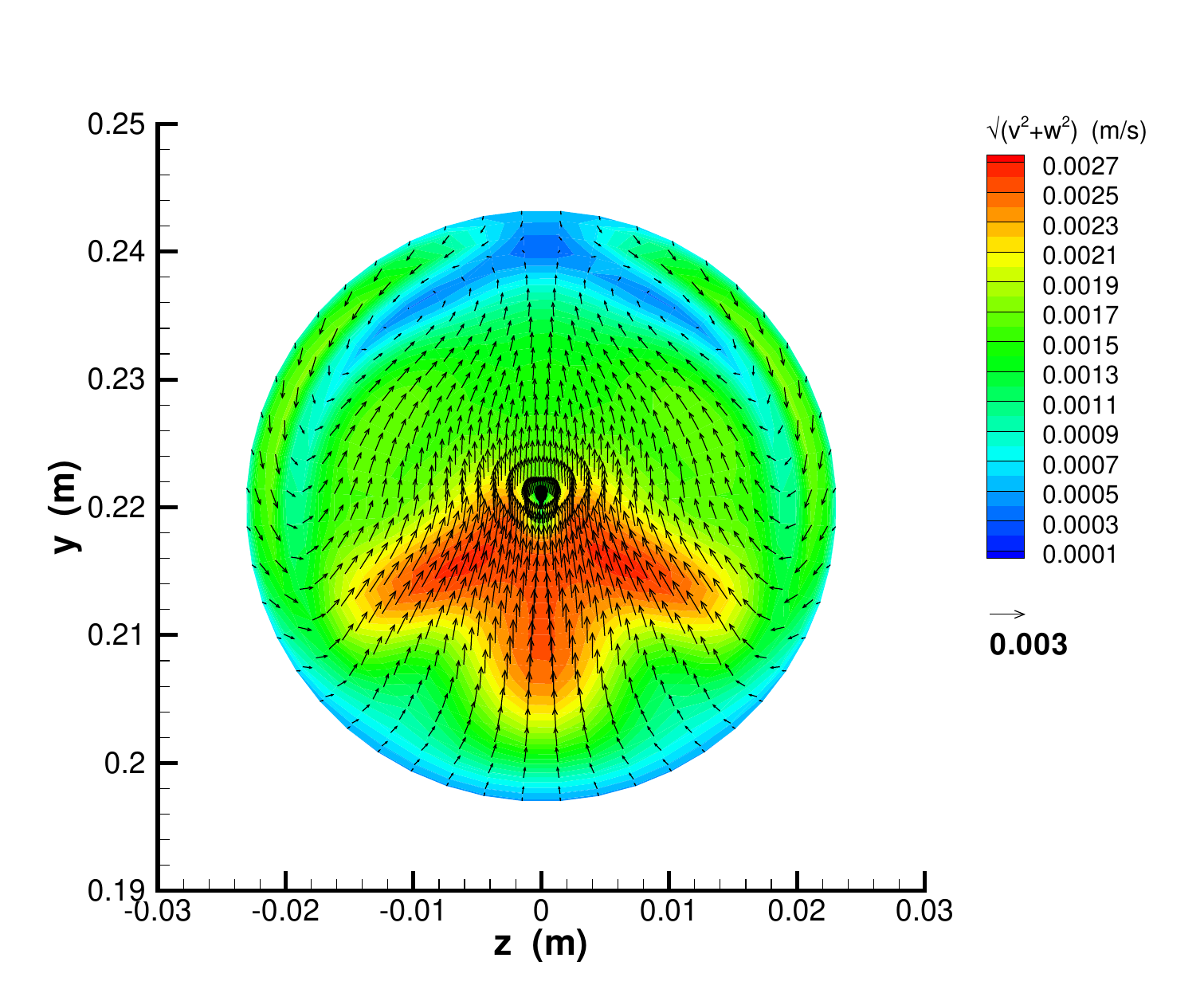}
						\caption{$\omega t = 360^\circ$}\label{fig:Compliance360B1}
					\end{subfigure} 
					\begin{subfigure}[]{0.40\textwidth}
						\centering
						\includegraphics[width=\textwidth]{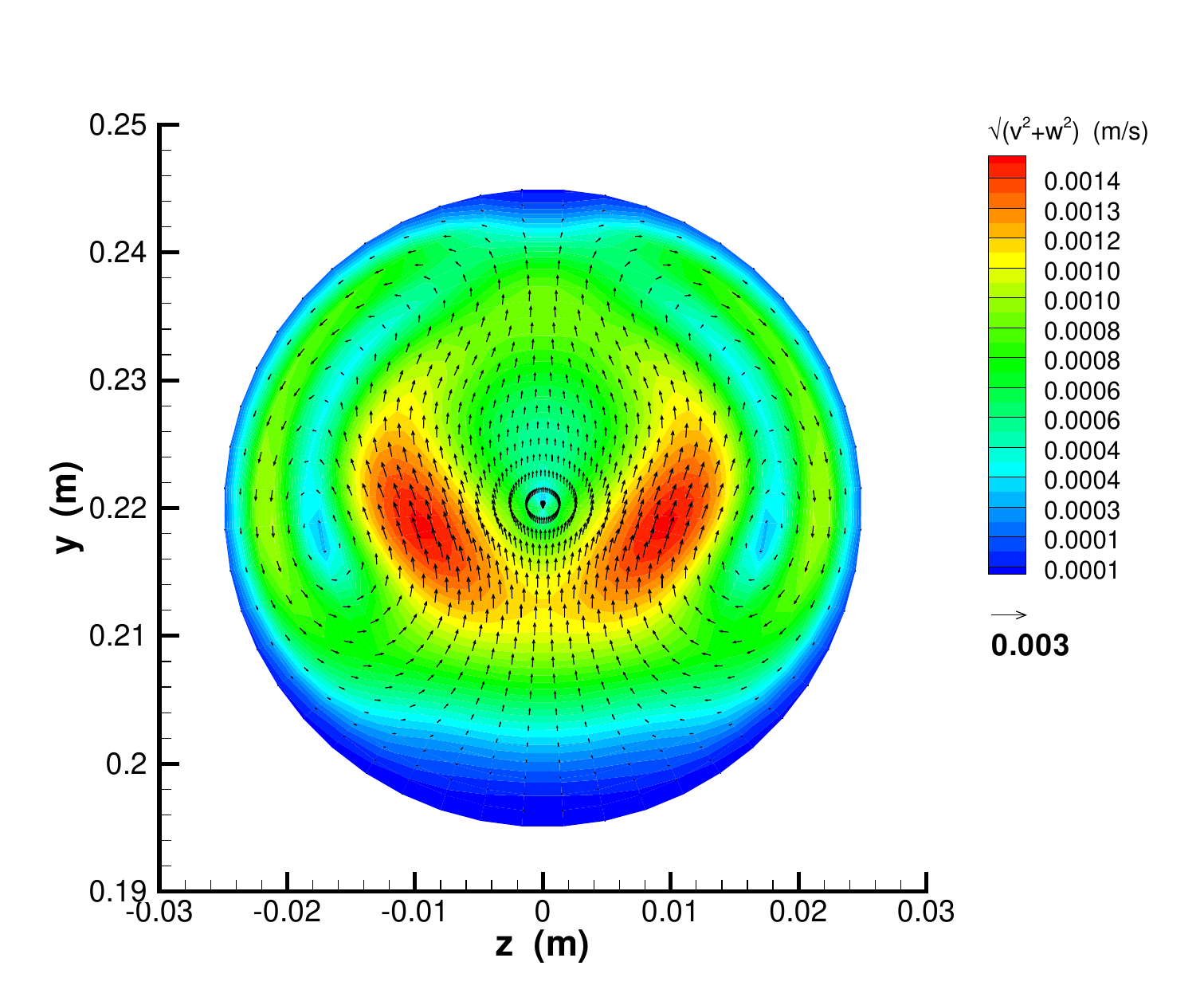}
						\caption{$\omega t = 360^\circ$}
					\end{subfigure}
\caption{Velocity vector field, tangential to the 
plane, interpolated along the $90^\circ$ cross section 
at different times and different rigidity coefficient, for $R_e = 600$, $\alpha_W = 17.17$, $\beta=1$ (left) and $\beta=10^{12}$ (right).}\label{fig:Compliance} 
\end{figure}

\begin{figure}[p]
\centering 
					\begin{subfigure}[]{0.40\textwidth}
						\centering
						\includegraphics[width=\textwidth]{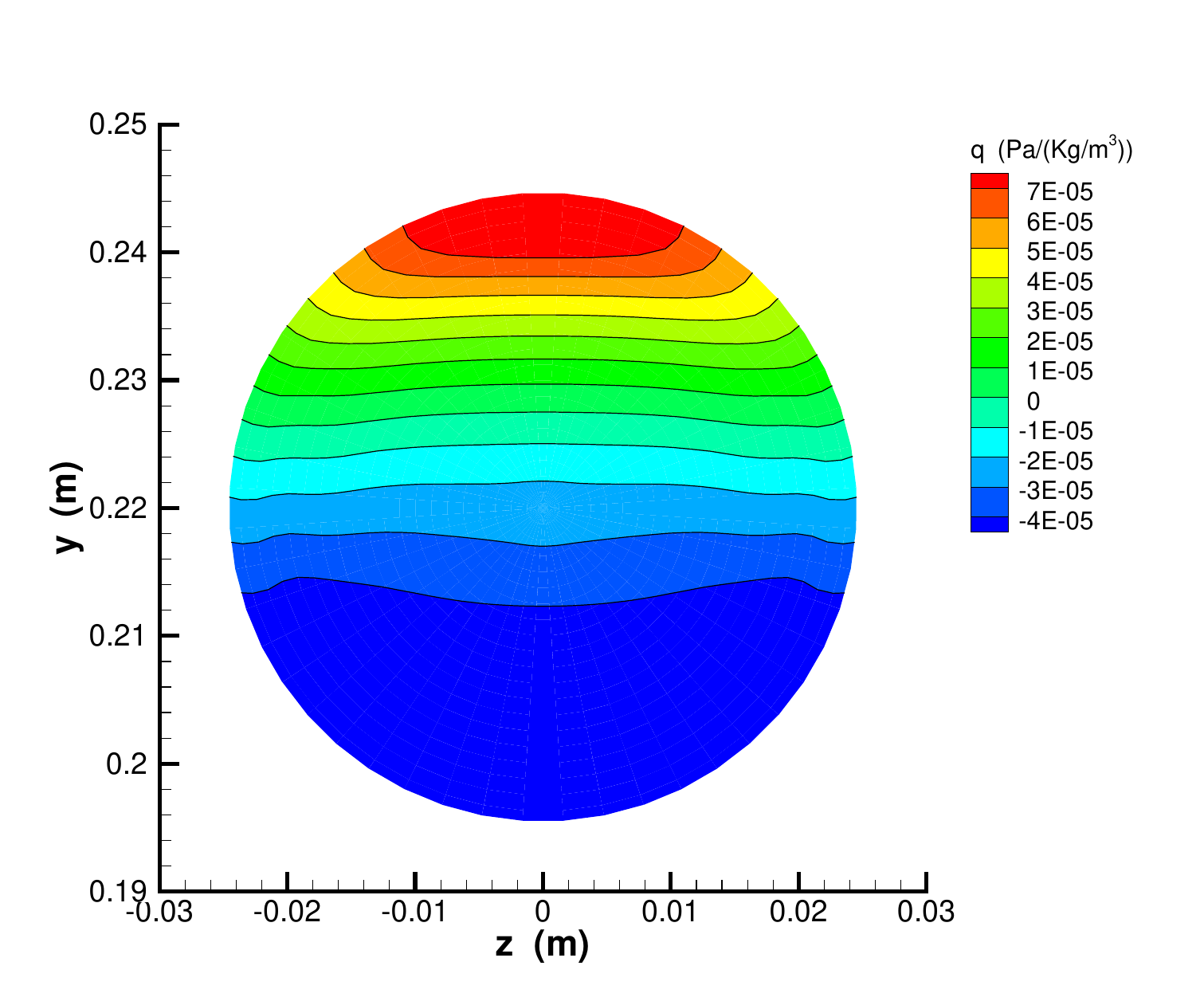}
						\caption{$\omega t = 90^\circ$}\label{fig:Compliance90B1Pnh} 
					\end{subfigure} 
					\begin{subfigure}[]{0.40\textwidth}
						\centering
						\includegraphics[width=\textwidth]{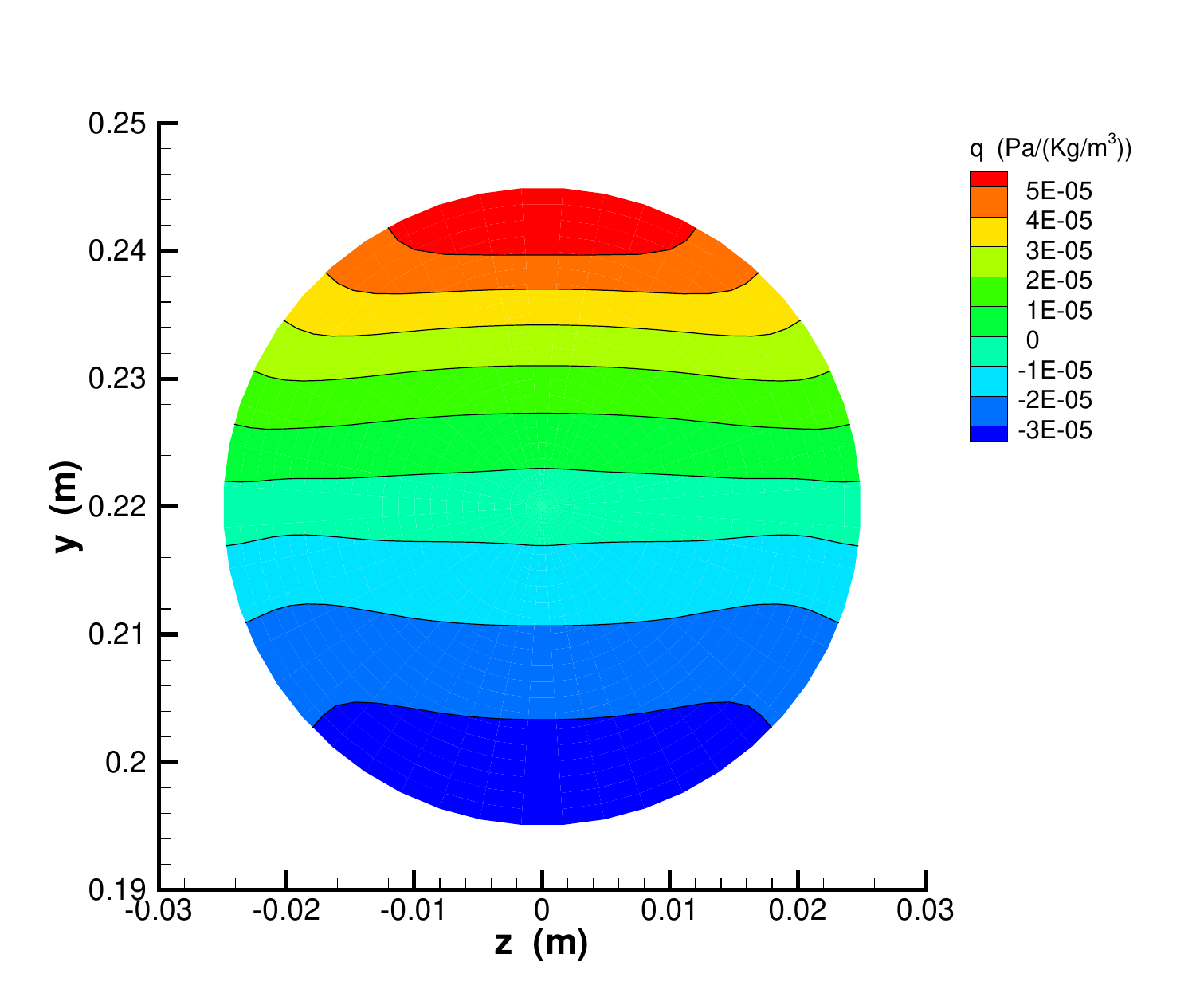}
						\caption{$\omega t = 90^\circ$}
					\end{subfigure}\\
					\begin{subfigure}[]{0.40\textwidth}
						\centering
						\includegraphics[width=\textwidth]{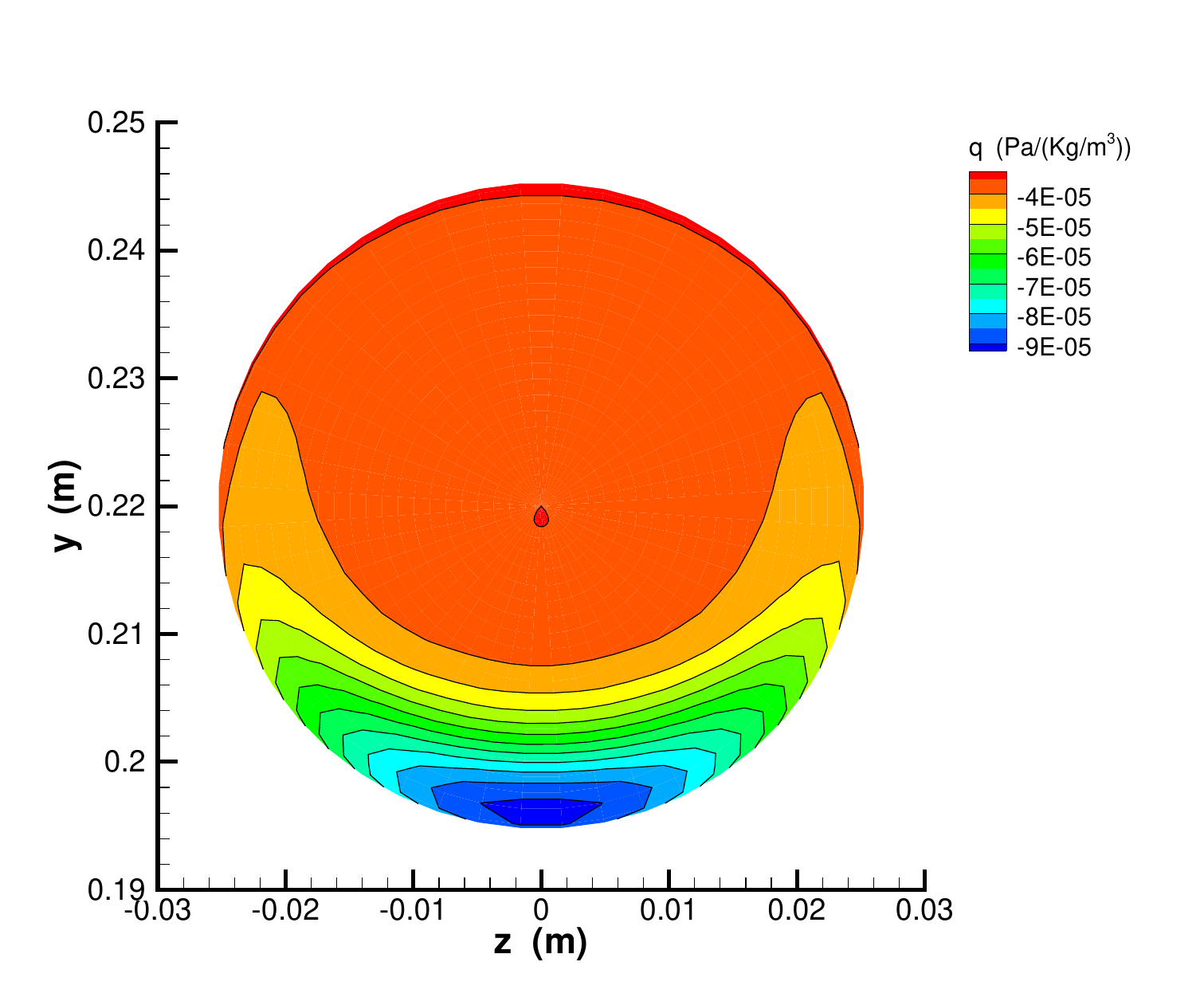}
						\caption{$\omega t = 180^\circ$}\label{fig:Compliance180B1Pnh} 
					\end{subfigure} 
					\begin{subfigure}[]{0.40\textwidth}
						\centering
						\includegraphics[width=\textwidth]{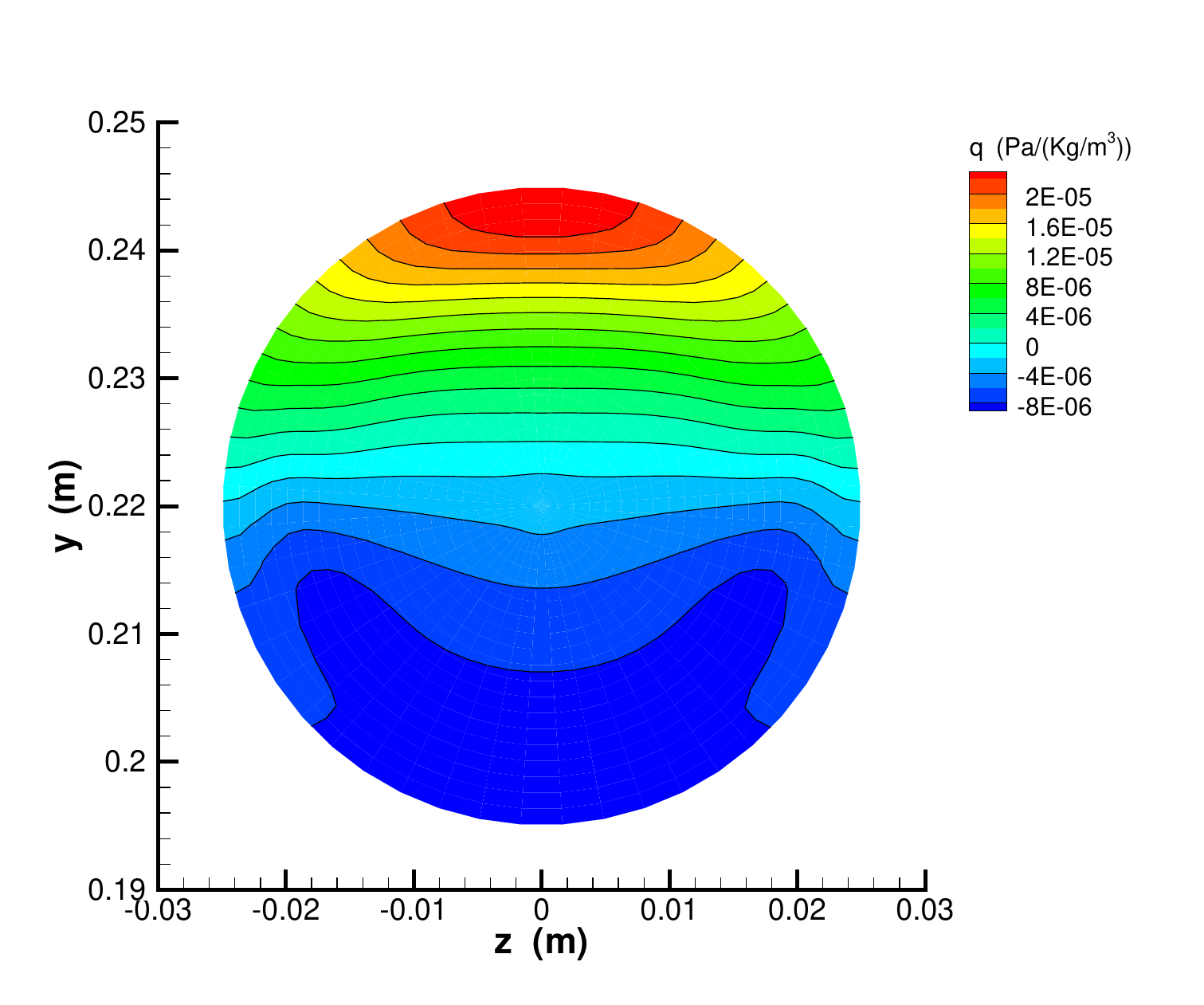}
						\caption{$\omega t = 180^\circ$}
					\end{subfigure}\\
					\begin{subfigure}[]{0.40\textwidth}
						\centering
						\includegraphics[width=\textwidth]{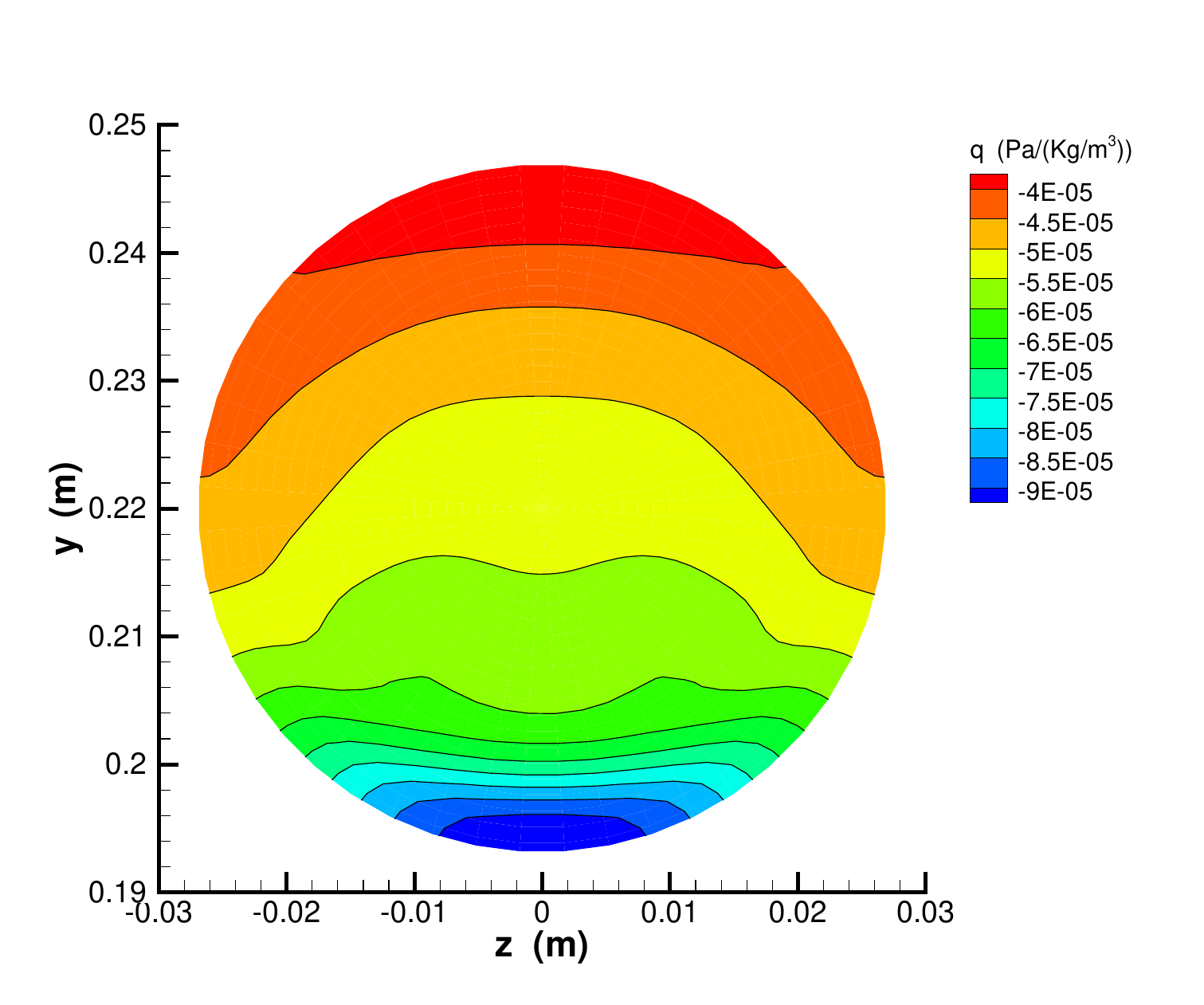}
						\caption{$\omega t = 270^\circ$}\label{fig:Compliance270B1Pnh} 
					\end{subfigure} 
					\begin{subfigure}[]{0.40\textwidth}
						\centering
						\includegraphics[width=\textwidth]{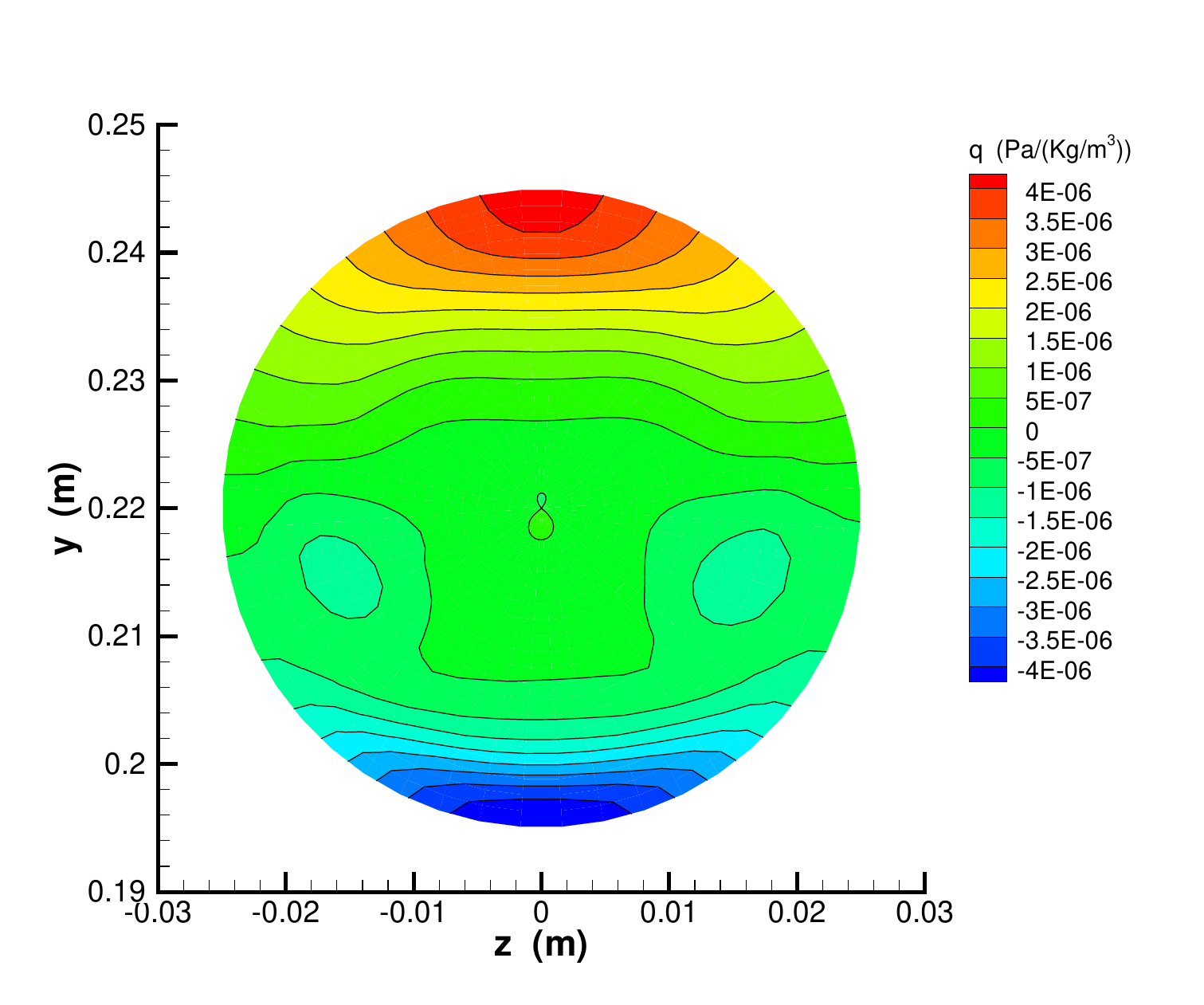}
						\caption{$\omega t = 270^\circ$}
					\end{subfigure}\\
					\begin{subfigure}[]{0.40\textwidth}
						\centering
						\includegraphics[width=\textwidth]{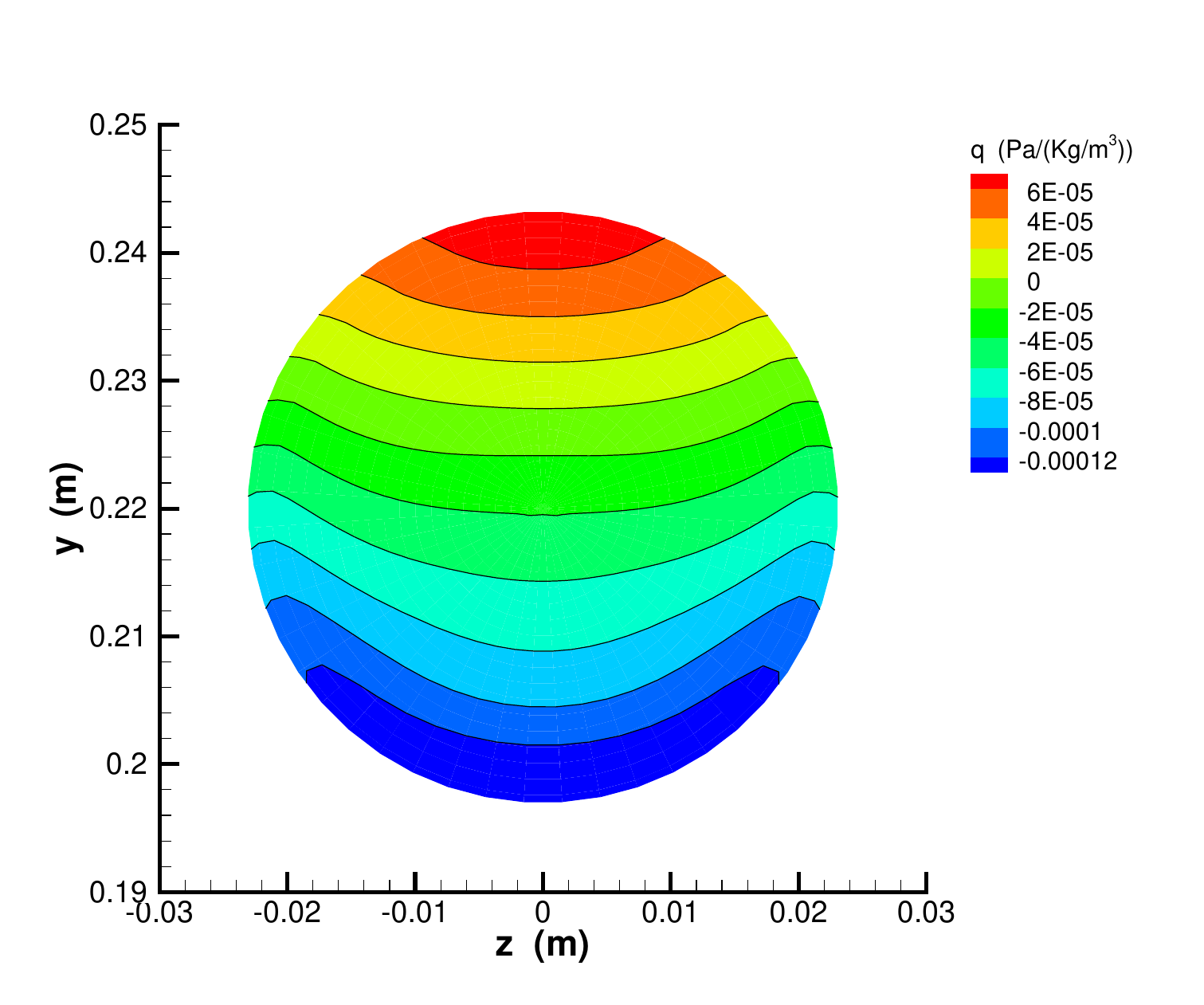}
						\caption{$\omega t = 360^\circ$}\label{fig:Compliance360B1Pnh}
					\end{subfigure} 
					\begin{subfigure}[]{0.40\textwidth}
						\centering
						\includegraphics[width=\textwidth]{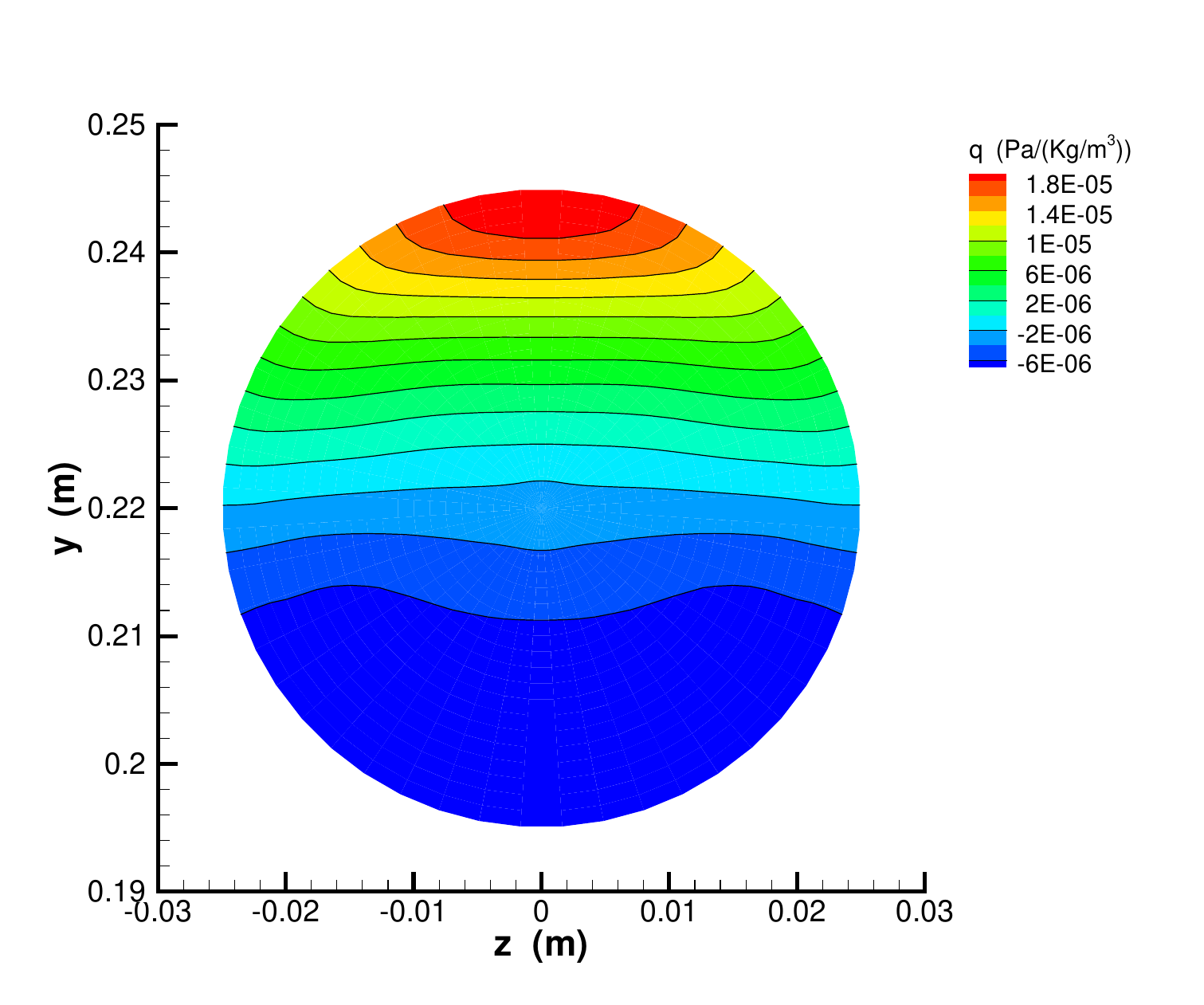}
						\caption{$\omega t = 360^\circ$}
					\end{subfigure}
\caption{Non hydrostatic pressure components interpolated along the $90^\circ$ cross section at different times and different rigidity coefficient, for $R_e = 600$, $\alpha_W = 17.17$, $\beta=1$ (left) and $\beta=10^{12}$ (right).}\label{fig:CompliancePnh} 
\end{figure}

Figure \ref{fig:Compliance} shows the velocity field tangential to the plane interpolated along the exit cross sections at different times, throughout one oscillation, for two different values of the rigidity coefficients, $\beta=1$ (compliant walls) and $\beta=10^{12}$ (rigid walls). The patterns of the secondary flow and the non-hydrostatic pressure components are significantly different from the rigid case (see Figures \ref{fig:Compliance}   and \ref{fig:CompliancePnh}). 
Principally, compliance has a double effect on the pulsatile flow: the axially reversed flow is amplified when the pressure gradient reaches its maximum (see Figure \ref{fig:ComplianceProfile}); secondly, it is shown that new  stagnation zones (with reference to the cross sectional motion) arise and that classical Dean vortices are consequently broken (see Figure \ref{fig:Compliance}). This is because the fluid is incompressible and aims to fill the  whole cross section, following the motion of the compliant vessel walls when they are stretching or tightening. Physically, three principal forces act at the same time and generate new flow patterns: the centrifugal force, the  centripetal pressure gradient and the radial acceleration due to the moving boundaries. In particular, when the velocity of the tube walls is large ($\omega t= 90^{\circ}$) the centrifugal flux is hindered in the inner part of the curve, and boosted in the outer; on the contrary, at the negative peak of the wall velocity ($\omega t= 270^{\circ}$) the centrifugal flux is boosted in the outer side, and hindered in the inner.
Note further that when the tube walls are stretching, the mean velocity increases, and when the tube is tightening, it becomes smaller. 
Unfortunately, we did not find available experimental reference data for this test problem to compare with.

\subsection{Wallclock time}%

In order to give a quantitative measure of the efficiency of the model and the demanded computational effort, the wallclock times needed with respect to a test problem of reference are tabulated. The Womersley problem presented above in Section \ref{sec.wom} is chosen as reference with respect to different discretization number of the radial and angular mesh, in order to compare the present three-dimensional non-hydrostatic model with the two-dimensional hydrostatic model of \cite{Casulli:2012}. Simulations are carried out with $t_e =2$, $\Delta t =0.01$,  by using $\theta = 0.5$, $\theta' = 1.0$, $N_x = 100$, $N_\varphi = 30$, and a sequence of successively refined meshes with $N_z=20$, $40$, $80$, $100$, $160$, $320$, $640$, $1000$ in order to see the computational cost linked to the radial discretization number (see Table \ref{tab:cap4wallclock}). Then, simulations are carried out by using $N_z = 50$ and a sequence of  successively refined meshes with $N_\varphi= 3$, $20$, $40$, $80$, $100$, $160$, $320$, $640$, $1000$ in order to see the computational cost linked to the angular discretization number  (see Table \ref{tab:cap4wallclock2}). 

The calculations have been performed with an Intel \emph{i}$7$ CPU having $1.70$ GHz clock frequency and $4$ GB of RAM. Problem sizes and wallclock times are listed in Tables \ref{tab:cap4wallclock} and \ref{tab:cap4wallclock2}, with distinction between the three-dimensional fully \emph{non-hydrostatic} approach of above and the three-dimensional \emph{hydrostatic} one that is obtained by completely skipping the second fractional step. These results confirm that the proposed numerical method is very accurate, highly efficient, and a good candidate for large-scale simulations of the human cardio-vascular system. With respect to the two-dimensional model, the three-dimensional model is obviously more time-consuming, but it is actually however a good alternative to other methods.  
The main causes of computing time with respect to the two-dimensional model of \cite{Casulli:2012} are actually (i) the increasing number of unknowns (the three velocity components and the non-hydrostatic pressure components in the \emph{three-dimensional space}) that are defined along (ii) an increasing number of control volumes and (iii) the resolution of the second-fractional step. Note, in particular, that this last step involves the fully three-dimensional coupling between velocities and pressures that comes from the discretization of the three-dimensional Navier-Stokes equations, the kinematic boundary condition and the conservation of volumes.

The \emph{hydrostatic splitting} obtained trough the fractional-step method allows a faster resolution of the fully non-hydrostatic problem. In order to give a quantitative measure of the \emph{saved} computing time,  in Table \ref{tab:cap4wallclock3} are reported also the wall-clock times needed for the resolution of the Womersley problem by a direct method obtained by setting the hydrostatic pressure components equal to zero ($p=q$), $\theta=0.5$ and  $\theta' = 0.5$. The results confirm that the semi-implicit fractional step approach proposed in this article is significantly cheaper form a computational point of view compared to the direct solution of the three-dimensional  incompressible Navier-Stokes equations. This is essentially due to the fact that the new numerical method proposed here can benefit from a dominant hydrostatic part of the pressure, while a direct approach cannot. 

\begin{table}
\begin{center}
\begin{tabular}{c c c c | c} 
\hline
Radial  			 & $N^\circ$ of control 		& Wallclock time  		 & Wallclock time  		& Wallclock time  		\\
resolution   	 & volumes 									& 3D non-hydrostatic   &  3D hydrostatic	  &  2D model	        \\
							 &  ($N_x N_z N_\varphi$)   &  (s) 							 &  (s)						 		&  (s)						 		\\
\hline
$N_z =$ $20\phantom{\,000}$    &    $\phantom{0\,0}60\,000$    & $8.081 E+00$   & $5.226 E+00$& $3.120 E-02$\\
$N_z =$ $40\phantom{\,000}$    &    $\phantom{0\,}120\,000$    & $1.587 E+01$   & $1.042 E+01$& $6.240 E-02$\\
$N_z =$ $80\phantom{\,000}$ 	 &    $\phantom{0\,}240\,000$    & $3.276 E+01$   & $2.051 E+01$& $1.248 E-01$ \\
$N_z =$ $100\phantom{\,00}$    &    $\phantom{0\,}300\,000$    & $3.933 E+01$   & $2.555 E+01$& $1.560 E-01$\\
$N_z =$ $160\phantom{\,00}$    &    $\phantom{0\,}480\,000$    & $6.301 E+01$   & $4.053 E+01$& $2.496 E-01$\\
$N_z =$ $320\phantom{\,00}$    &    $\phantom{0\,}960\,000$    & $1.310 E+02$   & $8.212 E+01$& $4.992 E-01$ \\
\hline
\end{tabular}
\end{center}
\caption{Problem size and wallclock time needed for the simulation of the Womersley problem: comparison between the three-dimensional non-hydrostatic model, the hydrostatic model and the two-dimensional approach of \cite{Casulli:2012}. $N_\varphi=30$, $N_x=100$.}
\label{tab:cap4wallclock}
\end{table}

\begin{table}
\begin{center}
\begin{tabular}{c c c  c | c} 
\hline
Angular  			 & $N^\circ$ of control 		& Wallclock time  		  & Wallclock time  	& Wallclock time  	\\
resolution   	 & volumes 									& 3D non-hydrostatic    &  3D hydrostatic	  &  2D model	        \\
							 &  ($N_x N_z N_\varphi$)   &  (s) 							    &  (s)						 	&  (s)						 	\\
\hline
$N_\varphi =$ $3\phantom{\,0000}$    &    $\phantom{0\,0}15\,000$    & $2.059 E+00$   & $1.420 E+00$   & $7.800 E-02$\\
$N_\varphi =$ $20\phantom{\,000}$    &    $\phantom{0\,}100\,000$    & $1.346 E+01$   & $1.034 E+01$   & \\
$N_\varphi =$ $40\phantom{\,000}$    &    $\phantom{0\,}200\,000$    & $2.742 E+01$   & $2.120 E+01$   & \\
$N_\varphi =$ $80\phantom{\,000}$ 	 &    $\phantom{0\,}400\,000$    & $6.363 E+01$   & $4.471 E+01$   & \\ 
$N_\varphi =$ $100\phantom{\,00}$    &    $\phantom{0\,}500\,000$    & $7.906 E+01$   & $5.546 E+01$   & \\
$N_\varphi =$ $160\phantom{\,00}$    &    $\phantom{0\,}800\,000$    & $1.380 E+02$   & $9.922 E+01$   & \\
$N_\varphi =$ $320\phantom{\,00}$    &    $\phantom{}1\,600\,000$    & $3.068 E+02$   & $2.166 E+02$   & \\
\hline
\end{tabular}
\end{center}
\caption{Problem size and wallclock time needed for the simulation of the Womersley  problem: comparison between the three-dimensional non-hydrostatic model, the hydrostatic model and the two-dimensional approach of \cite{Casulli:2012}. $N_z=50, N_x=100$.}
\label{tab:cap4wallclock2}
\end{table}

\begin{table}
\begin{center}
\begin{tabular}{c c c c } 
\hline
Radial  			 & $N^\circ$ of control 		& Wallclock time  		       & Wallclock time  		\\
resolution   	 & volumes 									& (fractional-step method) &  (direct approach)	\\
							 &  ($N_x N_z N_\varphi$)   &  (s) 							       &  (s)		\\
\hline
$N_z =$ $20\phantom{\,000}$    &    $\phantom{0\,0}60\,000$    & $8.081 E+00$   & $5.498 E+02$\\
$N_z =$ $40\phantom{\,000}$    &    $\phantom{0\,}120\,000$    & $1.587 E+01$   & $1.180 E+03$\\
$N_z =$ $80\phantom{\,000}$ 	 &    $\phantom{0\,}240\,000$    & $3.276 E+01$   & $9.020 E+03$\\
\hline
\end{tabular}
\end{center}
\caption{Problem size and wallclock time needed for the simulation of the Womersley problem: comparison between the fractional-step method of this article and a direct method (without \emph{hydrostatic splitting}). $N_\varphi=30$, $N_x=100$.}
\label{tab:cap4wallclock3}
\end{table}

\section{Conclusions}

A very efficient second order accurate semi-implicit finite difference and finite volume method for solving the governing equations of three-dimensional non-hydrostatic flows in compliant arterial vessels has been outlined. The computation is divided in two successive fractional steps. 
Some terms are selected for an implicit discretization and others for an explicit one. The grid is chosen in such a way that the resulting algebraic systems admit a fast and stable resolution. 
The chosen fractional-step structure, consisting in the \emph{splitting} of the hydrostatic from the non-hydrostatic pressure, improved significantly the computational efficiency. The more hydrostatic is the problem, the faster is  the algorithm. Further to this, the use of an Eulerian-Lagrangian method for treating advective terms allowed us to circumvent the usual CFL stability restriction that comes along with explicit schemes used to discretize 
the non-linear convection terms. 
The resulting semi-implicit algorithm is relatively simple, robust, efficient, and applies to non-hydrostatic three-dimensional, as well as hydrostatic three-, two- or one-dimensional  flow problems as a particular case. In particular, the model is locally and globally mass conservative. These features are illustrated on nontrivial test cases for flows in vessels with circular or elliptical cross section, where the exact analytical solution is known.  The moving boundaries, pressure and velocity fields are shown to be computed very accurately. Test cases of steady and pulsatile flow in uniformly curved pipes have also been presented. Axial velocity development and secondary flows have been shown and were compared  with previously published results in literature. 
This model could potentially be used in the future within realistic \emph{embedded multi-scale models} of the  human cardio-vascular system, 
see \cite{Formaggia:2001,Sherwin:2003,Quarteroni:2004,Formaggia:2007,Formaggia:2009,Quarteroni:2009}.

Future work will consider the modelling of three-dimensional junctions, an extension to the venous system and the implementation of higher order semi-implicit methods (see \cite{Dumbser:2013,TavelliDumbser}) that may 
become necessary for the accurate resolution of small-scale flow features that arise in the simulation of turbulent flows. 

\bibliography{Bibliography}{}
\bibliographystyle{plain}

\end{document}